\documentclass{article}

\usepackage{microtype}
\usepackage{graphicx}
\usepackage{subfigure}
\usepackage{booktabs} 

\usepackage{hyperref}

\usepackage[accepted]{icml2024}

\usepackage{amsmath}
\usepackage{amssymb}
\usepackage{mathtools}
\usepackage{amsthm}
\usepackage{makecell}
\usepackage{bbding}
\usepackage{enumerate}

\usepackage{graphicx}
\usepackage{times}  
\usepackage{helvet}  
\usepackage{courier}  
\usepackage{xcolor}
\usepackage{booktabs}

\usepackage{multirow}
\usepackage{bm}
\usepackage{colortbl}
\usepackage{amsmath}
\usepackage{wrapfig} 

\usepackage[capitalize,noabbrev]{cleveref}

\theoremstyle{plain}
\newtheorem{theorem}{Theorem}[section]

\newtheorem{lemma}[theorem]{Lemma}

\theoremstyle{definition}

\newtheorem{assumption}[theorem]{Assumption}
\theoremstyle{remark}
\newtheorem{remark}[theorem]{Remark}

\usepackage[textsize=tiny]{todonotes}

\icmltitlerunning{	Moreau Envelope for Nonconvex Bi-Level Optimization:  
	A Single-loop and Hessian-free Solution Strategy}

\begin{document}

\twocolumn[
\icmltitle{Moreau Envelope for Nonconvex Bi-Level Optimization: 
	\\A Single-loop and Hessian-free Solution Strategy}



\icmlsetsymbol{equal}{*}

\begin{icmlauthorlist}
	\icmlauthor{Risheng Liu}{DUT,Pazhou}
	\icmlauthor{Zhu Liu}{DUT}
	\icmlauthor{Wei Yao}{sustech,sustech2}
	\icmlauthor{Shangzhi Zeng}{uvic,sustech}
	\icmlauthor{Jin Zhang}{sustech2,sustech}
\end{icmlauthorlist}

\icmlaffiliation{DUT}{School of Software Technology, Dalian University of Technology, Dalian, China}

\icmlaffiliation{Pazhou}{Pazhou Lab, Guangzhou, China}

\icmlaffiliation{sustech}{National Center for Applied Mathematics Shenzhen, Southern University of Science and Technology, Shenzhen, China}

\icmlaffiliation{uvic}{Department of Mathematics and Statistics, University of Victoria, Victoria, British Columbia, Canada}

\icmlaffiliation{sustech2}{Department of Mathematics, Southern University of Science and Technology, Shenzhen, China}

\icmlcorrespondingauthor{Jin Zhang}{zhangj9@sustech.edu.cn}
\icmlkeywords{Machine Learning, ICML}

\vskip 0.3in
]



\printAffiliationsAndNotice{}  

\begin{abstract}
This work focuses on addressing two major challenges in the context of large-scale nonconvex Bi-Level Optimization (BLO) problems, which are increasingly applied in machine learning due to their ability to model nested structures. These challenges involve ensuring computational efficiency and providing theoretical guarantees. While recent advances in scalable BLO algorithms have primarily relied on lower-level convexity simplification, our work specifically tackles large-scale BLO problems involving nonconvexity in both the upper and lower levels. We simultaneously address computational and theoretical challenges by introducing an innovative single-loop gradient-based algorithm, utilizing the Moreau envelope-based reformulation, and providing non-asymptotic convergence analysis for general nonconvex BLO problems. Notably, our algorithm relies solely on first-order gradient information, enhancing its practicality and efficiency, especially for large-scale BLO learning tasks. We validate our approach's effectiveness through experiments on various synthetic problems, two typical hyper-parameter learning tasks, and a real-world neural architecture search application, collectively demonstrating its superior performance.
\end{abstract}

\section{Introduction}

Bi-Level Optimization (BLO) addresses the challenges posed by nested optimization structures that 
arise in a wide range of machine learning applications, such as
hyper-parameter optimization \cite{pedregosa2016hyperparameter,franceschi2018bilevel,mackay2018self}, 
meta learning
\cite{zugner2018adversarial,rajeswaran2019meta,KaiyiJi2020ConvergenceOM}, neural architecture search~\cite{liu2018darts,chen2019progressive,elsken2020meta}, etc. 
Refer to recent survey papers \cite{liu2021investigating,zhang2023introduction} for more applications of BLO in machine learning, computer vision and signal processing.
\begin{table*}[t]
	\caption{
		Comparison of our method MEHA with closely related works for addressing {\bf nonconvex-nonconvex BLO} 
		(
		IAPTT-GM
		\cite{liu2021towards},
		BOME!    
		\cite{NeurIPS2022-Liu}, 
		V-PBGD
		\cite{ShenC23},
		GALET
		\cite{xiao2023generalized},
		SLM
		\cite{lu2023slm}
		). 
		Different methods employ distinct stationary measures, so we do not delve into complexity comparison here.
		Below, PL Condition represents the Polyak-Łojasiewicz (PL) condition;
		$L$-Smooth means the Lipschitz continuous gradient condition; Bounded and Gradient-Bounded specify that $| F(x,y)|\leq C$ 
		and   
		$\|\nabla_{y}F(x,y)\|\leq C$ for all $(x,y)$, respectively;
		$F$ and $f$ are UL and LL objectives, respectively. 
	}
	\label{BLO-table1}
	\begin{center}
		\begin{tabular}{cccccc} 
			\multicolumn{1}{c}{\bf \makecell{Method}}  &\multicolumn{1}{c}{\bf \makecell{Upper-Level\\ Objective}}
			&\multicolumn{1}{c}{\bf \makecell{Lower-Level\\ Objective}}
			&\multicolumn{1}{c}{\bf 
				\makecell{Hessian-free
				}
			}
			&\multicolumn{1}{c}{\bf \makecell{Single-loop}
			}
			&\multicolumn{1}{c}{\bf \makecell{
					Non-asymptotic}}
			\\ \hline 
			\makecell{
				IAPTT-GM
			}         
			&Smooth 
			&\makecell{
				$L$-Smooth\\
				\& Compactness
			}   
			&\XSolidBrush
			&\XSolidBrush
			&\XSolidBrush
			\\ \hline 
			\makecell{
				GALET 
			}       
			&\makecell{$L$-Smooth \\  
				\& Gradient-Bounded
			}
			&\makecell{PL Condition\\
				\& $L$-Smooth }        
			&\XSolidBrush
			&\XSolidBrush
			&  \CheckmarkBold
			\\ \hline 
			\makecell{
				BOME!    
			}             
			&\makecell{$L$-Smooth   \\ 
				\& Bounded \\
				\& Gradient-Bounded
			}
			&\makecell{PL Condition\\
				\& $L$-Smooth }        
			& \CheckmarkBold
			&\XSolidBrush
			& \CheckmarkBold
			\\ \hline 
			\makecell{
				V-PBGD 
			}       
			&\makecell{$L$-Smooth  \\  
				\& Gradient-Bounded
			}
			&\makecell{PL Condition\\
				\& $L$-Smooth }        
			&\CheckmarkBold
			&\XSolidBrush
			&  \CheckmarkBold
			\\ \hline 
			\makecell{
				SLM 
			}       
			&\makecell{$L$-Smooth  \\  
				\& Compactness
			}
			&\makecell{PL Condition\\
				\& $L$-Smooth }        
			&\CheckmarkBold
			&\XSolidBrush
			&  \CheckmarkBold
			\\ \hline 
			\makecell{	MEHA  \\
				(smooth case)         
			}
			&$L$-Smooth  
			&\makecell{$L$-Smooth}        
			& \CheckmarkBold
			& \CheckmarkBold
			& \CheckmarkBold
			\\ \hline 
			\makecell{	MEHA        
				\\
				(general case)   
			}
			&$L$-Smooth  
			&\makecell{$L$-Smooth Part\\ 
				+ Weakly Convex
				\\ \ \ \ \ Nonsmooth Part}        
			&\CheckmarkBold
			& \CheckmarkBold
			& \CheckmarkBold
			\\ \hline 
		\end{tabular}
	\end{center}
\end{table*}
The inherent nested nature gives rise to several difficulties and hurdles in effectively solving BLO problems. 
Over the past decade, a large number of BLO methods have emerged, 
with a primary emphasis on addressing BLO problems
featuring strongly convex lower-level (LL) objective. 
The LL strong convexity assumption ensures the uniqueness of LL minimizer (a.k.a., LL Singleton), which simplifies both the optimization process and theoretical analysis, see, e.g.,  \cite{franceschi2018bilevel,ghadimi2018approximation,grazzi2020iteration,ji2020bilevel,chen2021closing,ji2022will,hong2020two}. 
To mitigate the restrictive LL Singleton condition,
another line of research is dedicated to BLO with convex LL problems, which bring about several challenges such as the presence of multiple LL optimal solutions (a.k.a., Non-Singleton).
This may hinder the adoption of implicit-based approaches that rely on the implicit function theorem.
To address this concern, recent advances include:  
aggregation methods 
\cite{liu2020generic,li2020improved,liu2022general,pmlr-v202-liu23y};
difference-of-convex algorithm \cite{gao2022value,ye2023difference}; 
primal-dual algorithms 
\cite{sow2022constrained}; first-order penalty methods \cite{lu2023first}.

In this work, we study a BLO problem with a nonconvex LL problem:\begin{equation}\label{general_problem}
	\min_{x \in X, y \in Y}  \, F(x,y)  
	\quad 
	\text{s.t.} \quad   y \in S(x),
\end{equation}
where $S(x)$ denotes the set of optimal solutions for the LL problem given by
\begin{equation}\label{LL_problem}
	\min_{y \in Y}~  \varphi(x,y) := f(x,y) + g(x,y),
\end{equation}
where $X$ and $Y$ are closed convex sets in $ \mathbb{R}^n$ and $ \mathbb{R}^m$, respectively. The function $f(x,y): \mathbb{R}^n\times\mathbb{R}^m\rightarrow \mathbb{R}$ is smooth, and generally nonconvex, while
$g(x,y): \mathbb{R}^n\times\mathbb{R}^m\rightarrow \mathbb{R}$ is potentially nonsmooth  with respect to (w.r.t.) the LL variable $y$. For specific conditions governing $F$, $f$ and $g$, we refer the reader to Assumptions \ref{assump-UL}, \ref{assump-LL}. 

While the nonconvex-convex BLO has been extensively studied in the literature, efficient methods for nonconvex-nonconvex BLO remain under-explored.
Beyond LL convexity, 
\cite{liu2021towards} proposed an iterative differentiation-based BLO method; 
\cite{arbel2022non} extended the approximate implicit differentiation approach; 
\cite{liu2021value} 
firstly 
utilized the value function reformulation of BLO to develop algorithms in machine learning.
All of these works, however, tend to be complicated and impractical for large-scale BLO problems, 
and lack a non-asymptotic analysis.
When LL objective satisfies the Polyak-Łojasiewicz (PL) or local PL conditions,
\cite{NeurIPS2022-Liu} introduced a fully first-order value function-based BLO algorithm with non-asymptotic convergence analysis. 
Recently, while still considering the PL condition,
\cite{huang2023momentumbased} introduced a momentum-based BLO algorithm;  \cite{xiao2023generalized} proposed a generalized alternating method; \cite{lu2023slm} proposed a smoothed first-order Lagrangian method;
\cite{ShenC23} proposed a penalty-based fully first-order BLO algorithm. 
However, the existing methods still present two significant challenges: ensuring computational efficiency and offering theoretical guarantees in the absence of the PL condition.
A concise comparison with works closely related to ours is summarized in Table \ref{BLO-table1}.

\subsection{Main Contributions}

To the best of our knowledge, this work is the first study to utilize Moreau envelope-based reformulation of BLO, 
originally presented in \cite{gao2023moreau}, to design a single-loop and Hessian-free gradient-based algorithm with non-asymptotic convergence analysis for general BLO problems with potentially nonconvex and nonsmooth LL objective functions. 
This setting encompasses a wide range of machine learning applications, see, e.g., the recent surveys \cite{liu2021investigating,zhang2023introduction}. 
Conducting non-asymptotic analysis for our algorithm, which addresses nonconvex LL problem, poses substantial challenges. 
Existing single-loop gradient-based methods generally require the LL objective to either be strongly convex or satisfy the PL condition, as a mechanism to control the approximation errors incurred when utilizing a single gradient descent step to approximate the real LL optimal solution. 
Our approach mitigates this limitation by employing Moreau envelope-based reformulation, where the proximal LL problem may exhibit strong convexity even if the original LL problem is nonconvex. 
Consequently, this enables effective error control and facilitates the algorithm's non-asymptotic convergence analysis for nonconvex LL problem. 
We summarize our contributions as follows.
\begin{itemize}
	\item We propose the {\bf M}oreau {\bf E}nvelope based {\bf H}essian-free {\bf A}lgorithm (MEHA), 
	for general BLO problems with nonconvex and probably nonsmooth LL objective functions. 
	MEHA avoids second-order derivative approximations related to the Hessian matrix and can be implemented efficiently in a single-loop manner, enhancing its practicality and efficiency for large-scale BLOs. 
	
	\item We provide a rigorous analysis of the non-asymptotic convergence of MEHA under milder conditions, avoiding the need for either the convexity assumption or the PL condition on LL problem. In the context of the smooth BLO scenario, our assumption simplifies to UL and LL objective functions being $L$-smooth.
	
	\item We validate the effectiveness and efficiency of MEHA on various synthetic problems, two typicial hyper-parameter learning tasks and the real-world neural architecture search application. These experiments collectively substantiate its superior performance.
\end{itemize}

\subsection{Related Work}
\label{relatedwork}

We provide a brief review of recent works closely related to ours, with a detailed review available in Section \ref{relatedwork2}.

{\bf Nonconvex-Nonconvex BLO.} 
Beyond LL convexity, \cite{huang2023momentumbased} introduced a momentum-based implicit gradient approach with convergence under PL conditions and nondegenerate LL Hessian. \cite{xiao2023generalized} proposed a new stationary metric for nonconvex-PL BLOs, alongside a generalized alternating method. These methods, however, necessitate complex Hessian-related computations. Conversely, \cite{NeurIPS2022-Liu} offered a Hessian-free algorithm leveraging value function reformulation. \cite{lu2023slm}  proposed a smoothed first-order Lagrangian method focused on the relaxation of the value function reformulation, also avoiding Hessian computations. 
\cite{ShenC23} introduced a penalty-based algorithm applicable to both unconstrained and constrained BLOs. These Hessian-free BLO methods all have non-asymptotic convergence under PL conditions, but featuring double-loop structures. 

{\bf Moreau Envelope Based Reformulation for BLO.} 
The authors in the work \cite{gao2023moreau} pioneered the utilization of Moreau envelope of LL problem for investigating BLOs. They devised an inexact proximal Difference-of-Weakly-Convex algorithm in a double-loop manner for solving constrained BLO problems with convex LL problem.
Subsequently, through the application of Moreau envelope-based reformulation of BLO, the study \cite{bai2023optimality} derived new  optimality condition results. More recently, while focusing on the convex smooth LL problem scenario, \cite{yao2024constrained}  introduced a min-max version of Moreau envelope for convex constrained LL problem and developed a single-loop Hessian-free gradient-based algorithm for the constrained BLO problem.

\section{
	A Single-loop and Hessian-free Solution Strategy
}
\label{alg}

\subsection{Moreau Envelope Based Reformulation}

In this work, we focus on developing a single-loop algorithm for solving BLO with a nonconvex LL problem. Our approach builds upon the  Moreau envelope based reformulation of BLO, which is initially proposed for convex LL scenarios in \cite{gao2023moreau}. The reformulation is expressed as follows:
\begin{equation}\label{wVP}
	\begin{aligned}
		\min_{(x,y)\in X \times Y}  & \,F(x,y) 
		\quad
		\mathrm{s.t.} & \varphi(x,y)-v_\gamma (x,y)\leq 0,
	\end{aligned}
\end{equation} 
where $ v_{\gamma}(x,y) $ represents the Moreau envelope associated with the LL problem, defined as:
\begin{equation}\label{def_vg}
	v_{\gamma}(x,y)  
	:=\inf_{\theta \in Y} \left \{ \varphi(x, \theta) + \frac{1}{2\gamma}\|\theta-y\|^2 \right \},
\end{equation}
with $\gamma>0$.
It is important to note that $\varphi(x,y) \ge v_\gamma (x,y)$ holds for all $(x,y)\in X \times Y$. For the convex LL scenarios where $\varphi(x,y)$ is convex in $y \in Y$ for any $x \in X$, the equivalence between the reformulated and the original BLO problems is established in Theorem 1 of \cite{gao2023moreau}. 
In scenarios where $\varphi$ lacks convexity, yet $\varphi(x, \cdot)$ is $\rho_{\varphi_2}$-weakly convex 
\footnote{A function $h:\mathbb{R}^p\rightarrow\mathbb{R}\cup\{\infty\}$ is $\rho$-weakly convex if $h(z)+\frac{\rho}{2}\|z\|^2$ is convex. In the context where $z=(x,y)$, we always say $h$ is $(\rho_1,\rho_2)$-weakly convex if $h(x,y)+\frac{\rho_1}{2}\|x\|^2+\frac{\rho_2}{2}\|y\|^2$ is convex.}
on $Y$ and $\gamma \in (0, 1/\rho_{\varphi_2})$,
we establish in Theorem \ref{thm_reformulate-a} that the reformulation (\ref{wVP}) is equivalent to a relaxed version of BLO problem (\ref{general_problem}),
\begin{equation}\label{problem_relax}
	\min_{x \in X, y \in Y}  \, F(x,y)  
	\quad
	\mathrm{s.t.}  \quad  y \in 
	\tilde{S}(x),
\end{equation}
where $\tilde{S}(x):= \left\{ 
y ~|~ 0 \in \nabla_{y}f(x,y) + \partial_y g(x,y) 
+ \mathcal{N}_Y(y)
\right\}$,
$\partial_y g(x,y)$ denotes the partial Fr\'{e}chet (regular) subdifferential of $g$ w.r.t. LL variable at $(x,y)$, and $\mathcal{N}_Y(y)$ signifies the normal cone to $Y$ at $y$. 
The stationary condition characterizing $\tilde{S}(x)$ is the first-order optimality conditions of LL problem within the context of this work, specifically, Assumption \ref{assump-LL}. This can be validated through the application of subdifferential sum rules, see, e.g., Proposition 1.30, Theorem 2.19 in \cite{mordukhovich2018variational}.

Specifically, the reformulated problem (\ref{wVP}) becomes equivalent to the original BLO problem when the set $\tilde{S}(x)$ coincides with $S(x)$. This equivalence holds, for instances, when $\varphi(x, \cdot)$ is convex or $\varphi (x,y)\equiv f(x,y)$ and it satisfies the PL condition, that is, there exists $\mu>0$ such that for any $x \in X$, the inequality 
$
\| \nabla_{y} f(x,y) \|^2 \geq 2\mu
\big(
f(x,y)-\inf_{\theta\in\mathbb{R}^m} f(x, \theta)
\big)
$
holds for all $y \in \mathbb{R}^m$.

We next introduce an approximation problem to the Moreau envelope based reformulation (\ref{wVP}), constructed by penalizing the constraint $ \varphi(x,y)-v_\gamma (x,y)\leq 0$ in (\ref{wVP}). This approximation problem is formulated as:
\begin{equation}\label{problem_pen}
	\min_{(x,y)\in X \times Y} \psi_{c_k}(x,y) := F(x,y) + c_k \big( \varphi(x,y) - v_{\gamma} (x,y)\big),
\end{equation}	
where $\varphi(x,y) = f(x,y) + g(x,y)$, and $c_k$ acts as the penalty parameter.
Importantly, as $c_k \rightarrow \infty$, any limit point of the sequence of solutions to the approximation problem (\ref{problem_pen}), associated with varying values of $c_k$, is a solution to the Moreau envelope based reformulation  (\ref{wVP}). This convergence is formally established in Theorem \ref{Thm_app}.

\subsection{Gradient of Moreau Envelope $v_{\gamma}(x,y)$}

Before presenting our proposed method, we briefly review some relevant preliminary results related to $v_{\gamma}(x,y)$, with a special focus on its gradient.
Assuming that $\varphi(x, y)$ is $(\rho_{\varphi_1},\rho_{\varphi_2})$-weakly convex on $X\times Y$, we demonstrate that for $\gamma \in (0, \frac{1}{2\rho_{\varphi_2}})$,
the function $v_{\gamma}(x,y) + \frac{\rho_{v_1}}{2}\|x\|^2 + \frac{\rho_{v_2}}{2}\|y\|^2$
is convex on $X \times \mathbb{R}^m$ when $\rho_{v_1} \ge \rho_{\varphi_1}$ and $\rho_{v_2}  \ge 1/\gamma $. That is, $v_{\gamma}(x,y)$ is weakly convex, as detailed in Lemma~\ref{Lem1-a}. 
Lastly, we define  
\begin{equation}\label{solutionset1}
	{S}_{\gamma}(x, y) := \mathrm{argmin}_{\theta \in Y} \left \{ \varphi(x, \theta) + \frac{1}{2\gamma}\|\theta-y\|^2 \right \}.
\end{equation}
For $\gamma \in (0, \frac{1}{2\rho_{\varphi_2}})$, the solution set ${S}_{\gamma}( {x}, {y})=\{\theta_{\gamma}^*(x,y)\}$ is a singleton. Furthermore, when the gradients $\nabla_{x} f(x,y)$ and $\nabla_{x} g(x,y)$ exist, $v_{\gamma}(x,y)$ is differentiable and the gradient $\nabla v_\gamma$ of $v_{\gamma}(x,y)$ at $(x,y)$ can be expressed as follows, 
\begin{equation}\label{valuefinc}
	\nabla v_\gamma =
	\begin{pmatrix} 
		\nabla_x f( {x},   \theta_{\gamma}^*(x,y) ) + \nabla_{x} g( {x},   \theta_{\gamma}^*(x,y) )
		\\
		\left( {y} -  \theta_{\gamma}^*(x,y)\right)/ \gamma 
	\end{pmatrix},
\end{equation}
which is established in Lemma~\ref{lem2-a}.

\subsection{Single-loop {\bf M}oreau {\bf E}nvelope Based {\bf H}essian-free {\bf A}lgorithm (MEHA)}
In this part, we present a single-loop algorithm for the general BLO problem (\ref{general_problem}), via solving its approximation problem (\ref{problem_pen}). 
At each iteration, we intend to employ the alternating proximal gradient method to the approximation problem (\ref{problem_pen}) for updating the variables $(x,y)$, starting from the current iterate $(x^k,y^k)$. However, this process encounters certain challenges. Specifically, as outlined in (\ref{valuefinc}), to calculate the gradient of the Moreau envelope $v_{\gamma}(x,y)$ at a given iterate $(x^k, y^k)$, we have to know $\theta^*_\gamma(x^k, y^k)$, which is the unique solution to the proximal LL problem (\ref{def_vg}) with $(x,y) = (x^k, y^k)$. Exact computation of $\theta^*_\gamma(x^k, y^k)$ is computationally intensive. To mitigate this, we introduce a new iterative sequence $\{\theta^{k}\}$, where each $\theta^{k+1}$ approximates $\theta^*_\gamma(x^k, y^k)$. At each iteration, $\theta^{k+1}$ is generated by applying a single proximal gradient step to the proximal LL problem (\ref{def_vg}) with $(x,y) = (x^k, y^k)$, starting from the current $\theta^k$. This update is formalized as:
\begin{equation}\label{update_z}
	\theta^{k+1} = \mathrm{Prox}_{\eta_k\tilde{g}(x^k,\cdot)}\big( \theta^k - \eta_k \big(\nabla_y f(x^k,\theta^k)  + \frac{\theta^k - y^k}{\gamma}\big) \big),
\end{equation}
where $\eta_k$ is the stepsize, and $\tilde{g}(x,y) := g(x,y) + \delta_{Y}(y)$ represents the nonsmooth part of LL problem. Here, $\mathrm{Prox}_{h}(y)$ denotes the proximal mapping of a function $h:\mathbb{R}^m\rightarrow \mathbb{R}\cup \{\infty\}$, 
$	\mathrm{Prox}_{h}(y):=\arg\min_{\theta\in\mathbb{R}^m}\left\{h(\theta)+  \|\theta-y\|^2/2 \right\}$.
Leveraging the formula (\ref{valuefinc}), an approximation for $\nabla v_{\gamma}(x^k,y^k) $ is constructed using $\theta^{k+1}$ as an approximation for $\theta^*_{\gamma}(x^k,y^k) $,
$$( \nabla_x f(x^k, \theta^{k+1}) + \nabla_{x}g(x^k, \theta^{k+1}), (y^k - \theta^{k+1})/ \gamma).$$ 

Subsequently, for the update of variable $x$, we construct an approximation for $[\nabla_x \psi_{c_k}(x^k, y^k)]/c_k$, serving as the update direction $d_{x}^k$:
\begin{equation}\label{def_dx}
	\begin{aligned}
		d_{x}^k &:= \frac{1}{c_k} \nabla_x F(x^k,y^k) + \nabla_x f(x^k, y^k) + \nabla_{x}g(x^k, y^k) \\
		&\quad \ -  \nabla_x f(x^k, \theta^{k+1}) - \nabla_{x}g(x^k, \theta^{k+1}),  
	\end{aligned}
\end{equation}
and update $x^{k+1}$ using the direction $d_{x}^k$:
\begin{equation}\label{update_x}
	x^{k+1} = \mathrm{Proj}_X \left( x^k - \alpha_k d_x^k \right),
\end{equation}
where $\alpha_k >0$ is the stepsize, and $\mathrm{Proj}_X$ denotes the Euclidean projection operator.
Next,  for variable $y$, we construct an approximation for $[\nabla_y (\psi_{c_k} - g)(x^{k+1}, y^k)]/c_k$ as the update direction $d_{y}^k$:
\begin{equation}\label{def_dy}
	d_{y}^k := \frac{1}{c_k} \nabla_y F(x^{k+1},y^k) + \nabla_y f(x^{k+1}, y^k) - \frac{y^k - \theta^{k+1}}{\gamma}.
\end{equation}
And $y^{k+1}$ is updated via a proximal gradient step along $d_{y}^k$:
\begin{equation}\label{update_y}
	y^{k+1} = \mathrm{Prox}_{\beta_k\tilde{g}(x^{k+1}, \cdot)}\left(  y^{k}  - \beta_k  d_y^k\right),
\end{equation}
where $\beta_k >0$ is the stepsize.

The complete algorithm is outlined in Algorithm \ref{MEHA}, 
whose specific form for smooth BLOs, i.e., $g(x,y)\equiv0$, is provided in Algorithm \ref{MEHA-SC} in Appendix.
\begin{algorithm}[h]
	\caption{
		Single-loop Moreau Envelope based Hessian-free Algorithm (MEHA)
	}\label{MEHA}
	\hspace*{0.02in} {\bf Input:} Initialize $x^0,y^0,\theta^0$, stepsizes $\alpha_k, \beta_k, \eta_k$, proximal parameter $\gamma$, penalty parameter $c_k$;
	\begin{algorithmic}[1]
		\FOR{$k=0,1,\dots,K-1$} 
		
		\STATE update $\theta^{k+1}$ by (\ref{update_z});
		
		\STATE calculate $d_x^k$ as in \eqref{def_dx} and update $x^{k+1}$ by (\ref{update_x});
		\STATE calculate $d_y^k$ as in \eqref{def_dy} and update $y^{k+1}$ by (\ref{update_y}).
	
		\ENDFOR
	\end{algorithmic}
\end{algorithm}

\section{Theoretical Investigations}
\label{theory}

\subsection{General Assumptions}

Throughout this work, we assume the following standing assumptions on $F$, $f$ and $g$ hold.

\begin{assumption}\label{assump-UL}
	The UL objective $F$ is bounded below on $X \times Y$, denoted by $\underline{F} := \inf_{(x, y)\in X\times Y} F(x,y) > -\infty$. Furthermore, $F$ is $L_F$-smooth\footnote{ A function $h$ is said to be $L$-smooth on $X \times Y$ if $h$ is continuously differentiable and its gradient $\nabla h$ is $L$-Lipschitz continuous on $X \times Y$.} 
	on $X \times Y$.  The smooth component of LL objective $f(x,y)$ is $L_f$-smooth on $X \times Y$. 
\end{assumption}

\begin{assumption}\label{assump-LL}
	The nonsmooth component of LL objective, $g(x,y)$, satisfies one of the following conditions:
	\begin{enumerate}[(i)]
		\item $g(x,y)=\hat{g}(y)$ with $\hat{g}(y)$ being weakly convex on $Y$;
		
		\item $g(x,y) = x\| y\|_1$ with $X=\mathbb{R}_+$ and $Y = \mathbb{R}^m$;
		
		\item $g(x,y) = \sum_{j = 1}^J x_j\| y^{(j)}\|_2$, where $\{1, \ldots,m\}$ is divided into $J$ groups, $y^{(j)}$ denotes the corresponding $j$-th group of $y$, and $X=\mathbb{R}^J_+$, $Y = \mathbb{R}^m$;
		
		\item The nonsmooth component $g(x,y)$ is $( \rho_{g_1}, \rho_{g_2} )$-weakly convex on $X \times Y$, i.e., 
		$
		g(x,y) + \frac{\rho_{g_1}}{2}\|x\|^2 + \frac{\rho_{g_2}}{2}\|y\|^2
		$
		is convex on $X \times Y$. Additionally, the gradient $\nabla_{x} g(x,y)$ exists and is $L_g$-Lipschitz continuous on $X \times Y$.
		Moreover, let $\tilde{g}(x,y) := g(x,y) + \delta_{Y}(y)$,
		there exist positive constants $L_{\tilde{g}},  \bar{s}$ such that for any $x, x' \in X$, $\theta \in Y$ and $s \in (0, \bar{s}]$,
		\begin{equation}\label{ass-LL-eq}
			\left\| 
			\mathrm{Prox}_{s \tilde{g}(x,\cdot)} (\theta) - \mathrm{Prox}_{s \tilde{g}(x',\cdot)}(\theta)
			\right\| 
			\leq L_{\tilde{g}} \|x - x' \|.
		\end{equation}
	\end{enumerate}
\end{assumption}

The demand for weak convexity in Assumption \ref{assump-LL} (i) is relatively lenient; a broad spectrum of functions meet this requirement. This encompasses conventional nonconvex regularizers like the Smoothly Clipped Absolute Deviation (SCAD) and the Minimax Concave Penalty (MCP) (refer to Section 2.1 of \cite{bohm2021variable}).  
It is worth noting that Assumptions \ref{assump-LL} (i), (ii), and (iii) represent specific instances of Assumption \ref{assump-LL}(iv). Specifically, in Section \ref{A9-assump} of Appendix, we provide comprehensive proofs demonstrating that Assumptions \ref{assump-LL} (ii) and (iii) are special cases of Assumption \ref{assump-LL} (iv).

These assumptions considerably alleviate LL problem's smoothness requirements prevalent in BLO literature. 
Even within the context of smooth BLO, 
our assumptions only require that UL and LL objectives are both $L$-smooth, without imposing any conditions on the boundedness of $\nabla_{y} F(x,y)$, as illustrated in Table \ref{BLO-table1}. Consequently, our problem setting encompasses a broad range of practical scenarios, see, e.g., the learning models in  \cite{grazzi2020iteration}.

Finally, leveraging the descent lemma (e.g., Lemma 5.7 of \cite{beck2017first}), it can be obtained that any function featuring a Lipschitz-continuous gradient is inherently weakly convex. Thus, under Assumption \ref{assump-UL}, $f(x,y)$ is $( \rho_{f_1}, \rho_{f_2} )$-weakly convex over $X \times Y$, with $\rho_{f_1}=\rho_{f_2}=L_{f}$. 
This leads us to the following result.
\begin{lemma}
	Under Assumptions \ref{assump-UL} and \ref{assump-LL}, LL objective $\varphi(x,y)$ is $( \rho_{\varphi_1},\rho_{\varphi_2} )$-weakly convex on $X \times Y$, where $\rho_{\varphi_1} = \rho_{f_1} + \rho_{g_1} $ and $\rho_{\varphi_2} = \rho_{f_2} + \rho_{g_2}$.
\end{lemma} 

\subsection{Non-asymptotic Convergence Results}
\label{convergence}

We consider the following residual function, denoted as $R_k(x,y)$, for the approximation problem (\ref{problem_pen}),
\begin{equation}\label{residual}
	\begin{aligned}
		R_k:= \mathrm{dist} \big( 0, \big[\nabla F 
		+ c_k \big(\nabla f + \partial g- \nabla v_{\gamma} \big) 
		+ \mathcal{N}_{X \times Y}\big]\big).   
	\end{aligned}
\end{equation}
This residual function is a stationarity measure for the approximation problem (\ref{problem_pen}). In particular,
it is evident that $R_k(x,y) = 0$ if and only if $0 \in \partial \psi_{c_k}(x,y) + \mathcal{N}_{X \times Y}(x,y) $, i.e., the point $(x,y)$ is a stationary point to the approximation problem (\ref{problem_pen}).

\begin{theorem}\label{prop1}
	Under Assumptions \ref{assump-UL} and \ref{assump-LL}, 
	suppose $\gamma \in (0, \frac{1}{2\rho_{f_2} + 2\rho_{g_2} })$, $c_k =\underline{c} (k+1)^p$ with $p \in [0,1/2)$, $\underline{c} >0$ and $\eta_k \in [ \underline{\eta}, (1/\gamma - \rho_{f_2})/(L_f+1/\gamma)^2] \cap  [ \underline{\eta}, 1/\rho_{g_2})$ with $ \underline{\eta} > 0$, 
	then there exists $c_{\alpha}, c_\beta > 0$ such that when $\alpha_k \in ( \underline{\alpha}, {c_\alpha})$ and $\beta_k \in ( \underline{\beta}, c_\beta)$ with $\underline{\alpha}, \underline{\beta} > 0$, the sequence of $(x^k, y^k, \theta^k)$ generated by MEHA (Algorithm \ref{MEHA}) satisfies
	\[
	\min_{0 \le k \le K} \left\| \theta^{k} - \theta_{\gamma}^*(x^{k},y^{k}) \right\| = O\left(\frac{1}{K^{1/2}} \right),
	\]
	and 
	\[
	\min_{0 \le k \le K} R_k(x^{k+1}, y^{k+1}) = O\left(\frac{1}{K^{(1-2p)/2}} \right).
	\]
	Furthermore, if $p \in (0,1/2)$ and the sequence $ \psi_{c_k}(x^k, y^k) $ is upper-bounded,
	the sequence of $(x^k, y^k)$ satisfies
	\[
	\varphi(x^K, y^K) - v_\gamma(x^K, y^K) = O\left(\frac{1}{K^{p}} \right).
	\]
\end{theorem}
\begin{remark}
	In the case where the penalty parameter $c_k$ remains fixed at a constant value $\underline{c}$, applying $p = 0$ in Theorem \ref{prop1} yields a non-asymptotic $O(1/\sqrt{K})$ convergence rate result of MEHA in solving the approximation problem (\ref{problem_pen}), with $c_k = \underline{c}$, based on its stationarity measure.
\end{remark}

The norm of the hyper-gradient, denoted as $\nabla \Phi(x)$, where $\Phi(x):=F(x, y^*(x))$ and $y^*(x)$ representes the solution of the LL problem, is widely employed as a measure of stationary for BLO problems in the existing literature, see, e.g., \cite{ghadimi2018approximation,ji2020bilevel,chen2021closing,ji2022will,hong2020two,kwon2023fully}. Nevertheless, it is crucial to note that this stationarity measure may not be suitable for BLO problems investigated in this work, as $\nabla \Phi(x)$ may not be well-defined. If we focus on a specific class of  BLO problems, where $\nabla \Phi(x)$ is well-defined, 
we can establish the connection between the residual function $R_k(x,y)$ and the hypergradient $\nabla\Phi(x)$, as demonstrated in Lemma \ref{hypergradient-approx}.  
Consequently, Theorem \ref{prop1} offers a non-asymptotic convergence result of MEHA based on the norm of $\nabla \Phi(x)$. 

\begin{theorem}\label{prop2}
	Under the same conditions as detailed in Theorem \ref{prop1}, and with the additional assumptions that $X=\mathbb{R}^n$, $Y=\mathbb{R}^m$, $g(x,y) = 0$, $f(x,y)$ is $\mu$-strongly convex with respect to $y$ for any $x$, $F$ and $f$ satisfy certain smoothness conditions (detailed in Lemma \ref{hypergradient-approx}), $p \in (0,1/2)$, and $\gamma > 1/\mu$,
	the sequence of $(x^k, y^k, \theta^k)$ generated by MEHA (Algorithm \ref{MEHA}) satisfies
	\[
	\min_{0 \le k \le K}\|\nabla\Phi(x^{k+1})\|=O(\frac{1}{K^{1/2-p}} + \frac{1}{K^{p}}).
	\]
\end{theorem}

\begin{table*}[!h]
	\centering
	\footnotesize
	\renewcommand{\arraystretch}{1.0}
	\caption{Comparison of  total iterative time  with representative  BLO methods on the LL  non-convex case with different dimensions.}~\label{tab:numerical}\vspace{-0.2cm}
	\setlength{\tabcolsep}{4.5mm}{
		\begin{tabular}{|c|cccc|cccc|}
			\hline
			Category     & \multicolumn{4}{c|}{Dimension =1 }                                 & \multicolumn{4}{c|}{Dimension =10 }                                                    \\ \hline
			Methods & \multicolumn{1}{c|}{BVFIM} & \multicolumn{1}{c|}{BOME} & \multicolumn{1}{c|}{IAPTT} & \multicolumn{1}{c|}{MEHA} & \multicolumn{1}{c|}{BVFIM} & \multicolumn{1}{c|}{BOME} & \multicolumn{1}{c|}{IAPTT} & MEHA  \\ \hline
			Time (S) & \multicolumn{1}{c|}{53.50}   & \multicolumn{1}{c|}{5.062} & \multicolumn{1}{c|}{21.23}   & \multicolumn{1}{c|}{\textbf{0.774}} & \multicolumn{1}{c|}{224.37}   & \multicolumn{1}{c|}{6.348} & \multicolumn{1}{c|}{35.73} &   \textbf{2.342}\\ \hline
			Category     & \multicolumn{4}{c|}{Dimension =50}                                 & \multicolumn{4}{c|}{Dimension =100}                                                    \\ \hline
			Methods & \multicolumn{1}{c|}{BVFIM} & \multicolumn{1}{c|}{BOME} & \multicolumn{1}{c|}{IAPTT} & \multicolumn{1}{c|}{MEHA} & \multicolumn{1}{c|}{BVFIM} & \multicolumn{1}{c|}{BOME} & \multicolumn{1}{c|}{IAPTT} & MEHA  \\ \hline
			Time (S) & \multicolumn{1}{c|}{1106.45} & \multicolumn{1}{c|}{11.02} & \multicolumn{1}{c|}{154.54}   & \multicolumn{1}{c|}{\textbf{9.552}} & \multicolumn{1}{c|}{1846.55}  & \multicolumn{1}{c|}{16.48} & \multicolumn{1}{c|}{290.82} &   \textbf{12.26}\\ \hline
		\end{tabular}
	}

\end{table*}
\section{Experimental Results}
In this section, we verify the convergence behavior and applicable practicality of MEHA on a series of different synthetic problems. Then we compare the performances on three real-world applications, including few-shot learning, data hyper-cleaning and neural architecture search. 

\subsection{Synthetic Numerical Verification}

\begin{figure}[h!]
	\footnotesize
	\centering
	\setlength{\tabcolsep}{1pt}

	\begin{tabular}{cc}
		\includegraphics[width=0.24\textwidth]{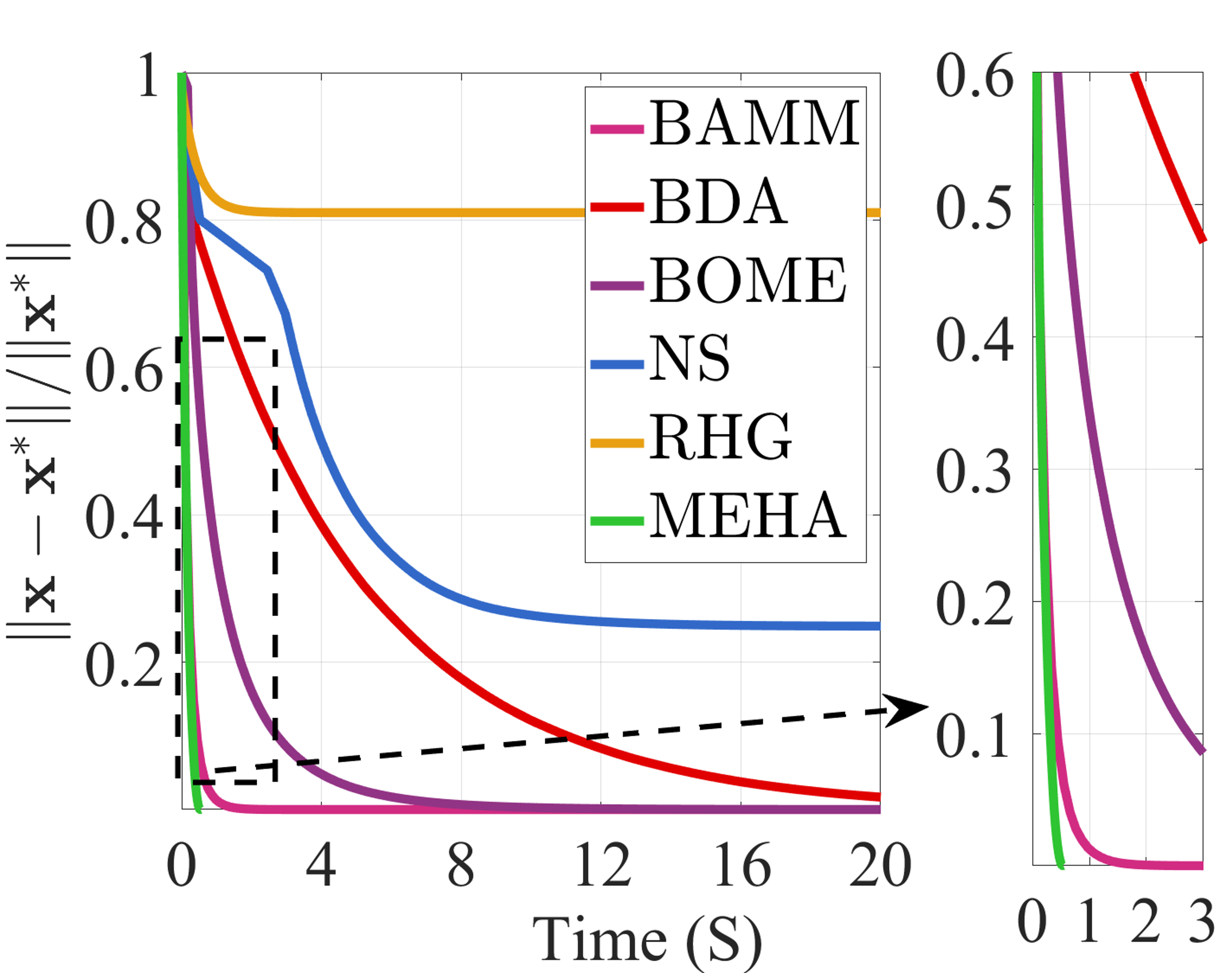}
		&  \includegraphics[width=0.24\textwidth]{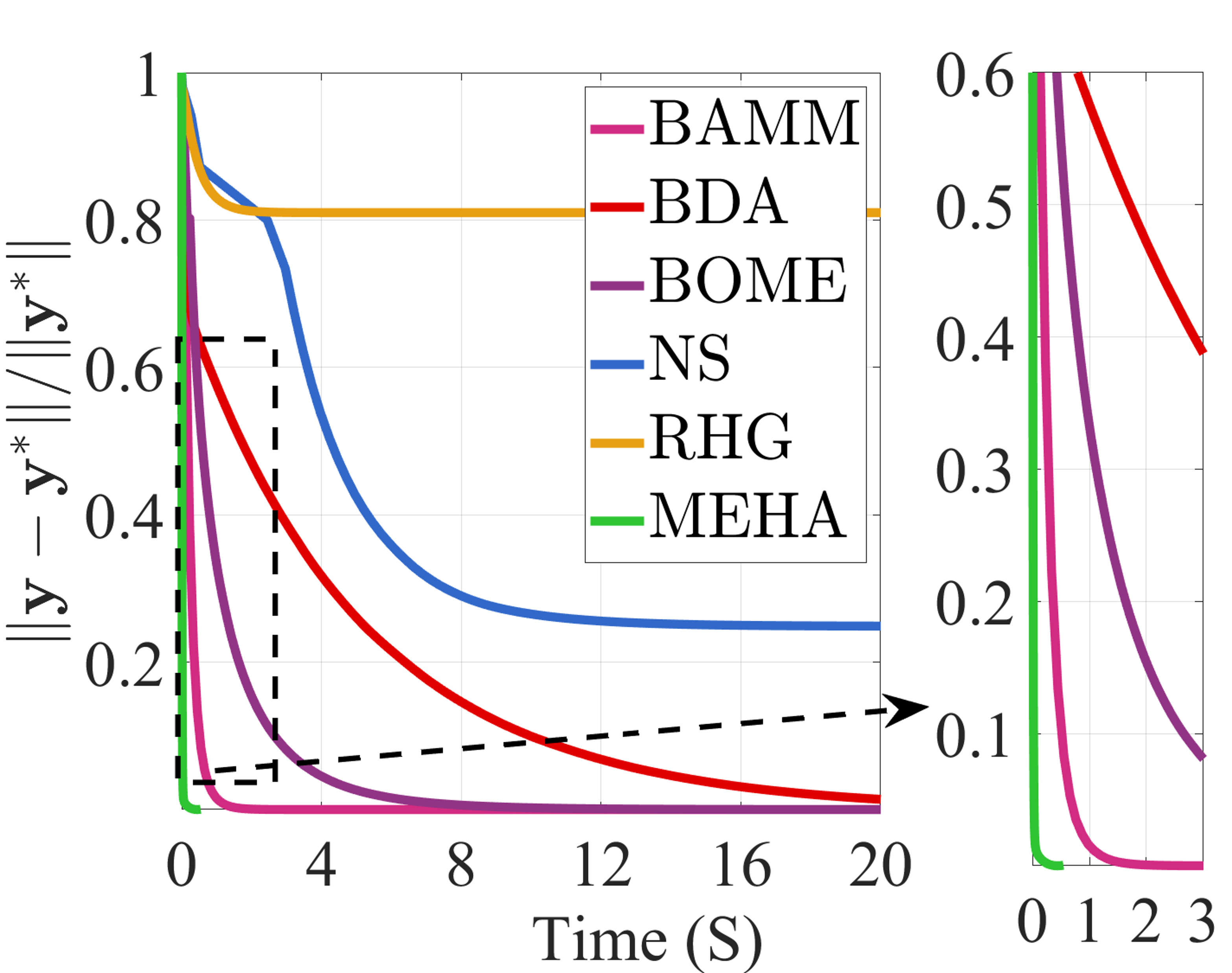}\\
	\end{tabular}
	\vspace{-0.2cm}
	\caption{ Illustrating the convergence curves of advanced BLO schemes  and MEHA by the criteria, $\|x-x^{*}\|/\|x^{*}\|$ and $\|y-y^{*}\|/\|y^{*}\|$, under  LL  merely convex case.}\label{fig:mc_trad}
\end{figure}
\textbf{LL Merely Convex Case.} 
We first demonstrate the high efficiency of the proposed method on a toy example in LL merely convex case, expressed as follows:
\begin{eqnarray}\label{eq:LLC}
	\begin{aligned}
		&\min_{x\in \mathbb{R}^n, y=\left(y_1, y_2\right)\in \mathbb{R}^{2n}} \frac{1}{2}\left\|x-y_2\right\|^2+\frac{1}{2}\left\|y_1-\mathbf{e}\right\|^2  \\
		&\mathrm { s.t. }\quad y=\left(y_1, y_2\right) \in \arg \min _{\left(y_1, y_2\right) \in \mathbb{R}^{2 n}} \frac{1}{2}\left\|y_1\right\|^2-x^{\top} y_1.\\
	\end{aligned}
\end{eqnarray} 
The unique solution of this BLO is $(\mathbf{e},\mathbf{e},\mathbf{e})$.
We plot the convergence curves on the Figure~\ref{fig:mc_trad}. 
Convergence criteria for numerical experiments are reported in Section~\ref{sec:detail}.
Our proposed MEHA achieves the fastest convergence compared with existing BLO scheme, while some other methods fail to converge to the correct solution.
\begin{figure}[h!]
	\footnotesize
	\centering
	\setlength{\tabcolsep}{1pt}

	\begin{tabular}{cc}
		\includegraphics[width=0.24\textwidth]{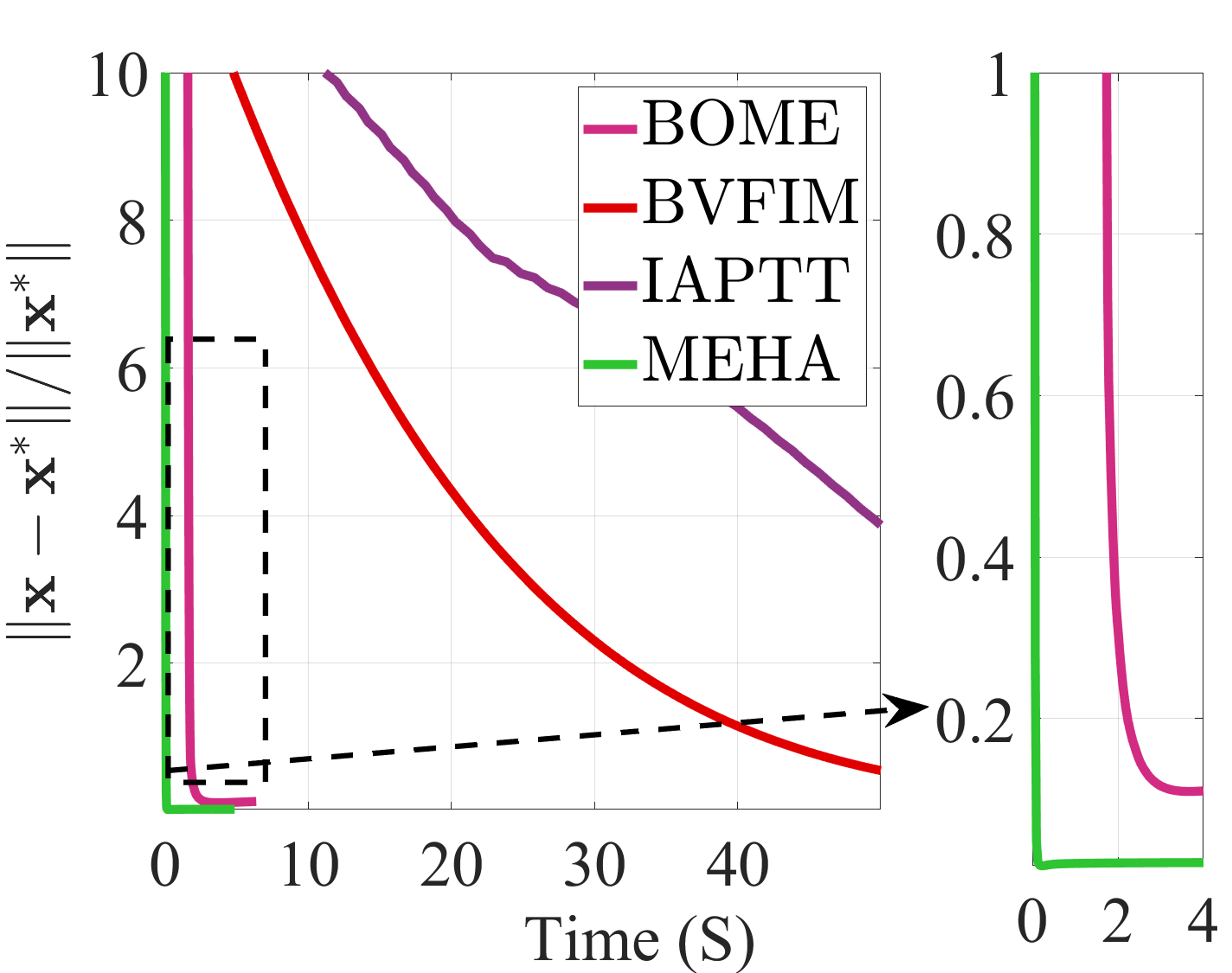}
		&  \includegraphics[width=0.24\textwidth]{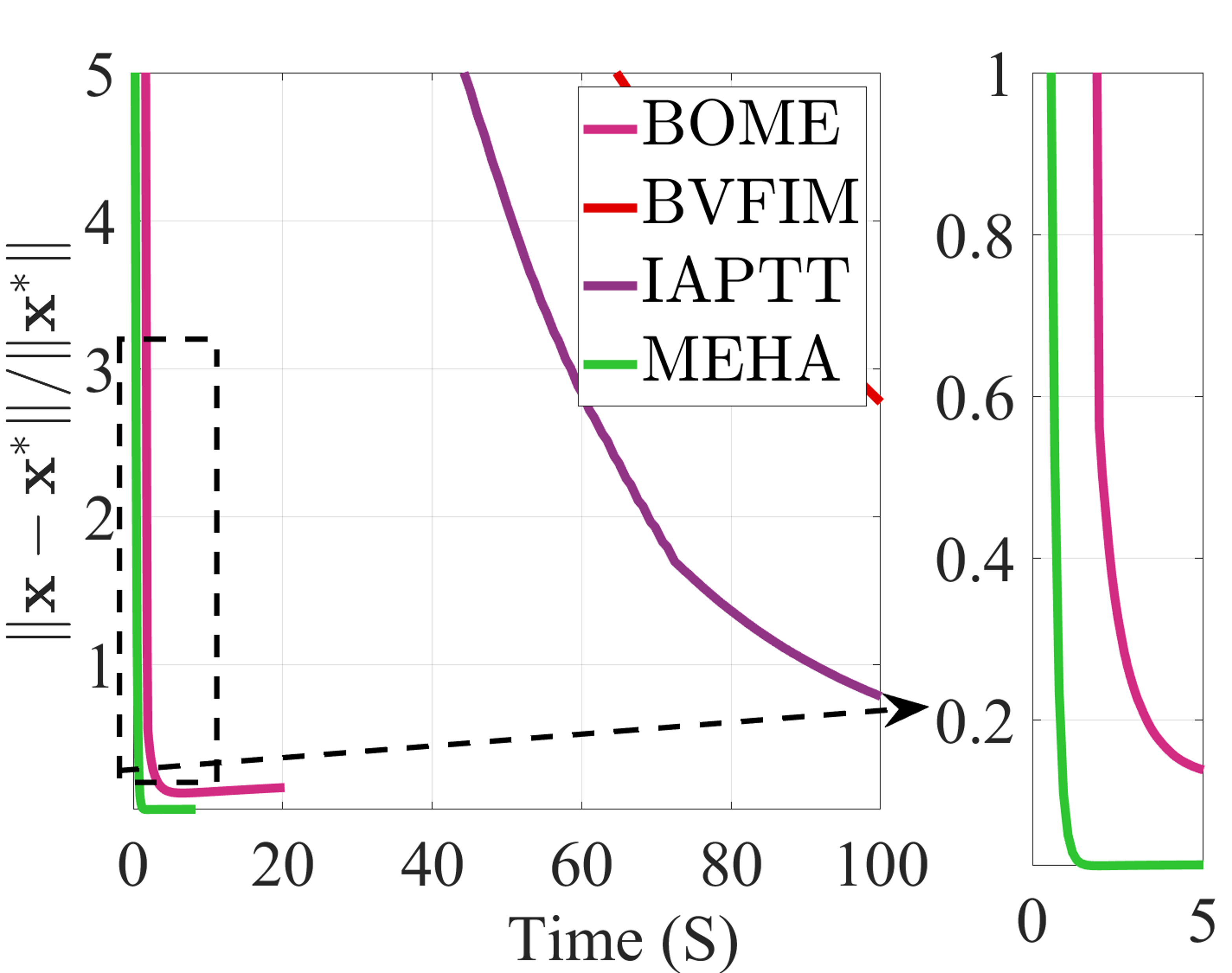}\\
		\includegraphics[width=0.22\textwidth]{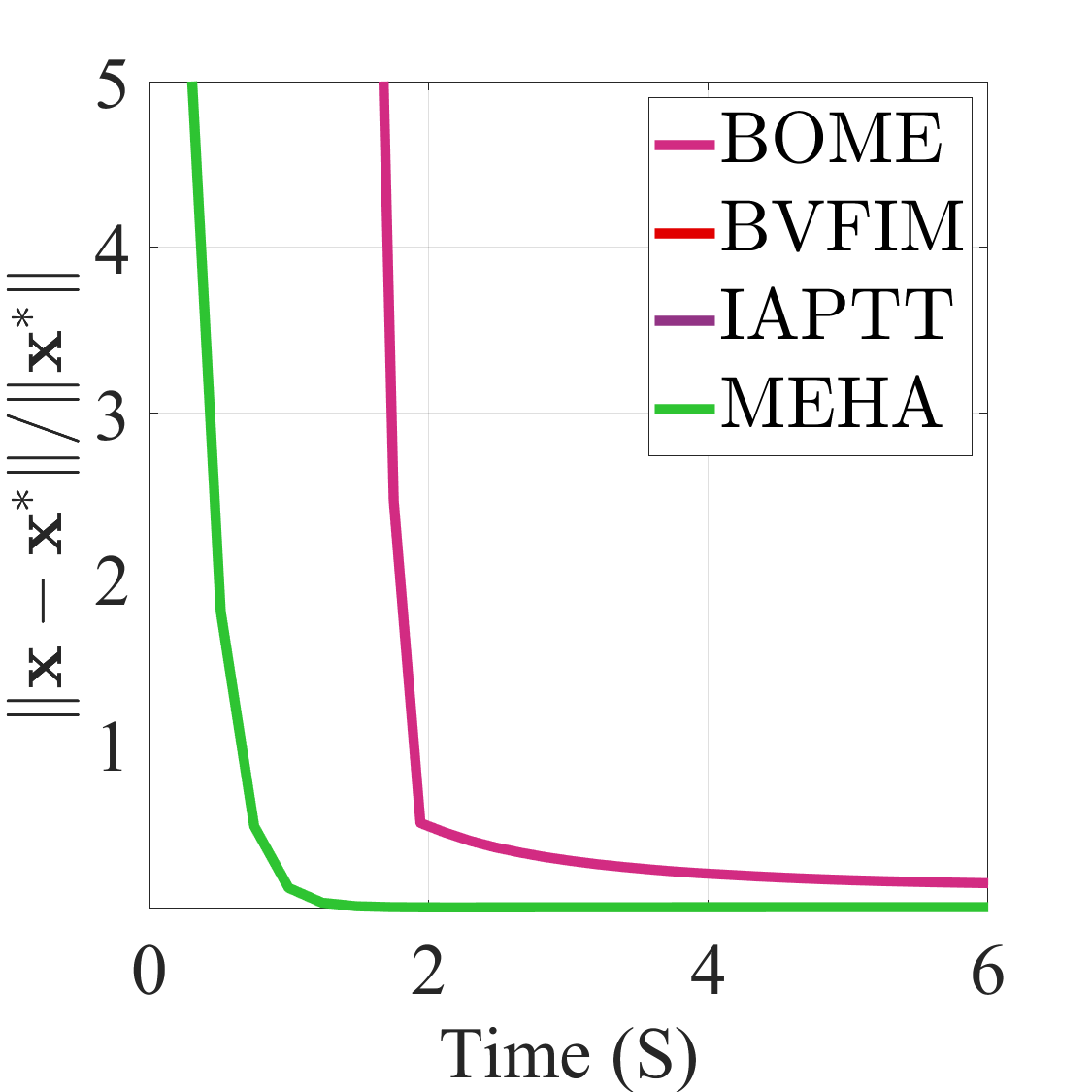}
		&\includegraphics[width=0.22\textwidth]{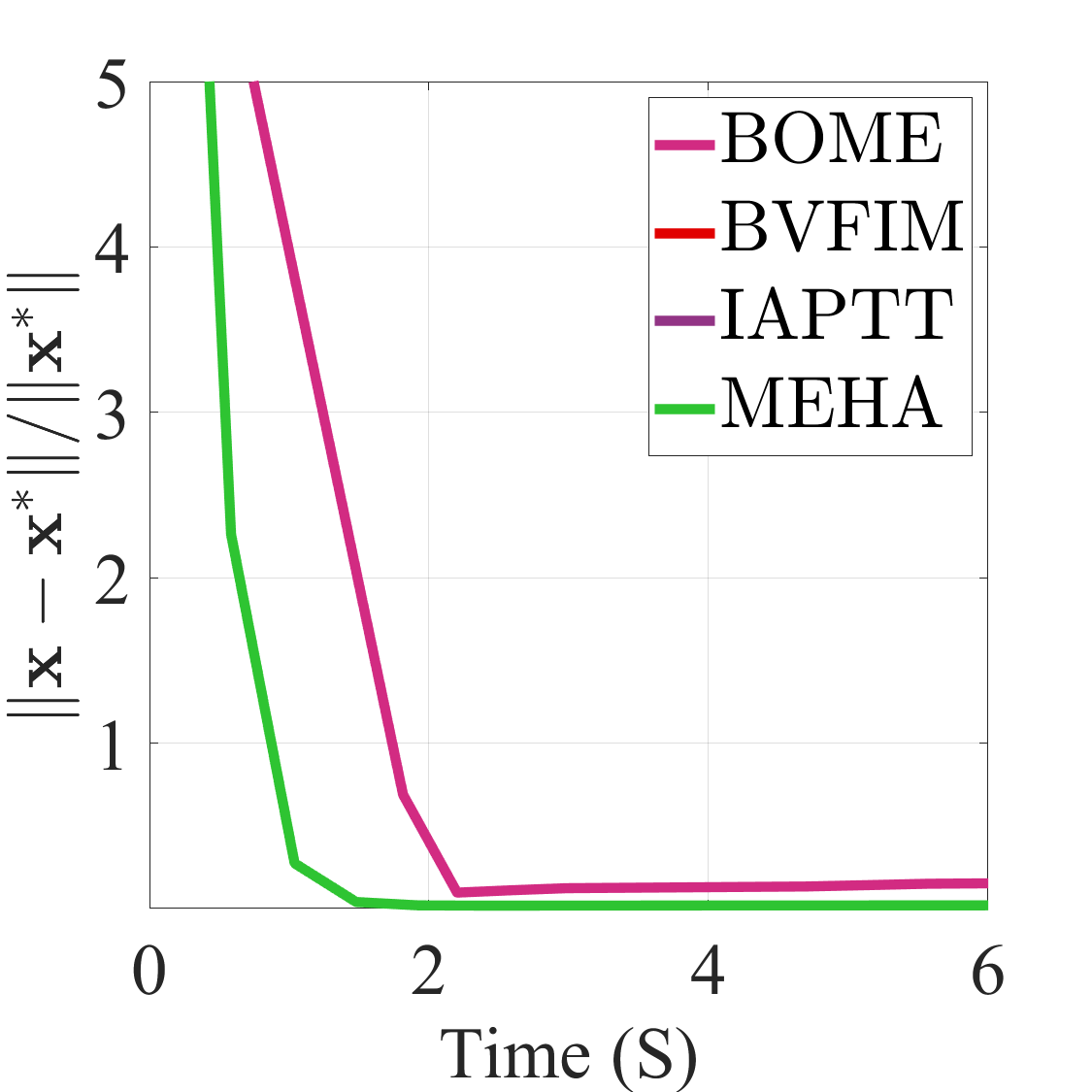}
		\\	
	\end{tabular}
	\vspace{-0.2cm}
	\caption{Convergence curves of advanced BLO methods and MEHA by the criteria, $\|x-x^{*}\|/\|x^{*}\|$, under the LL  non-convex case with different dimensions (10, 50, 100, and 200).}\label{fig:non_convex_dimension}
\end{figure}

\begin{table*}[h!]
	\centering
	\footnotesize
	\renewcommand{\arraystretch}{1.0}
	\caption{Comparison of LL non-smooth case utilizing lasso regression under diverse dimensions.}~\label{tab:non_smooth}	\vspace{-0.2cm}
	\setlength{\tabcolsep}{2.5mm}{
		\begin{tabular}{|ccccc|ccc|ccc|}
			\hline
			\multicolumn{5}{|c|}{ Dimension=2}                                                                                                                            & \multicolumn{3}{c|}{Dimension=100} & \multicolumn{3}{c|}{Dimension=1000}                                                                                                                            \\ \hline
			\multicolumn{1}{|c|}{}   & \multicolumn{1}{c|}{Grid}  & \multicolumn{1}{c|}{Random } & \multicolumn{1}{c|}{TPE}  & MEHA  & \multicolumn{1}{c|}{Random} & \multicolumn{1}{c|}{TPE}  & MEHA & \multicolumn{1}{c|}{Random} & \multicolumn{1}{c|}{TPE}  & MEHA \\ \hline
			\multicolumn{1}{|c|}{Time (S)}  & \multicolumn{1}{c|}{14.87}  & \multicolumn{1}{c|}{17.11} & \multicolumn{1}{c|}{3.32}  & \textbf{1.91} & \multicolumn{1}{c|}{86.27} & \multicolumn{1}{c|}{232.61}   & \textbf{2.34}  & \multicolumn{1}{c|}{700.74} & \multicolumn{1}{c|}{2244.44}   & \textbf{22.83} \\ \hline
			\multicolumn{1}{|c|}{Promotion} & \multicolumn{1}{c|}{$\times$36.97}   & \multicolumn{1}{c|}{$\times$8.96} & \multicolumn{1}{c|}{$\times$1.74}  & - & \multicolumn{1}{c|}{$\times$36.87}  & \multicolumn{1}{c|}{$\times$99.4}   & -  & \multicolumn{1}{c|}{$\times$30.69} & \multicolumn{1}{c|}{$\times$98.3}   & - \\ \hline
		\end{tabular}
	}
\end{table*}

\begin{table}[h!]
	\centering
	\footnotesize
	\renewcommand{\arraystretch}{1.0}
	\caption{Comparison of computation efficiency under large scale dimension on the LL non-convex case.}~\label{tab:large_scale}	\vspace{-0.2cm}
	\setlength{\tabcolsep}{2.5mm}{
		\begin{tabular}{|c|c|c|c|c|}
			\hline
			Dimension & BOME   & BVFIM   & IAPTT  & MEHA  \\ \hline
			200       & 44.47  &3729.74  & 622.77 &\textbf{33.41} \\ \hline
			300       & 48.73  & 4978.18&1164.16  & \textbf{37.85} \\ \hline
			500       & 52.83 &  6639.03& 1657.46 & \textbf{40.25} \\ \hline
			1000       & 68.65 &  8348.92&3537.19  & \textbf{47.45} \\ \hline
		\end{tabular}
	}

\end{table}

\textbf{LL Non-Convex Case.} 
We next demonstrate the effectiveness in handling the LL non-convex case by considering a toy example given in the form of
\begin{eqnarray}\label{eq:LNC}
	\begin{aligned}
		&\min_{x\in \mathbb{R}, y\in \mathbb{R}^n}\|x-a\|^2+\|y-a{\mathbf{e}}-{\mathbf{c}}\|^2 \\
		&\mathrm { s.t. } \quad y_i \in \underset{y_i \in \mathbb{R}}{\arg\min} \ \sin \left(x+y_i-{c}_i\right)\quad
		\forall i,\\
	\end{aligned}
\end{eqnarray} 
where constants $\mathbf{c} \in \mathbb{R}^{n}$ and $a \in \mathbb{R}$. As shown the literature~\cite{liu2021value}, the optimal solution is derived as $x^{*} = \frac{(1-n)a+nC}{1+n}$ and $y_i^{*} = C + c_{i} -x^{*}$ for $i=1,\dots,n$. Here, $C$ is defined as $C=\arg\min_{k}\left\{\left\|C_k-2 a\right\|: C_k=-\frac{\pi}{2}+2 k \pi, k \in \mathbb{Z}\right\}$. The optimal value is $F^{*} = \frac{n(C-2a)^{2}}{1+n}$.
When $n=1$, given $a = 2$ and $\mathbf{c}=2$, the concrete solution is $x^{*}= 3\pi/4$ and $y^{*}= 3\pi/4+2$. Given the initialization point ${x_{0},y_{0}} = (-6,0)$, we compare the performance of MEHA with several Bi-Level Optimization (BLO) schemes (\textit{i.e.,} BVFIM~\cite{liu2021value}, BOME, and IAPTT) across different dimensions in Figure \ref{fig:non_convex_dimension}. Numerical comparisons, particularly the iterative time with advanced competitors, are reported in Table \ref{tab:numerical}. Our method MEHA achieves the fastest convergence speed, outperforming BOME.

\textbf{Computation Efficiency under Large Scale.} We  present the convergence time under large scale setting in Table.~\ref{tab:large_scale}.   MEHA demonstrates rapid convergence across various high-dimensional non-convex scenarios. In 1000 dimensions,  MEHA  achieves a time reduction of 30.9\%.

\begin{table}[h!]
	\centering
	\footnotesize
	\renewcommand{\arraystretch}{1.0}
	\caption{Comparisons of sparse group lasso problem on the synthetic data. $m$ represents the feature dimension.}~\label{tab:lasso}	\vspace{-0.2cm}
	\setlength{\tabcolsep}{1mm}{
		\begin{tabular}{|c|cc|cc|}
			\hline
			\multirow{2}{*}{Settings} & \multicolumn{2}{c|}{$m$=600}                                                  & \multicolumn{2}{c|}{$m$=1200}                                                  \\ \cline{2-5} 
			& \multicolumn{1}{c|}{Test Err.}      & Time (S)       & \multicolumn{1}{c|}{Test Err.}      & Time (S)       \\ \hline
			Grid            & \multicolumn{1}{c|}{84.2$\pm$16.1}          &    4.1$\pm$0.0            & \multicolumn{1}{c|}{87.9$\pm$8.4}          &  11.8$\pm$0.1        \\\hline
			Random         & \multicolumn{1}{c|}{176.5$\pm$29.8}          &  4.5$\pm$1.1            & \multicolumn{1}{c|}{170.4$\pm$17.1}          & 10.7$\pm$0.1          \\ \hline
			TPE        &\multicolumn{1}{c|}{123.9$\pm$26.1}          &  11.7$\pm$2.6             & \multicolumn{1}{c|}{157.2$\pm$24.7}      &  58.2$\pm$6.9     \\ \hline
			IGJO           &\multicolumn{1}{c|}{96.7$\pm$15.6}          &  0.5$\pm$0.3            & \multicolumn{1}{c|}{107.8$\pm$12.2}      &  6.4$\pm$5.2         \\ \hline
			VF-iDCA         &\multicolumn{1}{c|}{74.7$\pm$11.0}          &  0.5$\pm$0.0            & \multicolumn{1}{c|}{85.2$\pm$8.3}      &  8.8$\pm$2.4         \\ \hline
			MEHA             & \multicolumn{1}{c|}{\textbf{71.9$\pm$10.5}}         &  \textbf{0.3$\pm$0.2}      & \multicolumn{1}{c|}{\textbf{84.3$\pm$8.3}} & \textbf{2.3$\pm$0.3} \\ \hline
			\multirow{2}{*}{Settings} & \multicolumn{2}{c|}{$m$=2400}                                                  & \multicolumn{2}{c|}{$m$=3600}                                                  \\  \cline{2-5}
			& \multicolumn{1}{c|}{Test Err.}      & Time (S)     & \multicolumn{1}{c|}{Test Err.}      & Time (S)      \\  \hline
			Grid              & \multicolumn{1}{c|}{89.5$\pm$11.6}          &    25.8$\pm$0.7          & \multicolumn{1}{c|}{91.3$\pm$17.0}          &  38.3$\pm$7.8      \\\hline
			Random             & \multicolumn{1}{c|}{141.9$\pm$29.9}          &    22.3$\pm$1.2            & \multicolumn{1}{c|}{ 132.1$\pm$18.7}          &  42.6$\pm$0.5      \\\hline
			TPE           & \multicolumn{1}{c|}{146.1$\pm$20.3}          &   89.8$\pm$15.9            & \multicolumn{1}{c|}{110.2$\pm$24.0}          &  157.8$\pm$2.4       \\\hline
			IGJO                  & \multicolumn{1}{c|}{115.5$\pm$15.6}          &    12.9$\pm$5.4            & \multicolumn{1}{c|}{120.4$\pm$14.9}          &  24.6$\pm$6.7      \\\hline
			VF-iDCA        & \multicolumn{1}{c|}{85.2$\pm$11.3}          &   16.9$\pm$0.5            & \multicolumn{1}{c|}{90.0$\pm$15.6}          &  49.2$\pm$3.7       \\\hline
			MEHA          &\multicolumn{1}{c|}{\textbf{85.1$\pm$11.4}} & \textbf{2.2$\pm$0.3}  & \multicolumn{1}{c|}{\textbf{89.5$\pm$15.8}} & \textbf{2.5$\pm$0.3} \\ \hline
		\end{tabular}
	}

\end{table}
\textbf{LL Non-Smooth Case.} We validate the computational efficiency of MEHA through a toy lasso regression example:
\begin{eqnarray}\label{eq:LNS}
	\begin{aligned}
		&\min_{x\in \mathbb{R}^{n}, 0 \leq x\leq 1,y\in \mathbb{R}^{n}} \sum_{i=1}^{n} y_{i}\\
		&\mathrm { s.t. }\quad y\in \arg\min_{y\in \mathbb{R}^{n}}  \frac{1}{2}\left\|y- \mathbf{a}\right\|^2+ \sum_{i=1}^{n}x_{i}\|y_{i}\|_{1},\\
	\end{aligned}
\end{eqnarray} 
where $\mathbf{a}:={\left(\frac{1}{n}, \frac{1}{n}, \cdots, \frac{1}{n}\right.}, {\left.-\frac{1}{n},-\frac{1}{n}, \cdots -\frac{1}{n}\right)}\in\mathbb{R}^{n}$. The number of positive and negative values is $\frac{n}{2}$. The optimal solution can be calculated as $x_{i} \in [\frac{1}{n},1]$, $y_{i}=0$ when
$i = 1, \cdots \frac{n}{2}$ and $x_{i} = 0 $, $y_{i}= - \frac{1}{n}$ when
$i = \frac{n}{2}+1, \cdots n$. The solving time of various methods with different dimensions is reported in Table~\ref{tab:non_smooth}.
In comparison to grid search, random search, and Bayesian optimization-based TPE~\cite{bergstra2013making}, our approach demonstrates superior efficiency, requiring the least search time to identify optimal solutions across varied dimensions, particularly notable in large-scale scenarios, such as 1000 dimensions. The convergence curves, depicting the accuracy of MEHA, are illustrated in Figure~\ref{fig:non_smooth} in Appendix~\ref{sec:aer}.
\begin{table}[h!]
	\centering
	\footnotesize
	\renewcommand{\arraystretch}{1.0}
	\caption{Analyzing the sensitivity of parameters combinations.}~\label{tab:sensitivity}	
	\vspace{-0.2cm}
	\setlength{\tabcolsep}{0.5mm}{
		\begin{tabular}{|c|c|c|c|c|c|c|c|c|}
			\hline
			Strategy                                  & $\alpha$ & $\beta$ & $\eta$ & $\gamma$ & $\underline{c}$ &$p$ & Steps & Time (S)\\ \hline
			Original                                  & 0.1        & 0.00001   & 0.1      & 10       & 2   &0.49            & 97    & 22.83    \\ \hline
			\multirow{3}{*}{ $\alpha$}      & 0.05       & 0.00001   & 0.1      & 10       & 2          &0.49     & 109   & 25.55    \\ \cline{2-9}
			& 0.5        & 0.00001   & 0.1      & 10       & 2           &0.49     & 88    & 13.66    \\ \cline{2-9} 
			& 0.8        & 0.00001   & 0.1      & 10       & 2           &0.49     & 86    & 9.65    \\ \hline
			\multirow{3}{*}{ $\beta$}       & 0.1        & 0.00002   & 0.1      & 10       & 2          &0.49     & 91    & 22.01    \\ \cline{2-9}
			& 0.1        & 0.00003   & 0.1      & 10       & 2          &0.49      & 85    & 13.41      \\ \cline{2-9} 
			& 0.1        & 0.000005   & 0.1      & 10       & 2          &0.49      & 105    & 18.24    \\ \hline
			\multirow{3}{*}{ $\eta$}        & 0.1        & 0.00001   & 0.5      & 10       & 2            &0.49   & 88    & 12.82    \\ \cline{2-9} 
			&  0.1        & 0.00001   & 0.8     & 10       & 2         &0.49      & 88   & 9.50    \\ \cline{2-9} 
			& 0.1        & 0.00001   & 0.01     & 10       & 2         &0.49      & 199   & 51.99    \\ \hline
			\multirow{3}{*}{ $\gamma$}        & 0.1        & 0.00001   & 0.1      & 2        & 2         &0.49      & 89    & 19.19    \\ \cline{2-9} 
			& 0.1        & 0.00001   & 0.1      & 20       & 2         &0.49      & 89    & 9.82    \\ \cline{2-9} 
			& 0.1        & 0.00001   & 0.1      & 100      & 2    &0.49           & 90    & 15.3     \\ \hline
			\multirow{3}{*}{ $\underline{c}$} & 0.1        & 0.00001   & 0.1      & 10       & 10   &0.49           & 105   & 17.52    \\ \cline{2-9}
			& 0.1        & 0.00001   & 0.1      & 10       & 50        &0.49      & 107   & 17.22     \\ \cline{2-9} 
			& 0.1        & 0.00001   & 0.1      & 10       & 100        &0.49      & 99   & 10.75    \\ \hline
			\multirow{3}{*}{ $p$}   & 0.1        & 0.00001   & 0.1      & 10       & 2         &0.05    & 956   & 10.53   \\ \cline{2-9} 
			& 0.1        & 0.00001   & 0.1      & 10       & 2         &0.15    & 95   & 10.42    \\ \cline{2-9} 
			& 0.1        & 0.00001   & 0.1      & 10       & 2    &0.30          & 93   & 10.13    \\ \hline
		\end{tabular}
	}

\end{table}
\begin{table*}[!ht]
	\centering
	\footnotesize
	\renewcommand{\arraystretch}{1.0}
	\caption{Comparison of convergence speed for few-shot learning (10-way and 20-way) and accuracy for  data hyper-cleaning tasks.}~\label{tab:meta}	\vspace{-0.2cm}
	\setlength{\tabcolsep}{3.5mm}{
		\begin{tabular}{|c|cc|cc|cc|cc|}
			\hline
			\multirow{2}{*}{Method} & \multicolumn{2}{c|}{10-Way} & \multicolumn{2}{c|}{20-Way} & \multicolumn{2}{c|}{ FashionMNIST} & \multicolumn{2}{c|}{MNIST}  \\ \cline{2-9} 
			& \multicolumn{1}{c|}{Acc. (\%)}    & Time (S) & \multicolumn{1}{c|}{Acc. (\%)}    & Time (S)    	& \multicolumn{1}{c|}{Acc. (\%)}    & F1 score &  \multicolumn{1}{c|}{Acc. (\%)} & \multicolumn{1}{c|}{F1 score}     \\ \hline
			
			RHG	& \multicolumn{1}{c|}{89.77}       &   557.71   &  \multicolumn{1}{c|}{90.17}       & \multicolumn{1}{c|}{550.51}      & \multicolumn{1}{c|}{80.87}       &  88.67  & \multicolumn{1}{c|}{86.93}  &  \multicolumn{1}{c|}{89.46}           \\ \hline
			BDA	& \multicolumn{1}{c|}{89.61}       & 869.98     & \multicolumn{1}{c|}{89.78}     &  1187.63 & \multicolumn{1}{c|}{80.93}       &  87.75   & \multicolumn{1}{c|}{85.83}   & \multicolumn{1}{c|}{89.20}                   \\ \hline
			CG	& \multicolumn{1}{c|}{89.56}       &  363.76 &  \multicolumn{1}{c|}{89.39}     &  602.24      	& \multicolumn{1}{c|}{80.15}       &  88.32   & \multicolumn{1}{c|}{84.02}   & \multicolumn{1}{c|}{81.39}         \\ \hline
			
			BAMM	& \multicolumn{1}{c|}{90.57}       & 180.48   & \multicolumn{1}{c|}{90.13}     &  255.99    & \multicolumn{1}{c|}{81.64}       &   88.16  & \multicolumn{1}{c|}{88.60}  & \multicolumn{1}{c|}{89.91}      \\ \hline
			IAPTT	& \multicolumn{1}{c|}{89.66}       &  299.27  &  \multicolumn{1}{c|}{89.57}     & 3450.80    	& \multicolumn{1}{c|}{81.87}       &  89.25   & \multicolumn{1}{c|}{87.57}   & \multicolumn{1}{c|}{91.10}         \\ \hline
			BOME	& \multicolumn{1}{c|}{89.76}       &  191.25 &  \multicolumn{1}{c|}{89.49}     &  273.19      	& \multicolumn{1}{c|}{82.16}       &  85.01   & \multicolumn{1}{c|}{88.19}   & \multicolumn{1}{c|}{87.95}         \\ \hline
			VPBGD	& \multicolumn{1}{c|}{90.28}       & 315.59 &  \multicolumn{1}{c|}{90.00 }     &  649.93     	& \multicolumn{1}{c|}{81.72}       &  89.93   & \multicolumn{1}{c|}{89.10}   & \multicolumn{1}{c|}{90.64}         \\ \hline
			MEHA	& \multicolumn{1}{c|}{91.17}       &   \textbf{130.81}   & \multicolumn{1}{c|}{89.96}     & \textbf{ 235.77}    & \multicolumn{1}{c|}{\textbf{83.19}}       &   \textbf{89.94}  & \multicolumn{1}{c|}{\textbf{89.55}}  & \multicolumn{1}{c|}{\textbf{91.26}}        \\ \hline
		\end{tabular}
	}

\end{table*}
\begin{table*}[!h]
	\centering
	\footnotesize
	\renewcommand{\arraystretch}{1.0}
	\caption{Comparing Top-1 accuracy in searching, inference, and final test stages for DARTS~\cite{liu2018darts}, P-DARTS~\cite{chen2019progressive}, PC-DARTS~\cite{xu2019pc1}, and typical  BLO schemes.}~\label{tab:nas}	\vspace{-0.2cm}
	\setlength{\tabcolsep}{6mm}{
		\begin{tabular}{|c|cc|cc|c|c|c|}
			\hline
			\multirow{2}{*}{Methods} & \multicolumn{2}{c|}{Searching} & \multicolumn{2}{c|}{Inference} & \multirow{2}{*}{Test} & \multirow{2}{*}{Params (M)} \\ \cline{2-5}
			& \multicolumn{1}{c|}{Train}    &  Valid    & \multicolumn{1}{c|}{Train}    &  Valid    &                       &                             \\ \hline
			DARTS & \multicolumn{1}{c|}{98.320}    &  88.940   & \multicolumn{1}{c|}{99.481}    &  95.639    &               95.569        &  1.277                           \\ \hline
			P-DARTS & \multicolumn{1}{c|}{96.168}    &  90.488   & \multicolumn{1}{c|}{99.802}    &  95.701    &      95.710                  &  1.359                           \\ \hline
			PC-DARTS & \multicolumn{1}{c|}{84.821}    &  83.516   & \multicolumn{1}{c|}{98.163}    &  95.630    &      95.540                  &  1.570                           \\ \hline
			RHG & \multicolumn{1}{c|}{98.448}    &  89.556   & \multicolumn{1}{c|}{99.688}    &  95.340    &      95.340                  &  1.359                           \\ \hline
			CG & \multicolumn{1}{c|}{99.126}    & 89.298    & \multicolumn{1}{c|}{98.909}    & 95.499     &    95.370                    &  1.268                           \\ \hline
			IAPTT & \multicolumn{1}{c|}{98.904}    &  99.512    & \multicolumn{1}{c|}{99.776}    &  95.840   &           95.809            &  1.963                           \\ \hline
			
			MEHA & \multicolumn{1}{c|}{99.060}    &  \textbf{99.764}   & \multicolumn{1}{c|}{99.419}    & \textbf{96.150}    &    \textbf{96.070}                   &    1.524                         \\ \hline
		\end{tabular}
	}
	
\end{table*}

Furthermore, we assess the algorithm's performance on the more challenging group lasso hyperparamter selection problem using synthetic data, formulated as follows:
\begin{eqnarray}\label{eq:lass}
	\begin{aligned}
		&\min_{x\in \mathbb{R}^{n}_+,y\in \mathbb{R}^m}   l_{val}(y)\\
		&\mathrm { s.t. }\quad y\in \arg\min_{y\in \mathbb{R}^{m}}   l_{tr}(y) + \sum_{i=1}^{n}x_{i}\|y^{(i)}\|_{2},\\
	\end{aligned}
\end{eqnarray} 
where $ l_{val}$ and $l_{tr}$ denote the loss $\frac{1}{2} \sum_{i} \left|b_{i}-y^{T}a_{i}\right|^2$ on the validation and training datasets, respectively. Here, $x$ denotes hyperparameters with n dimensions, while $\mathbf{a}$ and $\mathbf{b}$ stand for the inputs and their corresponding labels, respectively. Further details regarding data generation are elaborated in Section~\ref{sec:detail}. Competitive methods, including implicit differentiation IGJO~\cite{feng2018gradient} and VF-iDCA~\cite{gao2022value}, are included for comparison. Numerical results repeated five times are presented in Table~\ref{tab:lasso}.  Our proposed MEHA outperforms in high-dimensional tasks, yielding the most precise results.

\textbf{Sensitivity of Parameters.}  In investigating the sensitivity of the parameters in MEHA, we conducted additional experiments in the first instance of the LL non-smooth Case (\textit{i.e.,} Eq.\eqref{eq:LNS}), with a dimension of 1000. Convergence time and steps are detailed in Table \ref{tab:sensitivity} while altering the value of only one parameter, and keeping the values of others. It is crucial to emphasize that the algorithm achieves convergence and is stable across diverse parameter combinations. Additional results can be found in Table~\ref{tab:sensitivitymore} in Section \ref{sec:aer}.

\subsection{Real-world Applications}

\textbf{Few-Shot Learning.} The goal of N-way M-shot classification is to improve the adaptability of learnable model, facilitating rapid adaptation to new tasks. In our experiments using the Omniglot dataset~\cite{finn2017model}, specifically in 10-way 1-shot and 20-way 1-shot scenarios with a LL convex formulation, we provide a runtime comparison in Table~\ref{tab:meta} to achieve consistent accuracy levels (90\%) for both scenarios. Notably, our method attains comparable accuracy while significantly reducing computational time.

\begin{figure}[htb]
	\footnotesize
	\centering
	\setlength{\tabcolsep}{1pt}
	\begin{tabular}{cc}
		\includegraphics[width=0.23\textwidth]{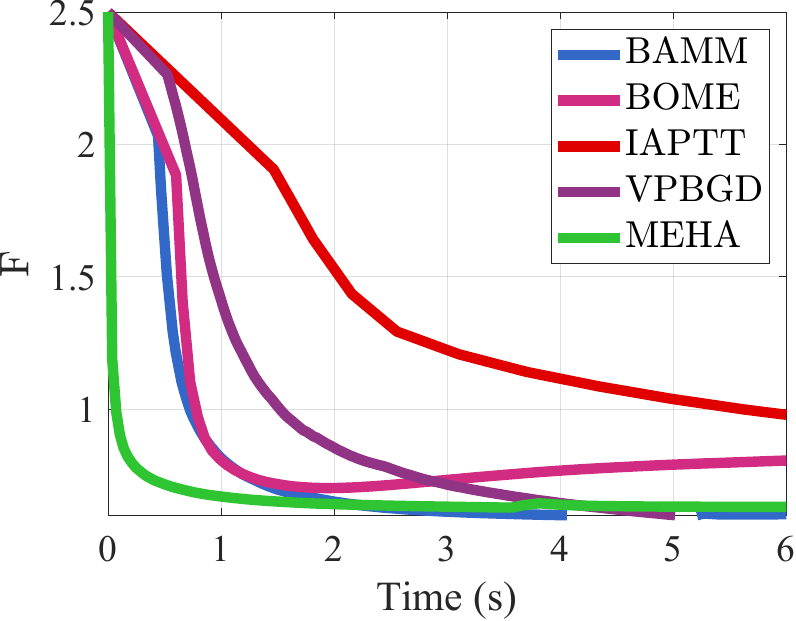}
		&\includegraphics[width=0.23\textwidth]{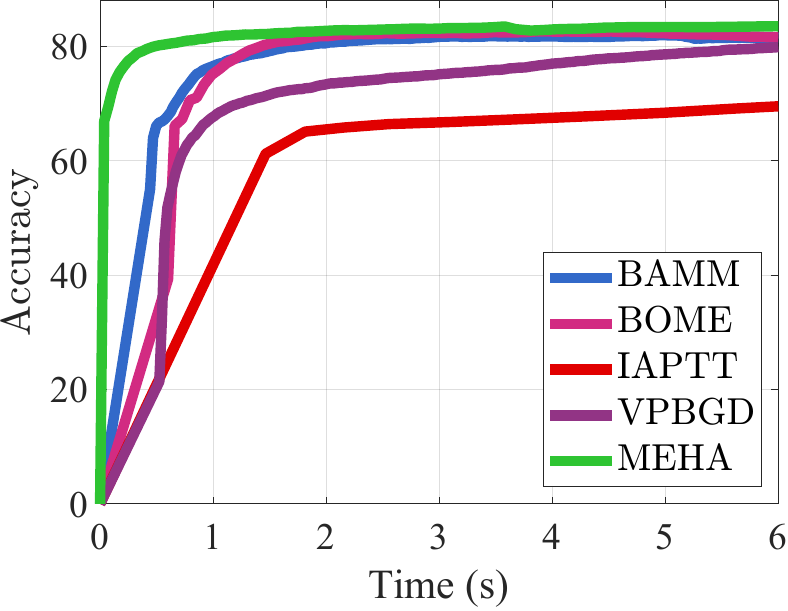}
		\\	
	\end{tabular}
	\vspace{-0.2cm}
	\caption{ Comparison of data hyper-cleaning on FashionMNIST.}
	\label{fig:data_clean}
\end{figure}

\textbf{Data Hyper-Cleaning.} 
In the right section of Table~\ref{tab:meta}, we depict the  accuracy and F1-score, for diverse methods  (LL nonconvex case) on FashionMNIST and MNIST datasets. Remarkably, our proposed method significantly attains the desired solution with the highest accuracy. Furthermore, Figure~\ref{fig:data_clean} visually represents the validation loss and test accuracy across various methods. It is evident that our approach not only demonstrates the fastest convergence rate but also maintains its superior performance.

{\bf Neural Architecture Search.}  Our specific focus  on differentiable NAS methods, representing a LL non-convex case. To showcase consistent performance, accuracy results are presented across different stages in Table~\ref{tab:nas}. Remarkably, MEHA consistently outperforms specialized designs for NAS, establishing its superiority. Additionally, a comparison with existing solutions of BLO methods is provided, demonstrating our superiority in addressing LL non-convex applications.  Note that, our theory can handle the neural architecture search problem with smooth activation functions, such as  Swish, where the results  are reported in the   \ref{sec:aer}.

\section{Conclusions}
By utilizing the Moreau envelope-based reformulation for general nonconvex and nonsmooth BLOs, we propose a provably single-loop and Hessian-free gradient-based algorithm, named MEHA. 
We validate the effectiveness and efficiency of MEHA for large-scale nonconvex-nonconvex BLO in both synthetic problems and various practical machine learning tasks. 
By leveraging the simplicity of our approach 
and integrating techniques such as variance reduction and momentum, 
we would be interested in studying the stochastic algorithms in the future.
\nocite{langley00}

\section*{Impact Statements}
This paper presents work whose goal is to advance the field of Machine Learning. There are many potential societal consequences of our work, none which we feel must be specifically highlighted here.

\section*{Acknowledgements}
Authors listed in alphabetical order.  This work is partially supported by the National Key R\&D Program of China (No. 2022YFA1004101), the National Natural Science Foundation of China (Nos. U22B2052, 12326605, 12222106, 12326605, 62331014, 12371305), Guangdong Basic and Applied Basic Research Foundation (No. 2022B1515020082).

\nocite{langley00}

\bibliography{example_paper}
\bibliographystyle{icml2024}

\newpage
\appendix
\onecolumn
\section{Appendix}

The appendix is organized as follows:
\begin{itemize}
	\item Expanded related work is provided in Section \ref{relatedwork2}.
	
	\item We highlight the proof sketch of Theorem \ref{prop1} in Section \ref{proofsketch}.

	\item Additional experimental results are provided in Section \ref{sec:aer}.
	\item Experimental details are provided in Section \ref{sec:detail}.
	\item The equivalent result of the reformulated problem \ref{wVP} is provided in Section \ref{equiv-a}.
	
	\item The convergence of approximation problem (\ref{problem_pen}) is established in Section \ref{conv_appro}.
	
	\item We prove the weakly convexity and derive the gradient formula of Moreau Envelope in Section \ref{propertiesME}.
	
	\item Some useful auxiliary lemmas are provided in Section \ref{lemmas}.
	
	\item The proof of Lemma~\ref{lem9} is given in Section \ref{proofV}.
	
	\item The proof of Theorem \ref{prop1} is provided in Section \ref{proofofthm}.
	
	\item The proof of of Theorem \ref{prop2} is provided in Section \ref{relation}.
	
	\item We verify Assumptions \ref{assump-LL}(ii) and (iii) are special cases of Assumption \ref{assump-LL}(iv) in Section \ref{A9-assump}.
	
	\item The single-loop Moreau envelope based Hessian-free Algorithm in the smooth case is provided in Section \ref{smoothcase}.
\end{itemize}

\subsection{Expanded Related Work}
\label{relatedwork2}

In this section, we provide an extensive review of recent studies closely related to ours.

{\bf Nonconvex-Convex BLO.} 
The LL strong convexity significantly contributes to the development of efficient BLO algorithms, see, e.g.,  \cite{ghadimi2018approximation,ji2020bilevel,chen2021closing,ji2022will,hong2020two,kwon2023fully}. 
It guarantees the uniqueness of the LL minimizer (Lower-Level Singleton), which facilitates the demonstration of asymptotic convergence for the iterative differentiation-based approach \cite{franceschi2018bilevel}. 
If further the LL objective is twice differentiable, the gradient of the UL objective (hyper-gradient) can be expressed using the implicit function theorem.
Then the uniformly LL strong convexity implies both the smoothness and the Lipschitz continuity properties of the LL solution mapping. 
These essential properties facilitates the demonstration of non-asymptotic convergence for both the iterative differentiation and the approximate implicit differentiation approaches with rapid convergence rates, see e.g., \cite{ghadimi2018approximation,ji2020bilevel,chen2021closing,sow2022convergence,ji2022will,ji2021lower,arbel2022amortized,li2022fully,dagreou2022framework,hong2020two,yang2023accelerating}. 
Due to the implicit gradient, the methods mentioned above necessitate costly manipulation involving the Hessian matrix, making them all second-order methods. 
Recently, \cite{kwon2023fully} developed  stochastic and deterministic fully first-order BLO algorithms based on the value function approach \cite{ye1995optimality}, and established their non-asymptotic convergence guarantees, while an improved convergence analysis is provided in the recent work \cite{chen2023near}. By using a projection-aided finite-difference Hessian/Jacobian-vector approximation, and momentum-based updates, a simple fully single-loop
Hessian/Jacobian-free stochastic BLO algorithm has been proposed in \cite{yang2023achieving} with an $\widetilde O(\epsilon^{-1.5})$ sample complexity. 

In the absence of strong convexity, additional challenges may arise, including the presence of multiple LL solutions (Non-Singleton), which can hinder the application of implicit-based approaches involved in the study of nonconvex- strongly-convex BLOs.
To tackle Non-Singleton,  sequential averaging methods (also referred to as aggregation methods) were proposed in \cite{liu2020generic,li2020improved,liu2022general,pmlr-v202-liu23y}.
Recent advances include value function based difference- of-convex algorithm \cite{gao2022value,ye2023difference}; primal-dual algorithms \cite{sow2022constrained}; first-order penalty methods using a novel minimax optimization reformulation \cite{lu2023first}.

{\bf Nonconvex-Nonconvex BLO.} 
While the nonconvex-convex BLO has been extensively studied in the literature, the efficient methods for nonconvex-nonconvex BLO remain under-explored.
Beyond the LL convexity, the authors in 
\cite{liu2021towards} develop a method with initialization auxiliary and pessimistic trajectory truncation; 
the study \cite{arbel2022non} extends implicit differentiation to a class of nonconvex LL functions with possibly degenerate critical points and then develops unrolled optimization algorithms.
However, these works requires second-order gradient information and do not provide finite-time convergence guarantees. 
Still with the second-order gradient information but providing non-asymptotic analysis, the recent works \cite{huang2023momentumbased} and \cite{xiao2023generalized} propose a momentum-based BLO algorithm and a generalized alternating method for BLO with a nonconvex LL objective that satisfies the Polyak-Łojasiewicz (PL) condition, respectively. 
In contrast to these methods discussed above,
the value function reformulation of BLO was firstly utilized in  \cite{liu2021value} to develop BLO algorithms in machine learning, using an interior-point method combined with a smoothed approximation. 
But it lacks a complete non-asymptotic analysis.
Subsequently, \cite{NeurIPS2022-Liu} introduced a fully first-order value function based BLO algorithm. 
They also established the non-asymptotic convergence results when the LL objective satisfies the PL or local PL conditions. 
Recently, \cite{ShenC23} proposed a penalty-based fully first-order BLO algorithm and established its finite-time convergence under the PL conditions. Notably, this work relaxed the relatively restrictive assumption on the boundedness of both the UL and LL objectives that was present in \cite{NeurIPS2022-Liu}. Other recent advances include penalty method and first-order stochastic approximation in \cite{kwon2023penalty}, which primarily explores BLO from the perspective of the hyper-objective, using a penalty value function-based approach; efficient adaptive projection-aid gradient methods based on mirror descent in \cite{huang2023adaptive} for both deterministic and stochastic BLO problems.

{\bf Nonsmooth BLO.} 
Despite plenty of research focusing on smooth BLOs, there are relatively fewer studies addressing nonsmooth BLOs, see, e.g., \cite{mairal2011task,okuno2018hyperparameter,bertrand2020implicit,bertrand2022implicit}. 
However, these works typically deal with special nonsmooth LL problems, e.g., 
task-driven dictionary learning with elastic-net (involving the $\ell_1$-norm) in \cite{mairal2011task}; 
the Lasso-type models (including the $\ell_1$-norm as well) for hyper-parameter optimization in \cite{bertrand2020implicit}; 
$\ell_p$-hyperparameter learning with $0<p<1$ in \cite{okuno2018hyperparameter};
non-smooth convex learning with separable non-smooth terms in \cite{bertrand2022implicit}. 
Recently, there is a number of works studying BLOs with general nonsmooth LL problems. 
By decoupling hyperparameters from the regularization, based on the value function approach, \cite{gao2022value} develop a sequentially convergent Value Function-based Difference-of-Convex Algorithm with inexactness for a specific class of bi-level hyper-parameter selection problems. 
\cite{gao2023moreau} introduces a Moreau envelope-based reformulation of BLOs and develops an inexact proximal Difference-of-weakly-Convex algorithm with sequential convergence, to substantially weaken the underlying assumption in \cite{ye2023difference} from lower level full convexity to weak convexity.
There is also a line of works devoted to tackle the nonsmooth UL setting, including: 
Bregman distance-based method in \cite{huang2021enhanced}; 
proximal gradient-type algorithm in \cite{chen2022fast}.

\subsection{Proof Sketch}\label{proofsketch}

The proof of Theorem \ref{prop1}, which establishes the non-asymptotic convergence, relies on the monotonically decreasing property of the merit function:
\begin{equation*}
	V_k := \phi_{c_k}(x^k, y^{k}) +  
	C_V
	\left\| \theta^{k} - \theta_{\gamma}^*(x^k,y^{k}) \right\|^2,
\end{equation*}
where $\theta_{\gamma}^*(x,y)$ is the unique solution to problem (\ref{solutionset1}), the constant $C_V:=(L_f + L_g)^2+ 1/\gamma^2 $, and 
\begin{equation}\label{def_varphi}
	\phi_{c_k}(x,y) := \frac{1}{ c_k}\big(F(x,y) - \underline{F} \big) +   (f+ g- v_\gamma) (x,y).
\end{equation}

\begin{lemma}\label{lem9}
	Under Assumptions \ref{assump-UL} and \ref{assump-LL}, 
	suppose $\gamma \in (0, \frac{1}{2\rho_{f_2} + 2\rho_{g_2} })$, $c_{k+1} \ge c_k$ and $\eta_k \in [ \underline{\eta}, (1/\gamma - \rho_{f_2})/(L_f+1/\gamma)^2] \cap  [ \underline{\eta}, 1/\rho_{g_2})$ with $ \underline{\eta} > 0$, 
	then there exists $c_{\alpha}, c_\beta, c_{\theta}> 0$ such that when $\alpha_k \in (0, {c_\alpha}]$ and $\beta_k \in (0, c_\beta]$, the sequence of $(x^k, y^k, \theta^k)$ generated by MEHA (Algorithm \ref{MEHA}) satisfies
	\begin{equation}\label{lem9_eq}
		\begin{aligned}
			V_{k+1} - V_k 
			\le \, & 
			- c_{\theta}
			\left\| \theta^{k} - \theta_{\gamma}^*(x^{k},y^{k}) \right\|^2\\ &
			-  \frac{\| x^{k+1} - x^k \|^2}{4\alpha_k} 
			- \frac{\| y^{k+1} - y^k \|^2}{4\beta_k} .
		\end{aligned}
	\end{equation}
\end{lemma}

We 
outline the pivotal steps leading to Lemma \ref{lem9}.

{\bf Step 1:} First, we consider the Moreau envelope $v_\gamma(x,y)$, focusing on two of its critical properties: its weak convexity (referenced in Lemma \ref{Lem1-a}) and the associated gradient formulas (outlined in Lemma \ref{lem2-a}). These properties enable us to derive an upper bound for the descent of the penalized objective value $\phi_{c_k}(x,y)$ with incorporating the error term $\left(\frac{\alpha_k}{2} (L_f + L_g)^2 + \frac{\beta_k}{\gamma^2} \right) \left\| \theta^{k+1} - \theta_{\gamma}^*(x^k,y^{k}) \right\|^2$, as formulated in Equation (\ref{lem6_eq}) in Lemma \ref{lem6}.

{\bf Step 2:} The subsequent step involves leveraging the Lipschitz continuity of the Moreau envelope solution $\theta_\gamma^*(x,y)$  (as per Lemma \ref{lem_theta}) along with the contraction properties of $\theta^k$ towards this solution (discussed in Lemma \ref{lem4}). This approach is instrumental in controlling the aforementioned error term $\left\| \theta^{k+1} - \theta_{\gamma}^*(x^k,y^{k}) \right\|^2$ as presented in Equation (\ref{lem6_eq}).

Ultimately, these steps culminate in confirming the monotonically decreasing property of the merit function $V_k$ as in Lemma \ref{lem9}. The detailed proof is provided in Section \ref{proofV}.

	\subsection{Additional Experimental Results}~\label{sec:aer}

	\textbf{LL Strong Convex Case.} 
	We initially present the convergence results using the toy numerical problem introduced in BDA~\cite{liu2020generic} with a lower-level convex objective, expressed as:
	\begin{equation}
		\min _{x\in \mathbb{R}^n} \frac{1}{2}\left\|x-z_0\right\|^2+\frac{1}{2} y^*(x)^{\top} {A} y^*(x) 
		\  \text { s.t. } 
		\  y^*(x)=\arg \min _{y\in \mathbb{R}^n} f(x, y)=\frac{1}{2} y^{\top} {A} y-x^{\top} y,
	\end{equation}
	where $x\in \mathbb{R}^{n}$ and $y\in \mathbb{R}^{n}$. We define ${A}$ has the positive-definite symmetric property and ${A} \in \mathbb{S}^{n\times n}$, ${z}_{0} \neq 0$ and $ {z}_{0} \in \mathbb{R}^{n}$.  Concretely, we set ${A} = {I}$ and ${z}_{0} = {e}$. Thus, the optimal solution is $x^{*} = y^{*} = {e}/2$, where $\mathbf{e}$ represents the vector containing all elements equal to one.
	
	This case notably aligns with several convergence assumptions of BLO methods, covering Explicit Gradient-Based Methods (EGBMs) like RHG~\cite{franceschi2017forward}, BDA, and IAPTT~\cite{liu2021towards}, Implicit Gradient-Based Methods (IGBMs) such as CG~\cite{pedregosa2016hyperparameter} and NS~\cite{rajeswaran2019meta}, and contemporary proposed methods (BRC~\cite{liu2023learningBRC}, BOME~\cite{NeurIPS2022-Liu}, F2SA~\cite{kwon2023fully}, BAMM~\cite{pmlr-v202-liu23y}). Detailed numerical comparisons are provided in Table~\ref{tab:strongconvex} in the Appendix, with visual comparisons in Figure~\ref{fig:sc_trad}. Analyzing the behaviors in Figure~\ref{fig:sc_trad}, our method exhibits the fastest convergence among EGBMs, IGBMs, and single-loop methods. From Table~\ref{tab:strongconvex}, it is evident that our method achieves two significant improvements. In comparison to the effective BAMM, our method demonstrates an 85.8\% improvement in inference time. Furthermore, the proposed scheme incurs the lowest computational cost, utilizing only 10.35\% of the memory used by BAMM.
	\begin{figure}[htb]
		\footnotesize
		\centering
		\setlength{\tabcolsep}{1pt}

		\begin{tabular}{cccc}
			\includegraphics[width=0.24\textwidth]{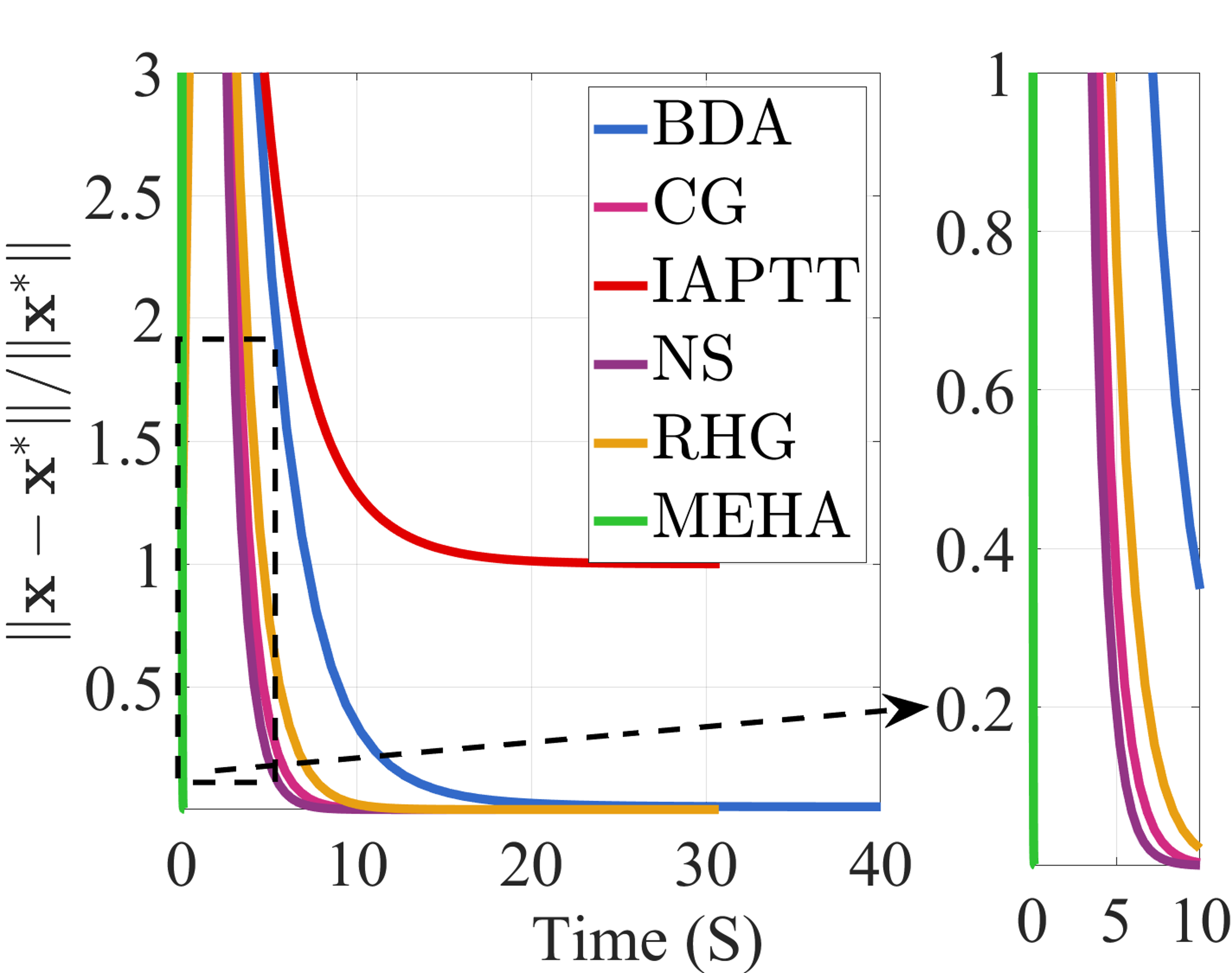}
			&\includegraphics[width=0.24\textwidth]{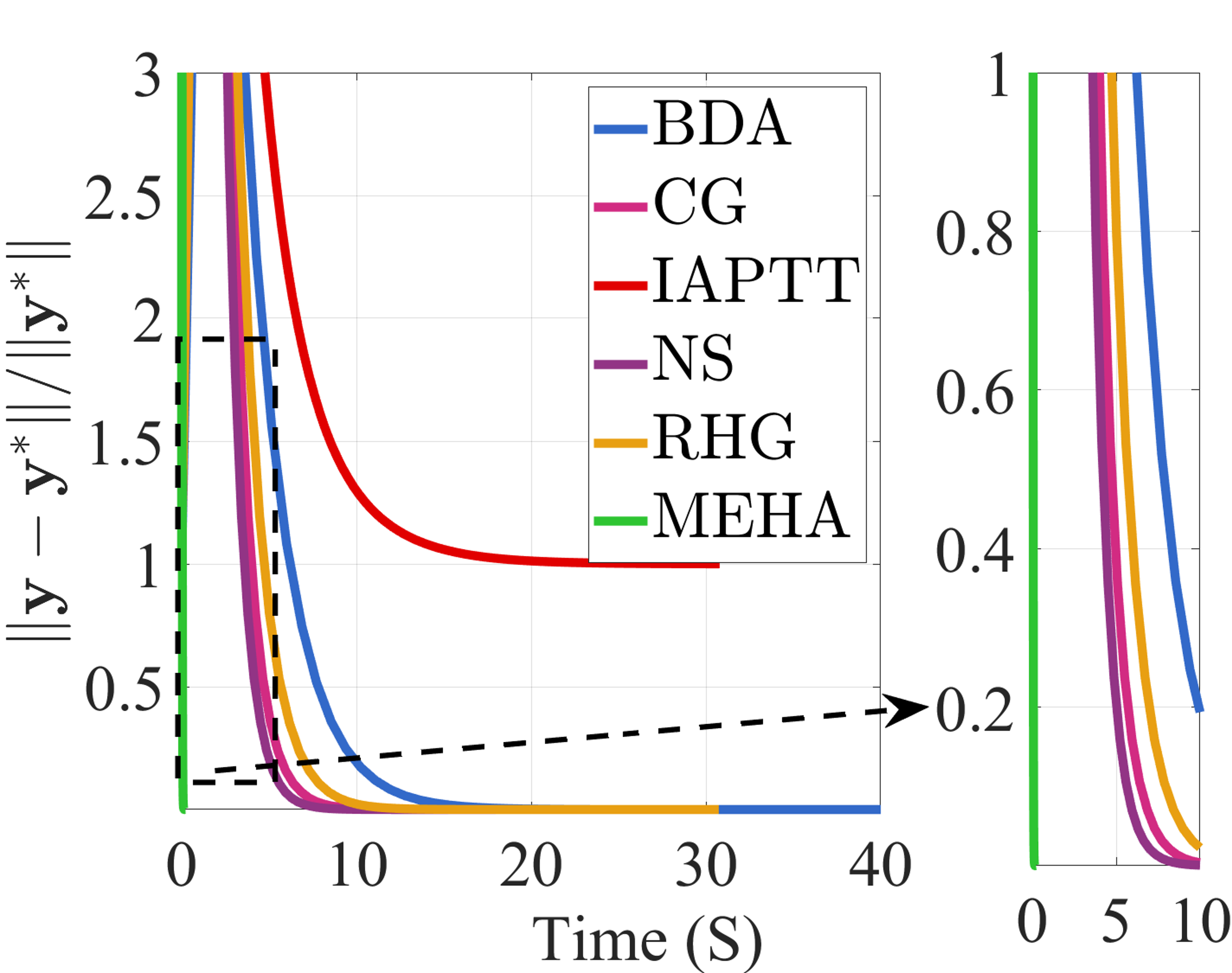}
			&	\includegraphics[width=0.24\textwidth]{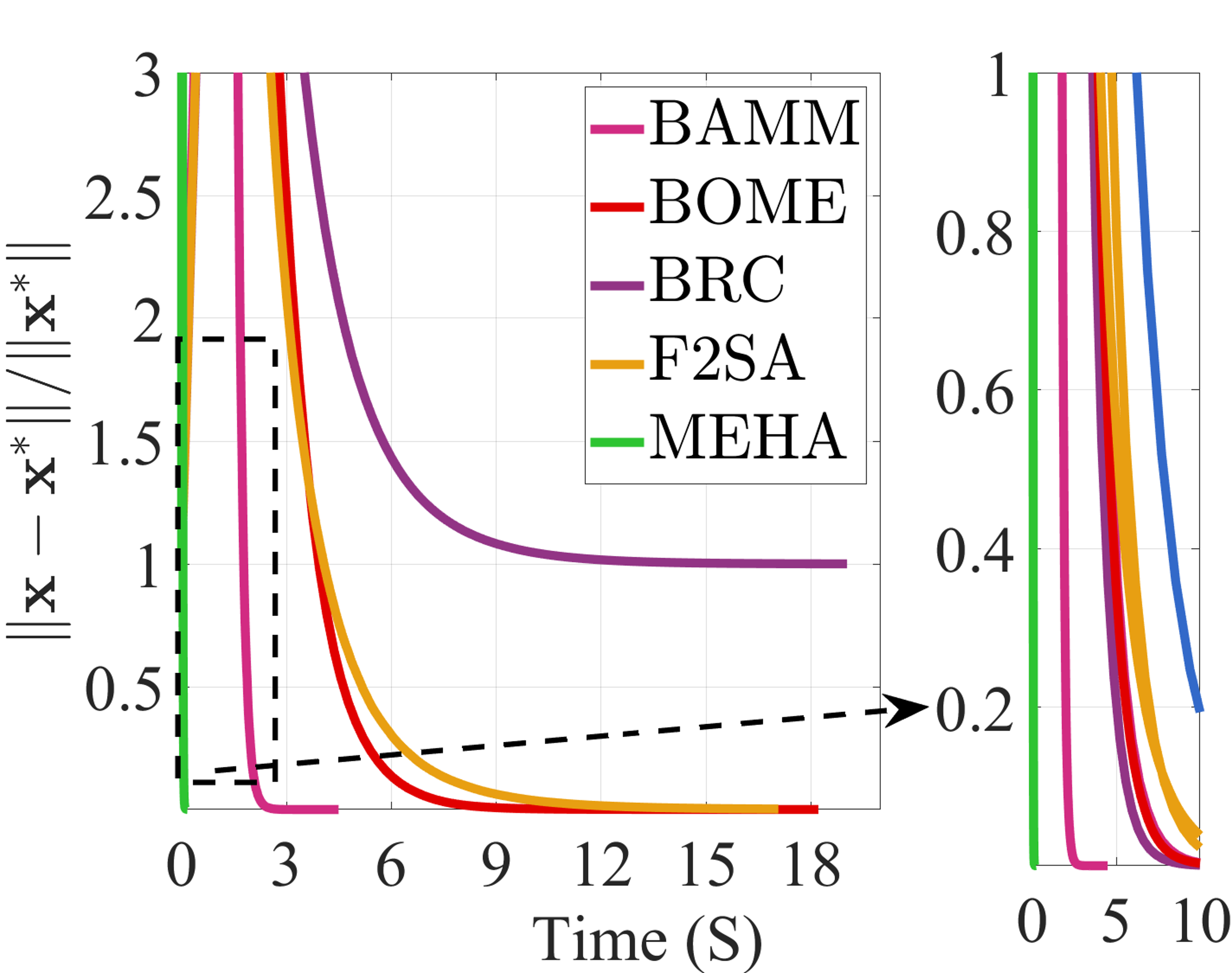}
			
			&\includegraphics[width=0.24\textwidth]{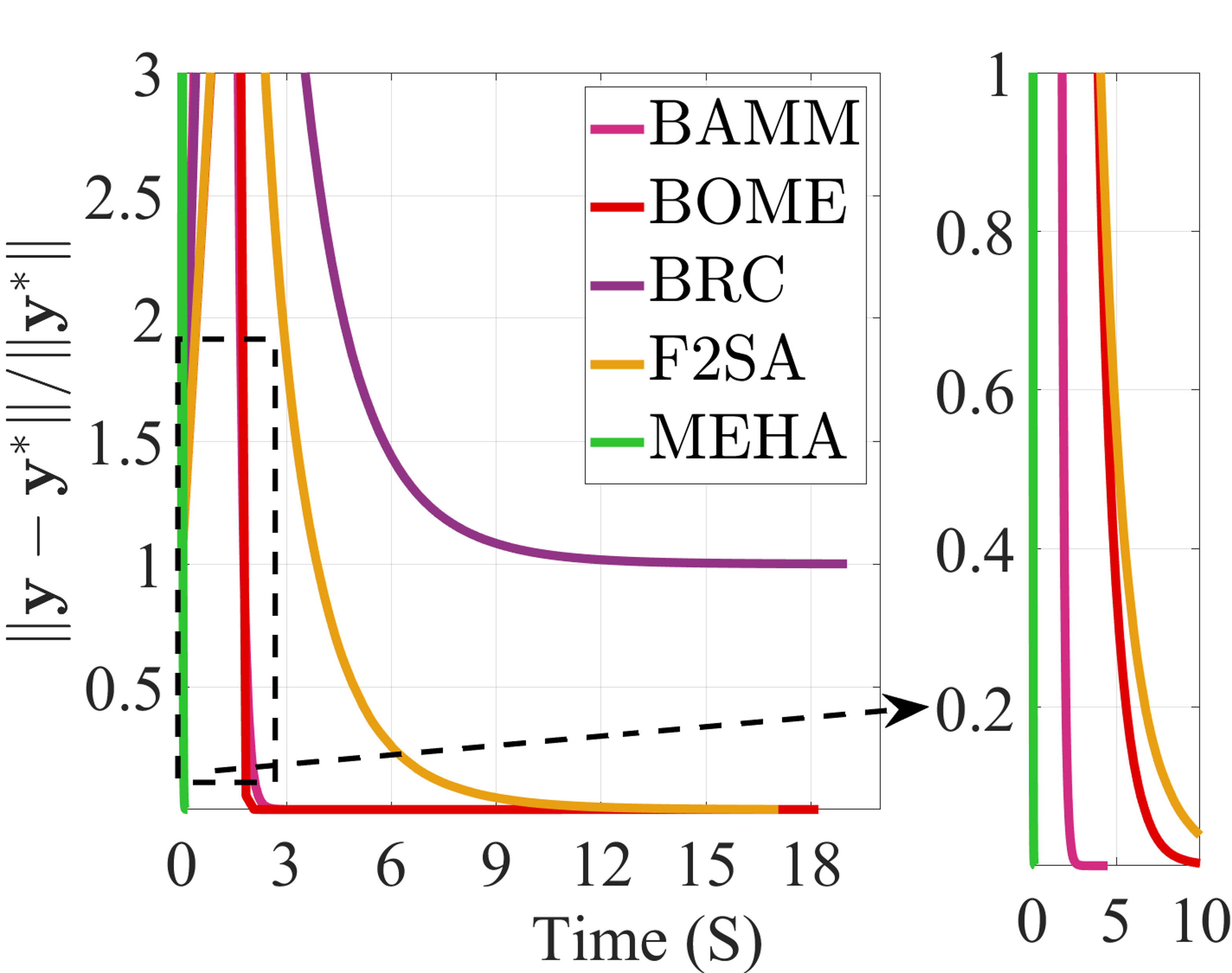}
			\\	
		\end{tabular}
		\vspace{-0.2cm}
		\caption{ Illustrating the convergence curves of advanced BLO methods and MEHA by the criterion of $\|x-x^{*}\|/\|x^{*}\|$ and $\|y-y^{*}\|/\|y^{*}\|$ under  LL strong  convex case.}\label{fig:sc_trad}
	\end{figure}
	\begin{table}[htb]
		\centering
		\footnotesize
		\renewcommand{\arraystretch}{1.0}
		\caption{Basic properties of   time and memory  under  LL strong convex case.}~\label{tab:strongconvex}	\vspace{-0.2cm}
		\setlength{\tabcolsep}{1.6mm}{
			\begin{tabular}{|c|ccc|cc|ccccc|}
				\hline
				Category     & \multicolumn{3}{c|}{EGBMs}                            & \multicolumn{2}{c|}{IGBMs}     & \multicolumn{5}{c|}{Others}                                                    \\ \hline
				Methods & \multicolumn{1}{c|}{RHG} & \multicolumn{1}{c|}{BDA} & IAPTT & \multicolumn{1}{c|}{CG} & NS & \multicolumn{1}{c|}{BRC} & \multicolumn{1}{c|}{BOME} & \multicolumn{1}{c|}{F$^2$SA} & \multicolumn{1}{c|}{BAMM} & MEHA  \\ \hline
				Time (S) & \multicolumn{1}{c|}{14.01} & \multicolumn{1}{c|}{50.13} & 31.98 & \multicolumn{1}{c|}{15.66} & 21.10  & \multicolumn{1}{c|}{19.48}  & \multicolumn{1}{c|}{11.18} & \multicolumn{1}{c|}{16.21} & \multicolumn{1}{c|}{2.650} &   \textbf{0.166}\\ \hline
				Memory   & \multicolumn{1}{c|}{160768} & \multicolumn{1}{c|}{212480} & 160768 & \multicolumn{1}{c|}{111104} & 110592  & \multicolumn{1}{c|}{12800}  & \multicolumn{1}{c|}{14848} & \multicolumn{1}{c|}{12288} & \multicolumn{1}{c|}{14848} &   \textbf{1536}\\ \hline
			\end{tabular}
		}
	\end{table}

	\begin{table}[htb]
		\centering
		\footnotesize
		\renewcommand{\arraystretch}{1.0}
		\caption{Computational efficiency comparison of  BLO schemes under LL strong-convex case.}~\label{tab:high_dimension}	\vspace{-0.2cm}
		\setlength{\tabcolsep}{1.0mm}{
			\begin{tabular}{|cccccc|cccccc|}
				\hline
				\multicolumn{6}{|c|}{Convergence Time (Dimension=10000)}                                                                                                                            & \multicolumn{6}{c|}{Convergence Time (Dimension=30000)}                                                                                                                            \\ \hline
				\multicolumn{1}{|c|}{}   & \multicolumn{1}{c|}{RHG} & \multicolumn{1}{c|}{CG} & \multicolumn{1}{c|}{NS} & \multicolumn{1}{c|}{BAMM} & MEHA & \multicolumn{1}{c|}{}   & \multicolumn{1}{c|}{RHG} & \multicolumn{1}{c|}{CG} & \multicolumn{1}{c|}{NS} & \multicolumn{1}{c|}{BAMM} & MEHA \\ \hline
				\multicolumn{1}{|c|}{Time (S)}      & \multicolumn{1}{c|}{70.59}    & \multicolumn{1}{c|}{36.97}   & \multicolumn{1}{c|}{40.53}   & \multicolumn{1}{c|}{5.890}     & 1.973   & \multicolumn{1}{c|}{Time (S)}  & \multicolumn{1}{c|}{110.26}    & \multicolumn{1}{c|}{64.21}   & \multicolumn{1}{c|}{65.73}   & \multicolumn{1}{c|}{5.949}     &  2.444      \\ \hline
				\multicolumn{1}{|c|}{Promotion} & \multicolumn{1}{c|}{$\times$35.7}    & \multicolumn{1}{c|}{$\times$18.7}   & \multicolumn{1}{c|}{$\times$20.5}   & \multicolumn{1}{c|}{$\times$3.0}     &   -   & \multicolumn{1}{c|}{Promotion} & \multicolumn{1}{c|}{$\times$31.2}    & \multicolumn{1}{c|}{$\times$26.3}   & \multicolumn{1}{c|}{$\times$26.9}   & \multicolumn{1}{c|}{$\times$2.43}     &  -     \\ \hline
			\end{tabular}
		}

	\end{table}

\begin{figure*}[h!]
		\footnotesize
		\centering
		\setlength{\tabcolsep}{12pt}

		\begin{tabular}{cccc}
			\includegraphics[width=0.24\textwidth]{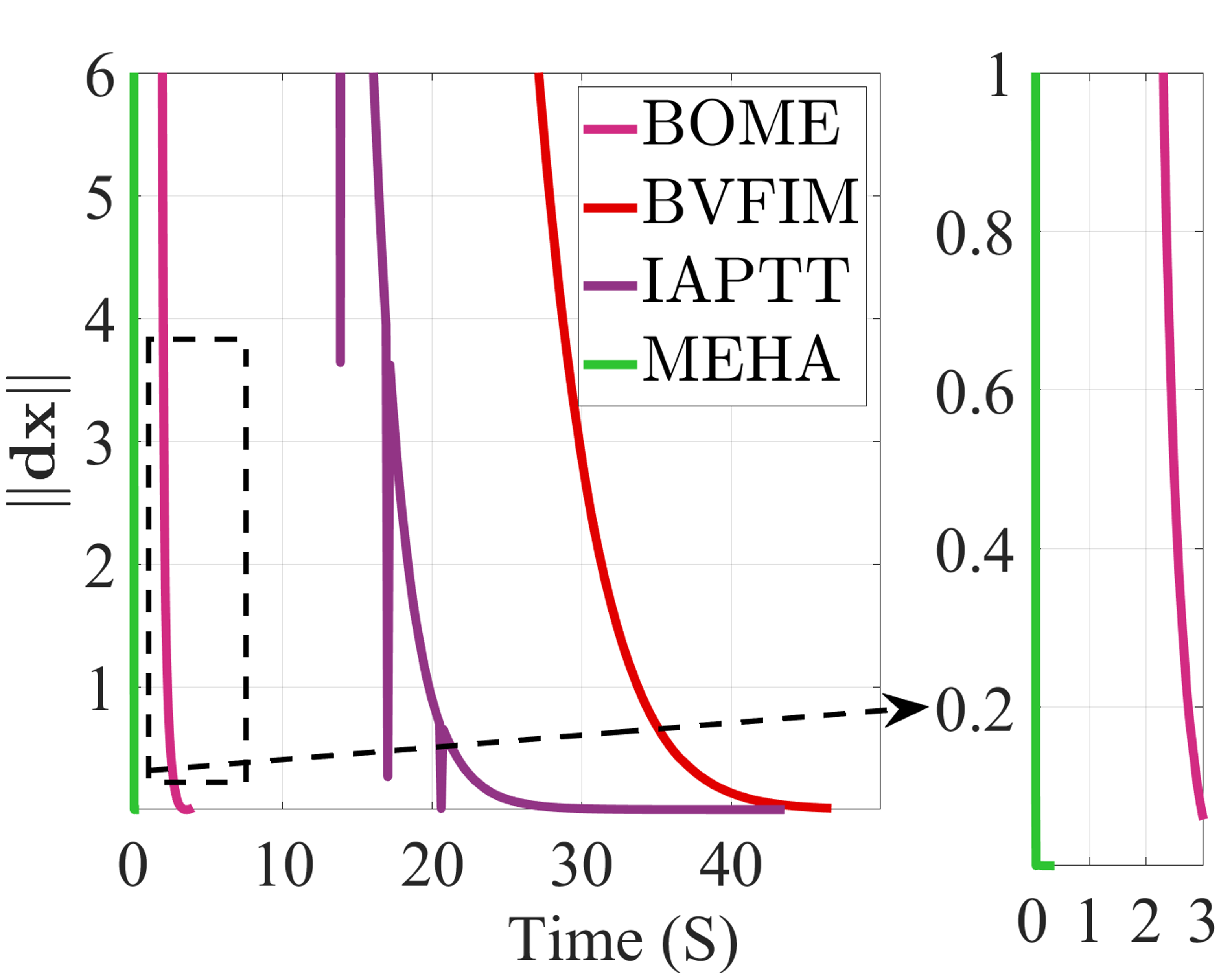}
			&  \includegraphics[width=0.24\textwidth]{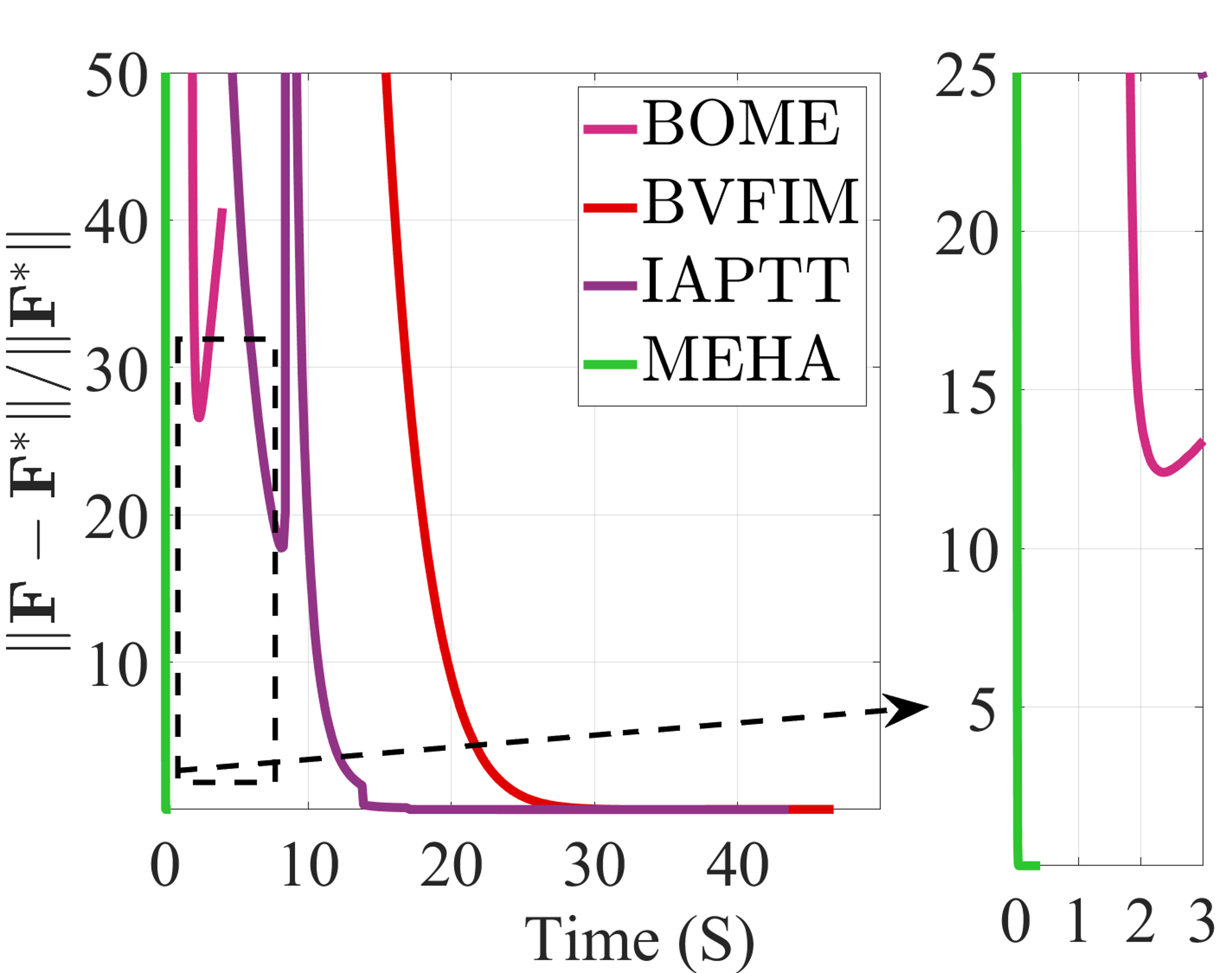}
			&\includegraphics[width=0.24\textwidth]{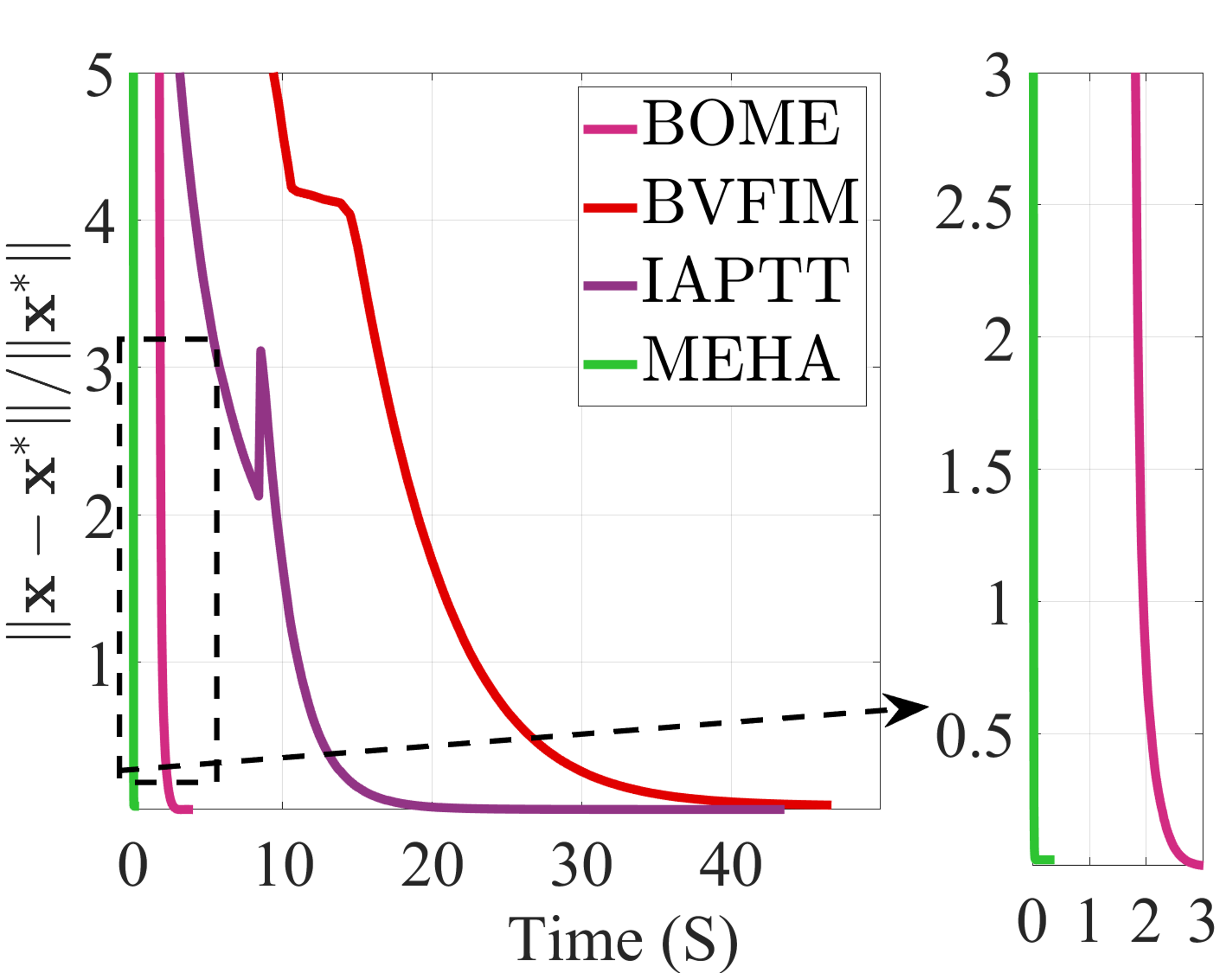}
			\\	
		\end{tabular}
		\vspace{-0.2cm}
		\caption{ Visualizing the convergence behaviours of BOME, BVFIM, IAPTT and MEHA under LL non-convex case with one dimension, using the metrics of  descent direction $\|\mathbf{d}x\|$, UL objective $F$ and reconstruction error with $x$.}\label{fig:non_convex_m}
	\end{figure*}

	\textbf{LL Non-convex Case.} 
	We depict the convergence curves utilizing advanced Bi-Level Optimization (BLO) methods in the non-convex scenario with diverse metrics in Figure~\ref{fig:non_convex_m}. 
	
	\textbf{Computational Efficiency Under Large Scale.} 
	We evaluate the computational efficiency to by increasing the dimension of $x$ and $y$ as $10^4$ and $3\times 10^4$ under  LL strong convex case. Table~\ref{tab:high_dimension} illustrates the superiority of our method to handle BLO problems with large dimension. 
	
	\begin{figure}[htb]
		\footnotesize
		\centering
		\setlength{\tabcolsep}{1pt}

		\begin{tabular}{cccc}
			\includegraphics[width=0.24\textwidth]{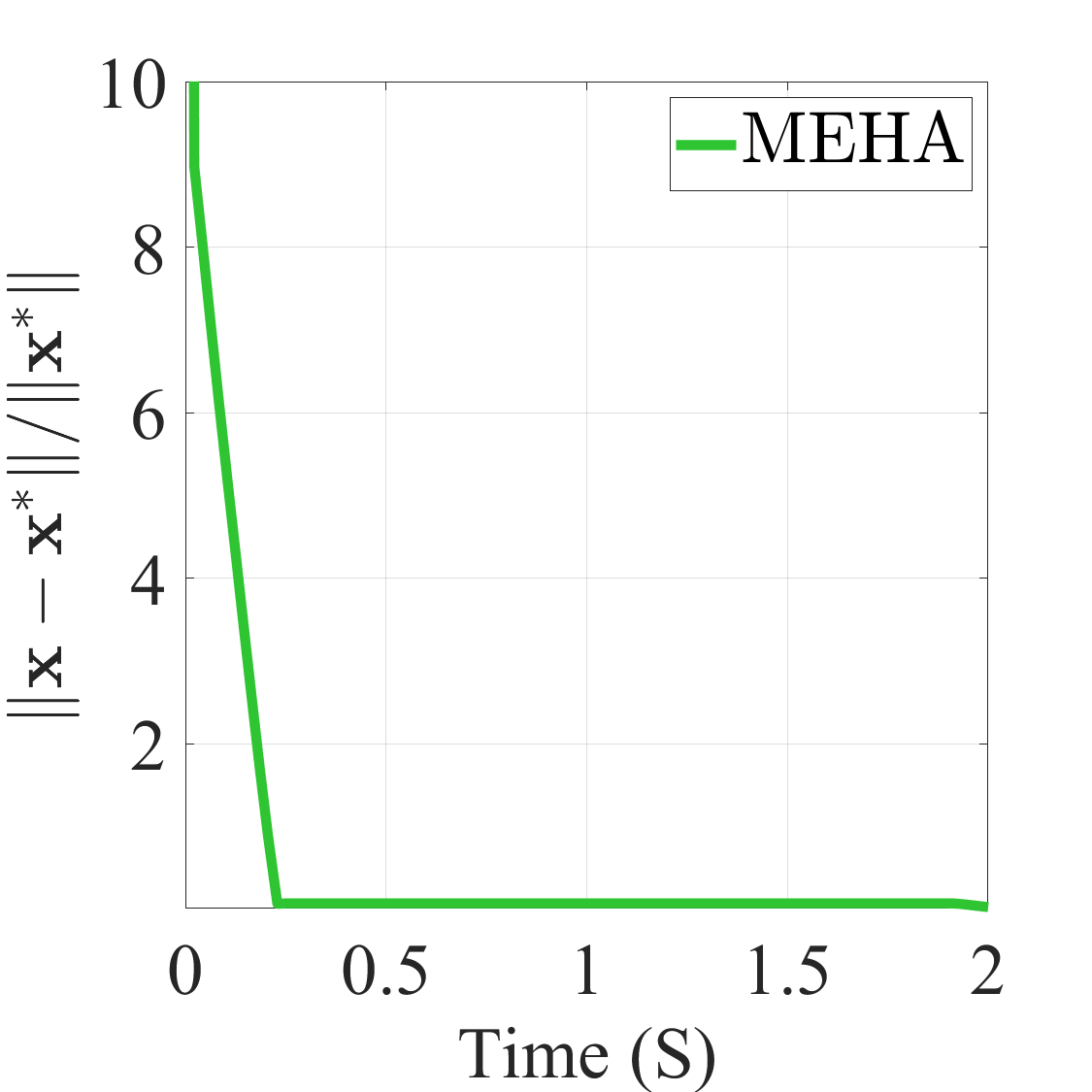}
			&\includegraphics[width=0.24\textwidth]{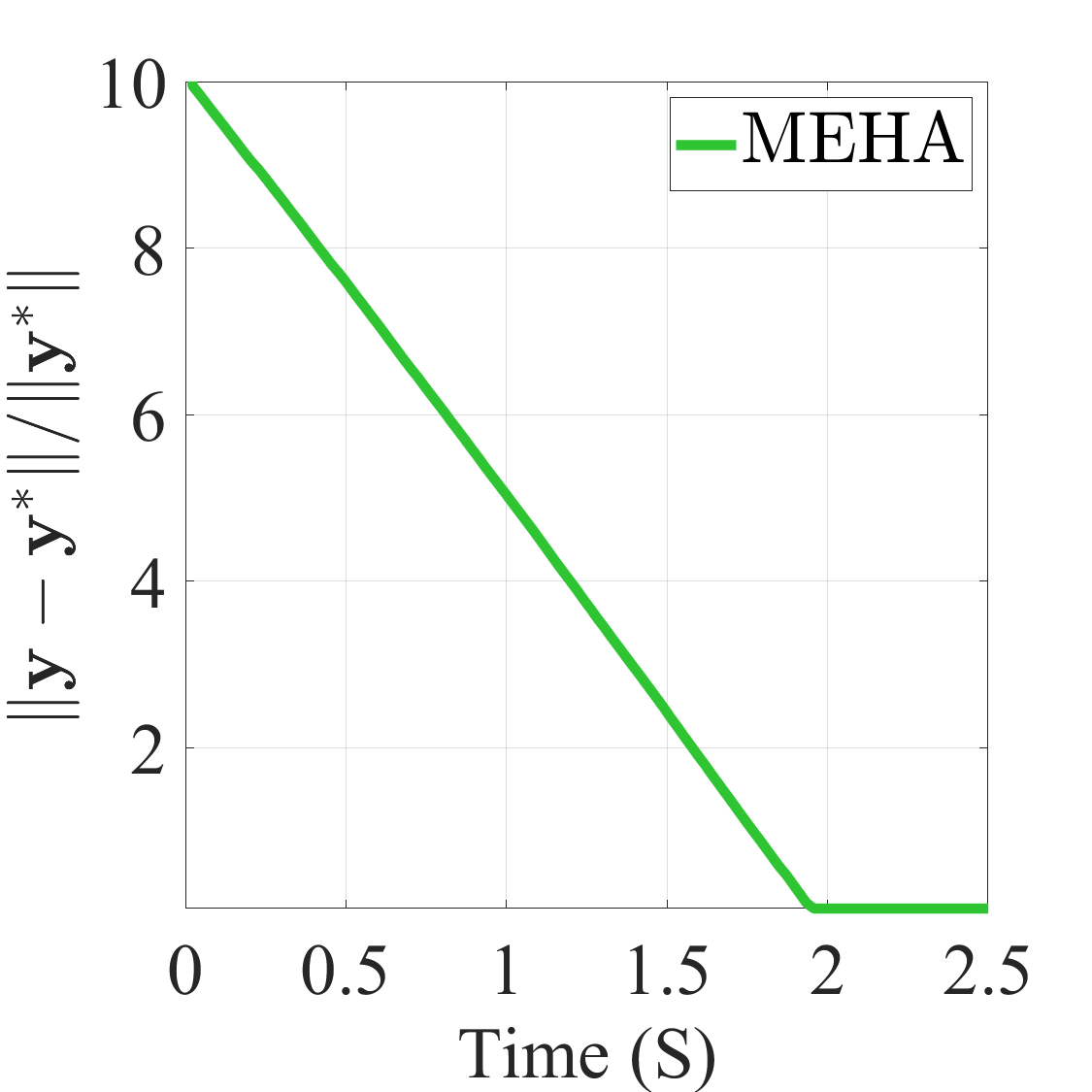}
			&	\includegraphics[width=0.24\textwidth]{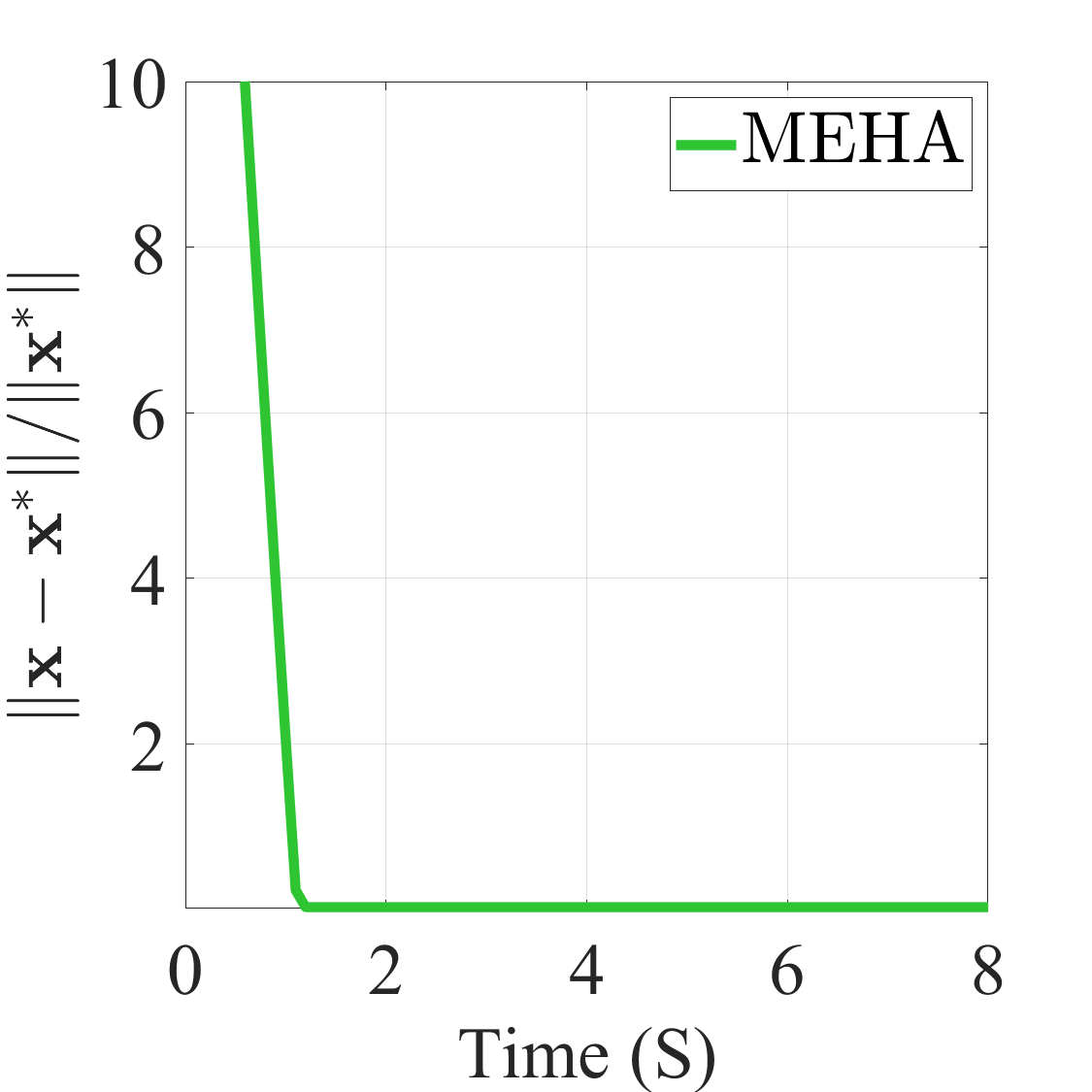}
			
			&\includegraphics[width=0.24\textwidth]{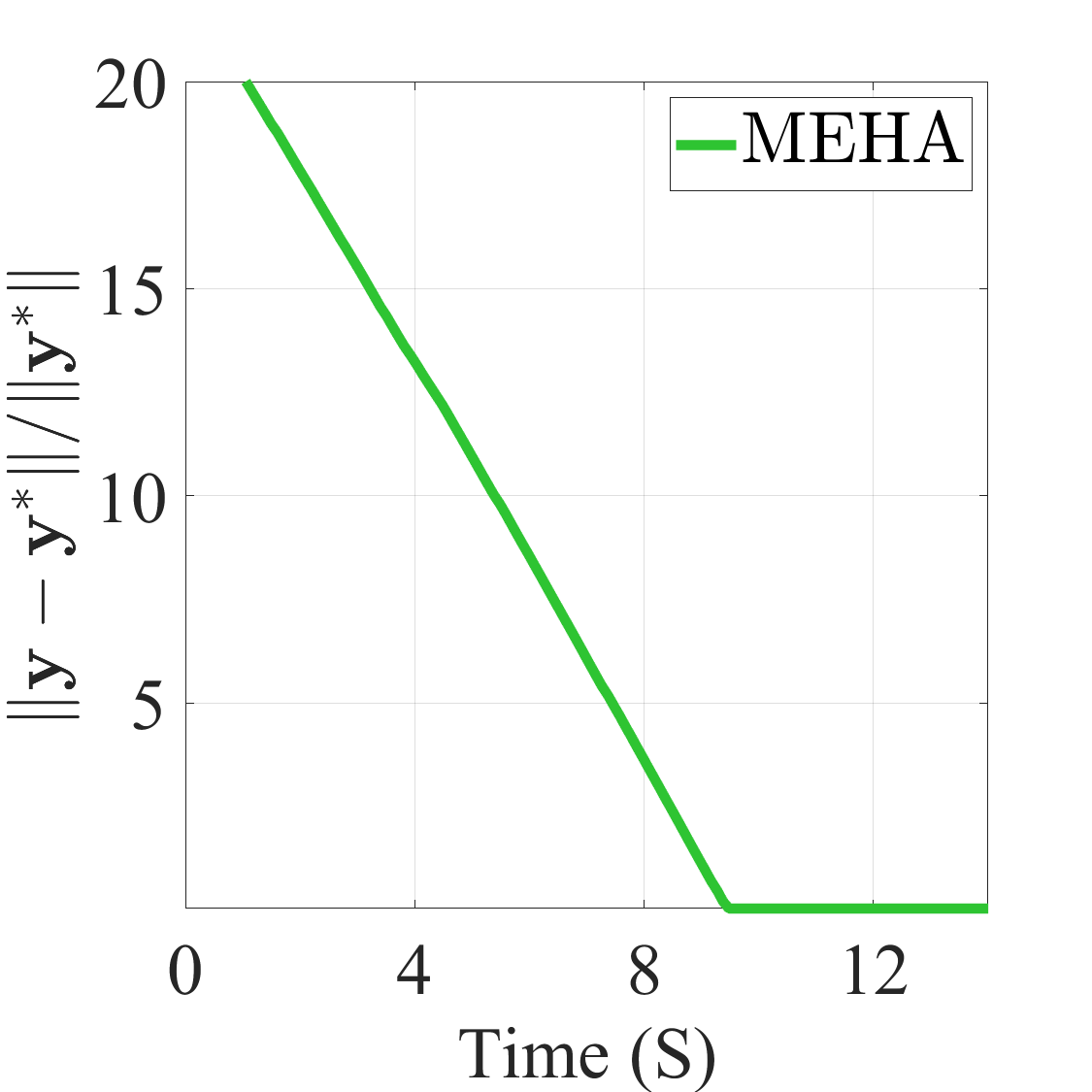}
			\\	
			\multicolumn{2}{c}{$n$= 100}  & 	\multicolumn{2}{c}{$n$= 500} \\
			\includegraphics[width=0.24\textwidth]{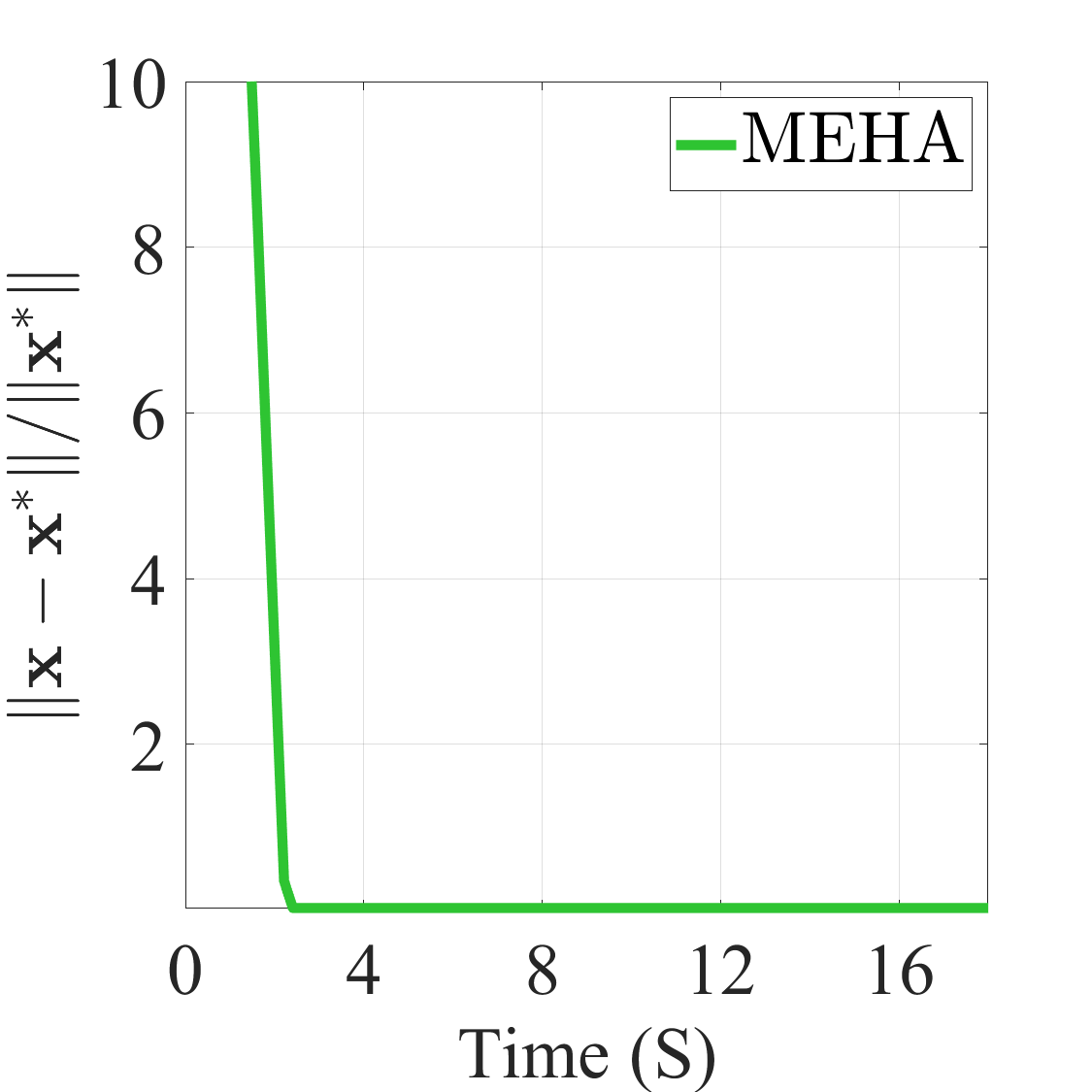}
			&  \includegraphics[width=0.24\textwidth]{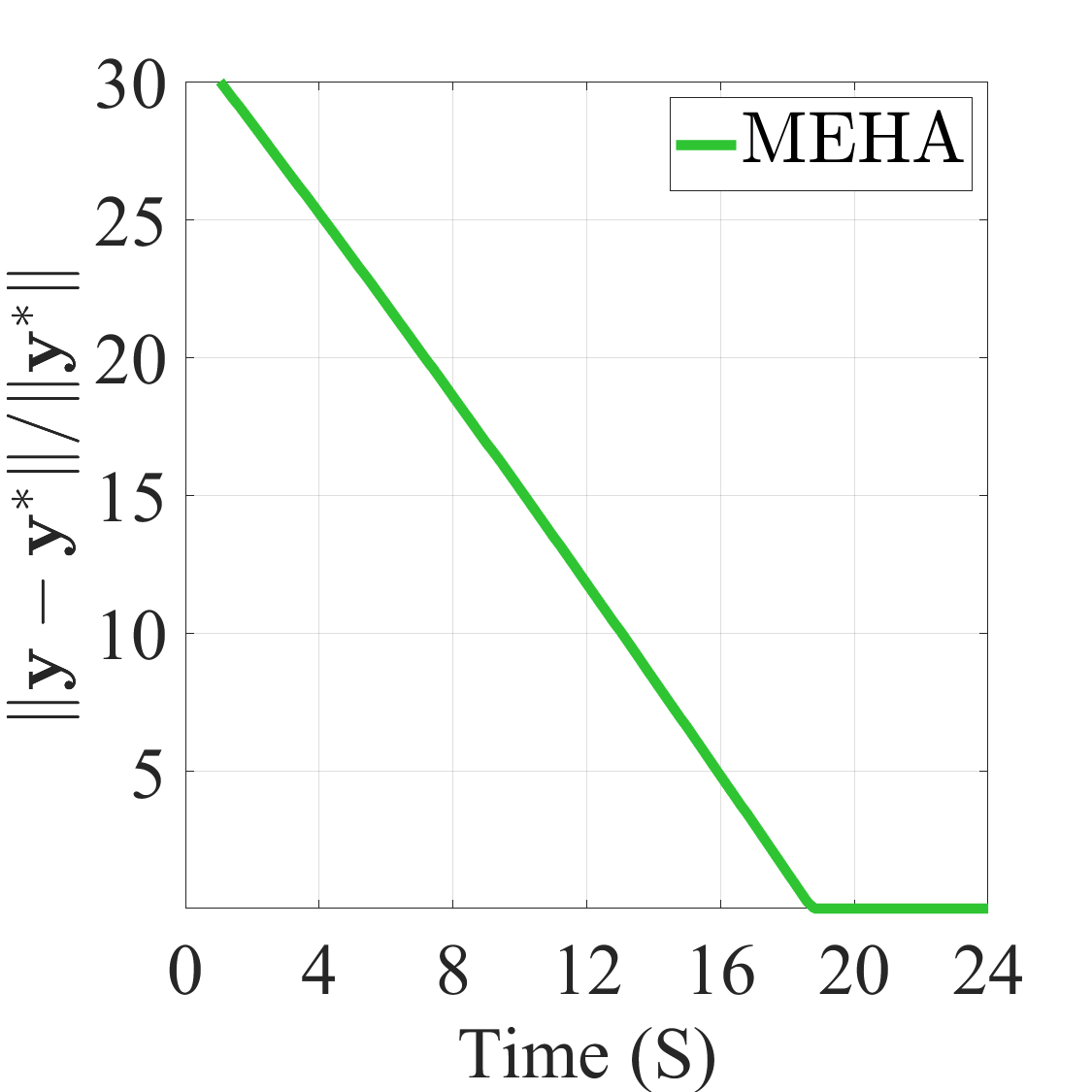}
			&\includegraphics[width=0.24\textwidth]{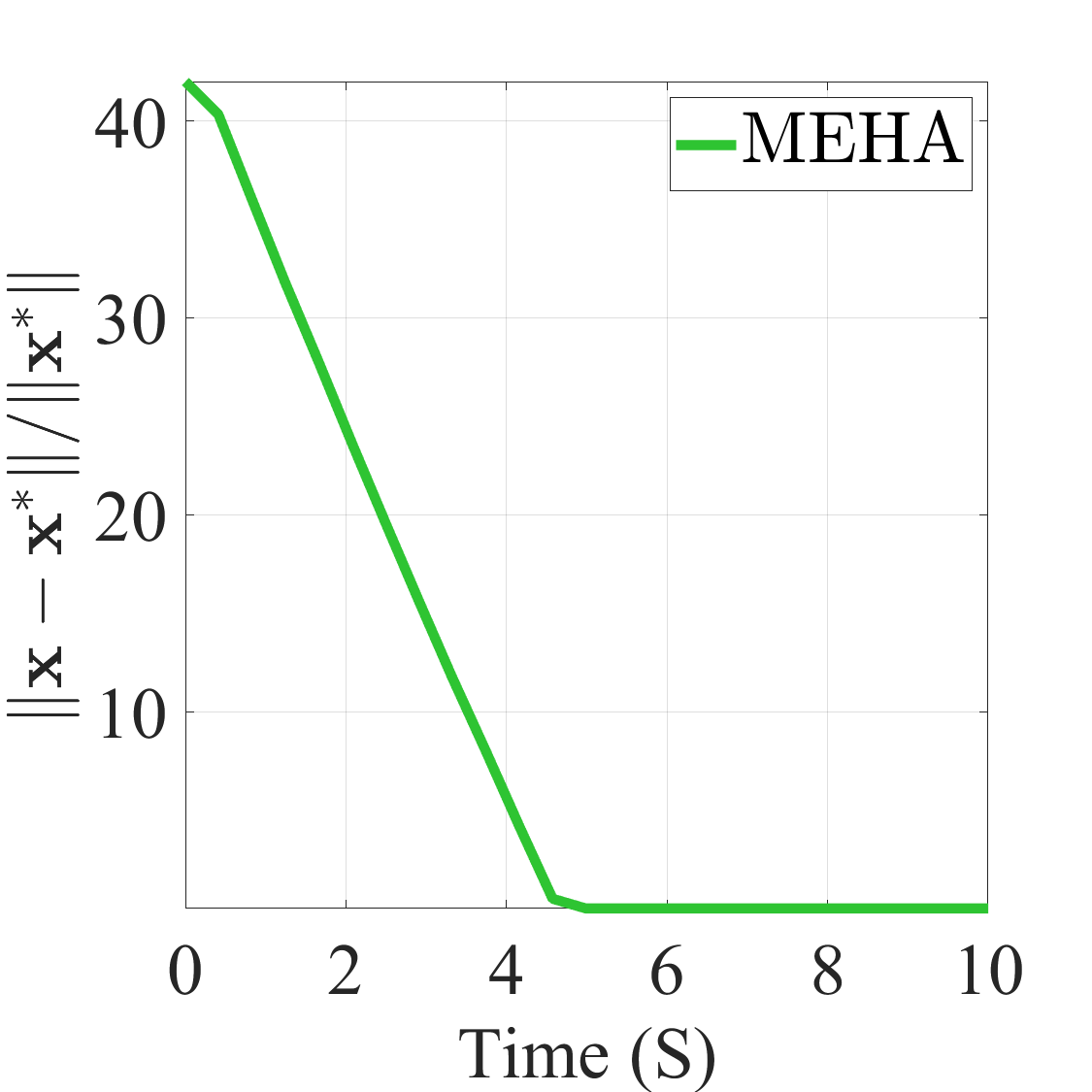}
			&\includegraphics[width=0.24\textwidth]{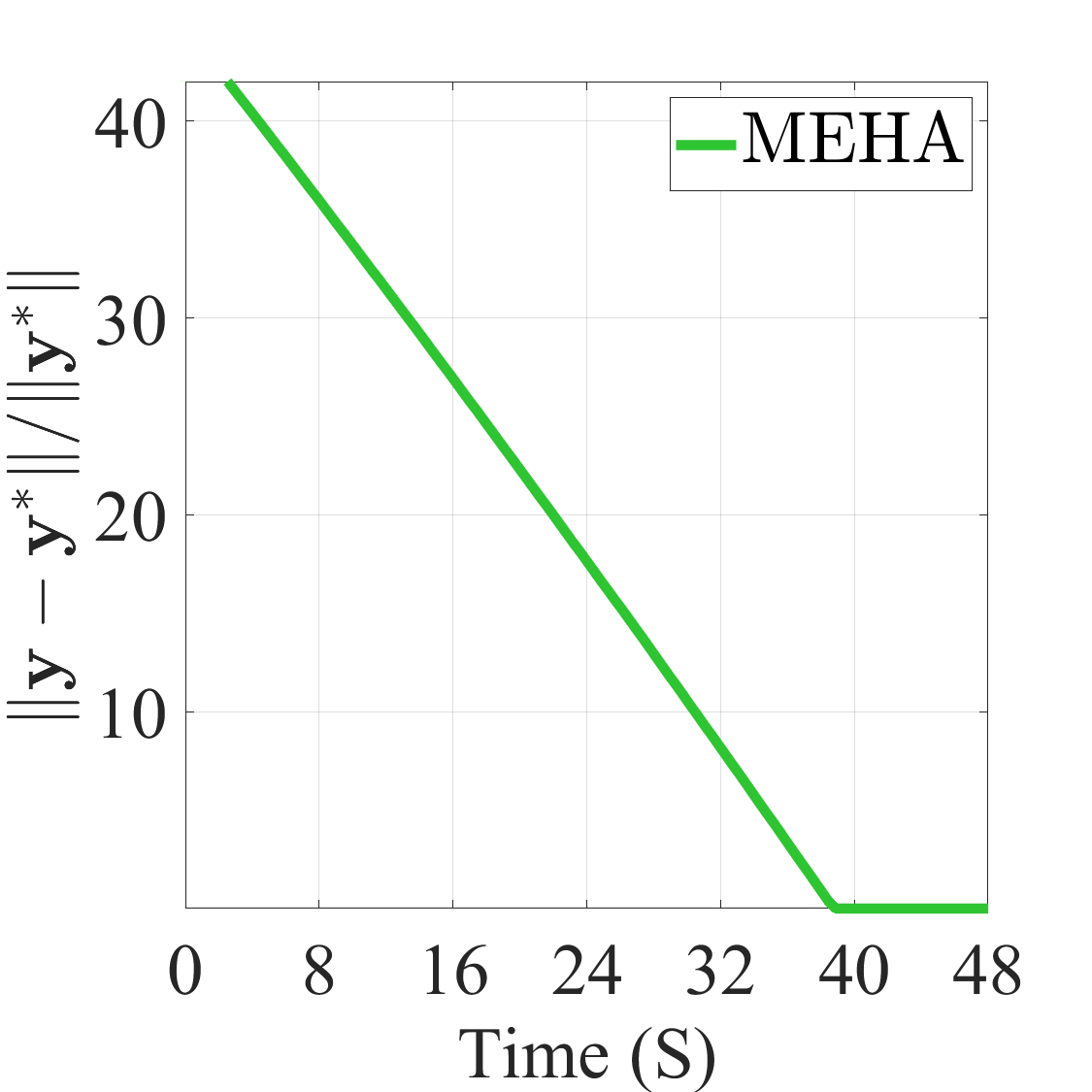}
			\\	
			\multicolumn{2}{c}{$n$= 1000}  & 	\multicolumn{2}{c}{$n$=2000} 
			
		\end{tabular}
		\vspace{-0.2cm}
		\caption{ Illustrating the convergence curves  by the criterion of $\|x-x^{*}\|/\|x^{*}\|$ and $\|y-y^{*}\|/\|y^{*}\|$ under the LL non-smooth case with different dimensions.}\label{fig:non_smooth}
	\end{figure}
	\textbf{LL Non-Smooth Case.}
	We present the convergence curves of our method, MEHA, across different dimensions in Figure~\ref{fig:non_smooth}. The results indicate that our method effectively identifies the  optimal hyper-parameter $x^*$  and the optimal solution $\{x^{*},y^{*}\}$ under diverse high dimensions.
	
	\textbf{Sensitivity of Parameters.} In Table~\ref{tab:sensitivitymore}, we report additional parameter combinations that showcase the adaptability of our scheme. These results emphasize that various combinations of parameters can achieve convergence in the algorithms.
	
	\begin{table}[h!]
		\centering
		\footnotesize
		\renewcommand{\arraystretch}{1.0}
		\caption{Analyzing the sensitivity of parameters by altering one parameter and keeping the others constant.}~\label{tab:sensitivitymore}	
		\setlength{\tabcolsep}{3mm}{
			\begin{tabular}{|c|c|c|c|c|c|c|c|c|}
				\hline
				Strategy                                  & $\alpha$ & $\beta$ & $\eta$ & $\gamma$ & $\underline{c}$ &$p$ & Steps & Time (s) \\ \hline
				Original                                  & 0.1        & 1e-5 & 0.1      & 10       & 2   &0.49            & 97    & 22.83    \\ \hline
				\multirow{4}{*}{ $\alpha$}      & 0.025       & 1e-5 & 0.1      & 10       & 2          &0.49     & 95   & 10.30    \\ \cline{2-9} 
				& 0.2        & 1e-5 & 0.1      & 10       & 2           &0.49     & 88    & 9.65    \\ \cline{2-9}
				& 0.4        & 1e-5 & 0.1      & 100      & 2    &0.49           & 88    & 12.55     \\ \cline{2-9}
				& 0.6       & 1e-5 & 0.1      & 10       & 2    &0.49          & 86   & 9.53   \\  \hline
				\multirow{3}{*}{$\beta$}       & 0.1        & 2e-6   & 0.1      & 10       & 2          &0.49     & 106    & 11.53  \\ \cline{2-9} 
				
				& 0.1        & 8e-6   & 0.1      & 100      & 2    &0.49           & 99   & 10.74    \\ \cline{2-9}
				& 0.1        & 3.5e-5   & 0.1      & 10       & 2    &0.49          & 82   & 11.71    \\ \hline
				\multirow{4}{*}{ $\eta$}        & 0.1        & 1e-5 & 0.02      & 10       & 2            &0.49   & 142    & 15.43    \\ \cline{2-9} 
				& 0.1        & 1e-5 & 0.08     & 10       & 2         &0.49      & 100   & 10.91    \\ \cline{2-9}
				& 0.1        & 1e-5 & 0.4     & 100      & 2    &0.49           & 89    & 9.59   \\ \cline{2-9}
				& 0.1        & 1e-5 & 0.6      & 10       & 2    &0.49          & 88   & 9.59   \\\hline
				\multirow{3}{*}{ $\gamma$}        & 0.1        & 1e-5 & 0.1      & 4        & 2         &0.49      & 90    & 9.84   \\ \cline{2-9} 
				
				& 0.1        & 1e-5 & 0.1      & 50       & 2    &0.49          & 89   & 9.71    \\ \cline{2-9}
				& 0.1        & 1e-5 & 0.1      & 200       & 2    &0.49          & 90   & 9.91    \\ \hline
				\multirow{3}{*}{ $\underline{c}$} & 0.1        & 1e-5 & 0.1      & 10     & 5   &0.49           & 94   & 10.20    \\ \cline{2-9} 
				& 0.1        & 1e-5 & 0.1      & 10       & 20        &0.49      & 97   & 10.65    \\ \cline{2-9}
				
				& 0.1        & 1e-5 & 0.1      & 10       & 80    &0.49          & 38   & 10.74    \\  
				\hline
				
				\multirow{3}{*}{ $p$}   
				& 0.1        & 1e-5 & 0.1      & 10       & 2    &0.1          & 97   & 10.59    \\ \cline{2-9}
				& 0.1        & 1e-5 & 0.1      & 10       & 2    &0.2          & 96   & 10.43    \\ \cline{2-9}
				& 0.1        & 1e-5 & 0.1      & 10       & 2    &0.4          & 90   & 9.69    \\ \hline
			\end{tabular}
		}

	\end{table}
	
	{\bf Neural Architecture Search.} 	Additionally, we have also conducted a NAS experiment to test the performance of our proposed method on architectures that use the smooth activation function (Swish). For comparison, the compared version was searched within the original search space that features non-smooth activation functions. This comparison is detailed in Table~\ref{tab:nas2} that includes two variants of our method: MEHA (non-smooth) and MEHA (smooth). 
	
	\begin{table*}[!h]
		\centering
		\footnotesize
		\renewcommand{\arraystretch}{1.0}
		\caption{Comparing Top-1 accuracy in searching, inference, and final test stages with diverse activated functions.}~\label{tab:nas2}	\vspace{-0.2cm}
		\setlength{\tabcolsep}{6mm}{
			\begin{tabular}{|c|cc|cc|c|c|c|}
				\hline
				\multirow{2}{*}{Methods} & \multicolumn{2}{c|}{Searching} & \multicolumn{2}{c|}{Inference} & \multirow{2}{*}{Test} & \multirow{2}{*}{Params (M)} \\ \cline{2-5}
				& \multicolumn{1}{c|}{Train}    &  Valid    & \multicolumn{1}{c|}{Train}    &  Valid    &                       &                             \\ \hline
				
				MEHA (non-smooth) & \multicolumn{1}{c|}{99.060}    &  {99.764}   & \multicolumn{1}{c|}{99.419}    & {96.150}    &    {96.070}                   &    1.524                         \\ \hline
				
				MEHA (smooth) & \multicolumn{1}{c|}{99.292}    &  {99.786}   & \multicolumn{1}{c|}{99.360}    & {95.870}    &    {95.870}                   &    1.479                        \\ \hline
			\end{tabular}
		}
		
	\end{table*}
	
	\subsection{Experimental  Details}~\label{sec:detail}
	
	First, by using the results in Section \ref{A9-assump}, it can be verified that all models employed in our experimental studies meet the Assumptions assumed.

	Second, we conducted the experiments on a PC with Intel i7-8700 CPU (3.2 GHz), 32GB RAM and NVIDIA RTX 1070 GPU. We utilized the PyTorch framework on the 64-bit Windows system.
	
	\textbf{Synthetic Numerical Examples.} The hyper-parameter settings for diverse numerical experiments are summarized in Table~\ref{tab:sne}.  
	We utilize the SGD optimizer to update the UL variable $x$. We uniformly utilize the $\|x-x^{*}\|/\|x^{*}\|\leq 10^{-3}$ criterion for LL  merely convex case. We utilize the $\|dx\|\leq  1e^{-8}$ for one-dimension LL non-convex case $\|dx\|\leq  1e^{-3}$ for other high dimension cases. The maximum of iteration steps is 800. The learning steps $\alpha$, $\beta$ and $\eta$ are fixed. Specifically, inverse power learning rate annealing strategy to dynamically adjust the learning rate  ($\alpha$ and $\beta$) for LL merely-convex case.
	
	As for the LL non-smooth case, We leverage the method~\cite{feng2018gradient} to generate the synthesized data for the group lasso hyper-parameter selection case, including 100 training, validation, and testing pairs, respectively. We set the response $b_{i}$ as $b_{i}=v^{T}a_{i}+\sigma\epsilon_{i}$, where $v = [v^{(1)},v^{(2)},v^{(3)}]$, $v^{(i)} = ({1,1,1,\cdots,1},0,0,0,\cdots,0)$. The number of elements with value 1 is 50. 
	The parameters $\sigma$ are sampled from a standard normal distribution, and $\epsilon$ is selected to achieve an SNR of 2.   The number of groups is set as 30 for 600 dimensions and 300 for other scenarios.  $\|x^{k}-x^{k-1}\|/\|x^{k}\|\leq 0.2$ is as the stop criterion. As for the experimental settings of other compared methods (including grid search, random search, TPE, IGJO and VF-iDCA), we follow the effective practice~\cite{gao2022value}. 
	\begin{table}[htb]
		\centering
		\footnotesize
		\renewcommand{\arraystretch}{1.1}
		\caption{Values for hyper-parameters of synthetic numerical experiments.}~\label{tab:sne} \vspace{-0.2cm}
		\setlength{\tabcolsep}{4mm}{
			\begin{tabular}{|c|c|c|c|c|c|c|}
				\hline
				Category & $\alpha_0$ & $\beta_0$ & $\eta_0$ & $\gamma$ & $\underline{c}$ & $p$\\ \hline
				LL strong-convex   &      1.5    & 0.8   &    0.8        & 10                &  33.3 &  0.49   \\ \hline
				LL merely-convex  &      0.012    & 0.1   &    0.009        & 5                &  0.167    &   0.49   \\ \hline
				LL non-convex   &      5e$-4$    & 5e$-4$    &    0.001        & 200                &  0.02    &  0.49   \\ \hline
				LL non-smooth (Eq.~\eqref{eq:LNS}) &      0.1    &   1e$-5$ &    0.1        & 10                &  2   & 0.49 \\ \hline
				LL non-smooth (Eq.~\eqref{eq:lass}) &      0.01    &   0.05&    0.05      & 100                &  20   & 0.48 \\ \hline
			\end{tabular}
		}
		
	\end{table}

	\textbf{Few-shot Learning.} As for this task, the upper-level variables $x$ represent the shared weights for feature extraction. $y:=\{y^j\}$ denotes the task-specific parameters. Leveraging the cross-entropy loss as the objective $\mathcal{L}$, we provide the bi-level formulation as:
	\begin{equation}
		\min_{x, y} \sum_{j}\mathcal{L}(x,y^{j};\mathcal{D}_\mathtt{val}^{j}) \quad \text { s.t. } \quad y\in \arg\min_{\tilde{y}} \sum_{j} \mathcal{L}(x,\tilde{y}^{j};\mathcal{D}_\mathtt{train}^{j}).
	\end{equation}
	Following with the practice~\cite{liu2023augmenting}, we utilize four layers of convolution blocks (ConvNet-4) to construct the backbone (\textit{i.e.,} $x$), which is widely utilized for few-shot learning tasks. The task-specific classifier $y$ is composited by fully-connection layers with softmax operation. 
	Adam and SGD optimizers are utilized to update $x$ and $y$ for all algorithms fairly.
	The concrete hyper-parameters of Alg.1 and other shared hyper-parameters  are summarized in Table~\ref{tab:meta1}.  We utilize the inverse power learning rate annealing strategy to dynamically adjust the learning rate  ($\alpha$ and $\beta$). $\eta$ and $\gamma$ are fixed.
	\begin{table}[htb]
		\centering
		\footnotesize
		\renewcommand{\arraystretch}{1.1}
		\caption{Values for hyper-parameters of few-shot learning.}~\label{tab:meta1} \vspace{-0.2cm}
		\setlength{\tabcolsep}{3mm}{
			\begin{tabular}{|c|c|c|c|c|c|c|c|c|}
				\hline
				Parameter & Meta batch size &  Hidden size & $\alpha_0$ & $\beta_0$ & $\eta_0$ & $\gamma$ & $\underline{c}$ &$p$\\ \hline
				Value  & 16 & 32 &     0.08    & 0.05   &   0.001       & 100                &  0.067      & 0.08   \\ \hline
				
			\end{tabular}
		}
		
	\end{table}
	
	\textbf{Data Hyper-Cleaning.}  The mathematical formulation can be written as:
	\begin{equation}
		\min_{x, y} \sum_{\mathbf{u}_{i},\mathbf{v}_{i}\in\mathcal{D}_\mathtt{val}}\mathcal{L}(y;\mathbf{u}_{i},\mathbf{v}_{i}) \quad \text { s.t. } \quad y\in  \arg\min_{\tilde{y}}\sum_{\mathbf{u}_{i},\mathbf{v}_{i}\in\mathcal{D}_\mathtt{train}}[\sigma(x)]_{i}\mathcal{L}(\tilde{y};\mathbf{u}_{i},\mathbf{v}_{i}),
	\end{equation}
	where the upper-level variable $x$ is a vector with the same dimension of the number of corrupted examples. $y$ denotes the target classification model. $\sigma(x)$ is a sigmoid function. $\{\mathbf{u},\mathbf{v}\}$ are the data pairs. In detail, we only utilize two-layer of fully-connection to define $y$ for the noncovex case. Two datasets FashionMNIST and MNIST are utilized to conduct the experiments. We randomly split these datasets to composite the training, validation and testing subsets with 5000, 5000, 10000 examples, respectively. Half of data in the training dataset is tampered.  The concrete hyper-parameters of Alg.1 are summarized in Table~\ref{tab:dhc}. SGD optimizer is utilized to update the UL variable $x$ fairly. We utilize the inverse power learning rate annealing strategy to dynamically adjust the learning rate ($\alpha$ and $\beta$). $\eta$ and $\gamma$ are fixed.
	
	\begin{table}[htb]
		\centering
		\footnotesize
		\renewcommand{\arraystretch}{1.1}
		\caption{Values for hyper-parameters of data hyper-cleaning.}~\label{tab:dhc} \vspace{-0.2cm}
		\setlength{\tabcolsep}{4mm}{
			\begin{tabular}{|c|c|c|c|c|c|c|}
				\hline
				Parameter & $\alpha_0$ & $\beta_0$ & $\eta_0$ & $\gamma$ & $\underline{c}$ & $p$ \\ \hline
				Value  &       0.1    & 0.15   &    0.001        & 100                &  0.2 &0.25  \\ \hline
				
			\end{tabular}
		}
		
	\end{table}
	\textbf{Neural Architecture Search.}
	The bi-level formulation of  neural architecture search is 
	\begin{equation}
		\min_{x, y} \mathcal{L}_\mathtt{val}(y,x;\mathcal{D}_\mathtt{val}) \quad \text { s.t. } \quad y\in  \arg\min_{\tilde{y}}\mathcal{L}_\mathtt{train}({\tilde{y}},x;\mathcal{D}_\mathtt{train}),
	\end{equation}
	where the architecture parameters are denoted as the upper-level variable $x$ and the lower-level variable $y$ represents the network weights. $\mathcal{L}_\mathtt{val}$ and $\mathcal{L}_\mathtt{train}$ are the losses on validation and training datasets. The definition of search space, cells, and experimental hyper-parameters settings are following with the literature~\cite{liu2018darts}. We leveraged the Cifar-10 dataset to perform the experiments of image classification. As for the super-network, we conducted the search procedure with three layers of cells for 50 epochs.
	The network for training is increased with 8 layers and trained from scratch with 600 epochs.  The concrete hyper-parameters of Alg.1 are summarized in Table~\ref{tab:nas_p}. We utilized the cosine decreasing learning rate annealing strategy to dynamically adjust the learning rate ($\alpha$, $\beta$ and $\eta$). $\eta$ and $\gamma$ are fixed.  
	
	 Our experiment for NAS was designed to evaluate the performance of our proposed algorithm within practical bilevel optimization problems, particularly those involving a nonconvex lower-level objective. To ensure a consistent and fair comparison against existing methods, we adopted the subgradient descent technique for handling the nonsmoothness in neural networks, aligning our approach with that of our competitors~\cite{liu2018darts}. Specifically, our method utilizes the same search space defined in the referenced works, which incorporates the nonsmooth ReLU activation function in certain operations. For the gradient computation of the ReLU activation function, we define a unit step function Heaviside(x) for ReLU. It takes the value of $x>0$ when 
	and $x\leq 0$ when 
	to approximate the gradient.
	\begin{table}[htb]
		\centering
		\footnotesize
		\renewcommand{\arraystretch}{1.1}
		\caption{Values for hyper-parameters of neural architecture search.}~\label{tab:nas_p} \vspace{-0.2cm}
		\setlength{\tabcolsep}{4mm}{
			\begin{tabular}{|c|c|c|c|c|c|c|}
				\hline
				Parameter & $\alpha_0$ & $\beta_0$ & $\eta_0$ & $\gamma$ & $\underline{c}$ & $p$ \\ \hline
				Value  &       8e-5     &  0.025   &    0.025         & 200                &  2    &  0.49 \\ \hline
				
			\end{tabular}
		}
	\end{table}

	\subsection{Equivalence of Moreau Envelope Based Reformulation}\label{equiv-a}
	
	The following theorem establishes the equivalence between the Moreau Envelope based reformulation problem (\ref{wVP}) and the relaxed bilevel optimization problem (\ref{problem_relax}). The proof is inspired by the one of Theorem 1 in \cite{gao2023moreau}.
	For the convenience of the reader, we restate problems (\ref{wVP}) and (\ref{problem_relax}) as follows:
	\begin{equation}\tag{\ref{wVP}}
		\begin{aligned}
			\min_{(x,y)\in X \times Y}  & F(x,y) 
			\quad
			\mathrm{s.t.} & \varphi(x,y)-v_\gamma (x,y)\leq 0,
		\end{aligned}
	\end{equation} 
	where $v_{\gamma}(x,y)  
	:=\inf_{\theta \in Y} \left \{ \varphi(x, \theta) + \frac{1}{2\gamma}\|\theta-y\|^2 \right \}$, $\varphi(x,y)=f(x,y)+g(x,y)$, and 
	\begin{equation}\tag{\ref{problem_relax}}
		\begin{aligned}
			\min_{x \in X, y \in Y}  & F(x,y)  
			\quad
			\mathrm{s.t.}  &   
			0 \in \nabla_{y}f(x,y) + \partial_y g(x,y) + \mathcal{N}_Y(y).
		\end{aligned}
	\end{equation}

	\begin{theorem}\label{thm_reformulate-a}
		Suppose that $\varphi(x, \cdot)$ is  $\rho_{\varphi_2}$-weakly convex on $Y$ for all $x$, i.e., $\varphi(x, \cdot)+\frac{\rho_{\varphi_2}}{2}\|\cdot\|^2$ is convex on $Y$ for all $x$.
		Then for $\gamma \in (0, 1/\rho_{\varphi_2})$, the Moreau Envelope based reformulation problem (\ref{wVP}) is equivalent to the relaxed BLO problem (\ref{problem_relax}).
	\end{theorem}
	\begin{proof}
		First, given any feasible point $(x,y)$ of problem (\ref{wVP}), it necessarily belongs to $X \times Y$ and satisfies
		\[
		\varphi(x,y) \le v_\gamma(x,y) := \inf_{\theta \in Y} \left\{ \varphi(x, \theta) + \frac{1}{2\gamma}\|\theta - y\|^2 \right\} \le \varphi(x,y).
		\]
		From which, it follows that $\varphi(x,y) = v_\gamma(x,y)$ and thus $y \in \mathrm{argmin}_{\theta \in Y}\left\{ \varphi(x, \theta) + \frac{1}{2\gamma}\|\theta - y\|^2 \right\} $.
		This leads to
		\[
		0 \in \nabla_y  f(x,y) + \partial_y g(x,y) + \mathcal{N}_Y(y),
		\]
		implying that $(x,y)$ is feasible for problem (\ref{problem_relax}).
		
		Conversely, consider that $(x,y)$ is an feasible point of problem (\ref{problem_relax}). This implies that $(x,y) \in X \times Y$, $0 \in \nabla_y  f(x,y) + \partial_y g(x,y) + \mathcal{N}_Y(y)$. 
		Given that $\varphi(x, \cdot) : \mathbb{R}^m \rightarrow \mathbb{R}$ is $\rho_{\varphi_2}$-weakly convex on $Y$, then when $\gamma \in (0, 1/\rho_{\varphi_2})$, the function $ \varphi(x,\cdot) + \frac{1}{2\gamma} \|\cdot - y\|^2 $ is convex on $Y$, making it lower regular. Clearly,  $\delta_{Y}(\cdot)$ is lower regular since $Y$ is a closed convex set.
		By leveraging  the subdifferential sum rules for two lower regular l.s.c. functions Theorem 2.19 of \cite{mordukhovich2018variational}, we arrive at
		$$  
		\partial \left( \varphi(x,\cdot) + \frac{1}{2\gamma} \|\cdot - y\|^2 + \delta_{Y}(\cdot)\right)
		= \nabla_y f(x,\cdot) + \partial_y g(x, \cdot) + (\cdot - y)/\gamma + \mathcal{N}_Y(\cdot).
		$$ 
		With the right-hand set-valued mapping at $y$ containing $0$, we can deduce that $$0\in \partial_\theta\left.\left( \varphi(x, \theta) + \frac{1}{2\gamma} \| \theta - y\|^2 + \delta_{Y}(\theta)\right)\right|_{\theta=y}. $$
		Thus, invoking the first-order optimally condition for convex functions, we infer 
		$$y \in \mathrm{argmin}_{\theta \in Y}\left\{ \varphi(x, \theta) + \frac{1}{2\gamma}\|\theta - y\|^2 \right\}. $$ This implies $\varphi(x,y) = v_\gamma(x,y)$, confirming $(x,y)$ as an feasible point to problem (\ref{wVP}).
	\end{proof}

	\subsection{Convergence of Approximation Problem (\ref{problem_pen})}\label{conv_appro}
	
	For the convenience of the reader, we restate the approxiimation problem  (\ref{problem_pen}) as follows:
	\begin{equation}\tag{\ref{problem_pen}}
		\min_{(x,y)\in X \times Y} \psi_{c_k}(x,y) := F(x,y) + c_k \Big( f(x,y) + g(x,y) - v_{\gamma} (x,y)\Big).
	\end{equation}	
	We will show that, as $c_k \rightarrow \infty$, any limit point of the sequence of solutions to the approximation problem (\ref{problem_pen}), associated with varying values of $c_k$, is a solution to the Moreau envelope based reformulation (\ref{wVP}). 
	Using the same proof technique as presented in Lemma 2 of \cite{liu2020generic}, we can establish the following result.
	\begin{lemma}\label{lem_USC_v}
		If $\varphi(\cdot,y)$ is continuous on $X$ for any $y \in Y$, then $v_\gamma (x,y)$ is upper semi-continuous on $X \times \mathbb{R}^m$.
	\end{lemma}
	We are now prepared to establish the convergence of the approximation problem (\ref{problem_pen}).
	\begin{theorem}\label{Thm_app}
		Assume that $X$ and $Y$ are closed and $F$, $f$ and $g$ are all continuous on $X \times Y$. Additionally, suppose $c_k \rightarrow \infty$ and let 
		\[
		(x_k,y_k) \in \underset{(x,y)\in X \times Y}{\mathrm{argmin}} \psi_{c_k}(x,y).
		\]
		Then, for any limit point $(\bar{x}, \bar{y})$ of the sequence $\{(x_k, y_k)\}$, $(\bar{x}, \bar{y})$ is a solution to the Moreau envelope based reformulation (\ref{wVP}).
	\end{theorem}
	\begin{proof}
		Let $(\bar{x}, \bar{y})$ be any limit point of the sequence $\{(x_k, y_k)\}$ and $\{(x_j, y_j)\} \subseteq \{(x_k, y_k)\}$ be the subsequence such that $(x_j, y_j) \rightarrow (\bar{x}, \bar{y})$. Due to the closedness of both $X$ and $Y$, we have $(\bar{x}, \bar{y}) \in X \times Y$. For any $\epsilon > 0$, consider $(x_\epsilon, y_\epsilon) \in X \times Y$ as a point that satisfies $\varphi(x_\epsilon, y_\epsilon) -v_\gamma (x_\epsilon, y_\epsilon) \leq 0 $ and
		\[
		F(x_\epsilon, y_\epsilon) < \inf_{x \in X, y \in Y} \left\{ F(x,y) ~~\mathrm{s.t.} ~~ \varphi(x,y)-v_\gamma (x,y)\leq 0\right\} + \epsilon.
		\]
		Then, since $(x_k, y_k) \in \arg\min_{(x,y)\in X \times Y} \psi_{c_k}(x,y)$, we have
		\begin{equation}\label{Thm_app_eq1}
			\begin{aligned}
				&F(x_k, y_k)+ c_k \Big( f(x_k, y_k) + g(x_k, y_k) - v_{\gamma} (x_k, y_k)\Big) \\ \le ~&F(x_\epsilon, y_\epsilon) 
				<  \inf_{x \in X, y \in Y} \left\{ F(x,y) ~~\mathrm{s.t.} ~~ \varphi(x,y)-v_\gamma (x,y)\leq 0\right\} + \epsilon.
			\end{aligned}
		\end{equation}
		By taking $k = j$ and letting $j \rightarrow \infty$ in  (\ref{Thm_app_eq1}), and as $c_k \rightarrow \infty$, we obtain
		\[
		\limsup\limits_{j \rightarrow \infty} ~~ f(x_j, y_j) + g(x_j, y_j) - v_{\gamma} (x_j, y_j) \le 0.
		\]
		Now, considering that, as shown in Lemma \ref{lem_USC_v}, $v_{\gamma} (x, y)$ is upper semi-continuous at $(\bar{x}, \bar{y})$, and $f(x,y)$, and $g(x,y)$ are continuous at $(\bar{x}, \bar{y})$, we can conclude that 
		\begin{equation}\label{Thm_app_eq2}
			f(\bar{x}, \bar{y}) + g(\bar{x}, \bar{y}) - v_{\gamma} (\bar{x}, \bar{y}) \le 0.
		\end{equation}
		From (\ref{Thm_app_eq1}), as $f(x_k, y_k) + g(x_k, y_k) - v_{\gamma} (x_k, y_k) \ge 0$, we have
		\[
		F(x_k, y_k) <  \inf_{x \in X, y \in Y} \left\{ F(x,y) ~~\mathrm{s.t.} ~~ \varphi(x,y)-v_\gamma (x,y)\leq 0\right\} + \epsilon.
		\]
		By taking $k = j$ and letting $j \rightarrow \infty$ in the above inequality, and considering that $F(x,y)$ is continuous at $(\bar{x}, \bar{y})$, we obtain
		\[
		F(\bar{x}, \bar{y}) \le  \inf_{x \in X, y \in Y} \left\{ F(x,y) ~~\mathrm{s.t.} ~~ \varphi(x,y)-v_\gamma (x,y)\leq 0\right\} + \epsilon.
		\]
		Due to the arbitrariness of $\epsilon$, we can conclude that
		\[
		F(\bar{x}, \bar{y}) \le  \inf_{x \in X, y \in Y} \left\{ F(x,y) ~~\mathrm{s.t.} ~~ \varphi(x,y)-v_\gamma (x,y)\leq 0\right\},
		\]
		and then the conclusion follows from $(\bar{x}, \bar{y}) \in X \times Y$ and (\ref{Thm_app_eq2}).
	\end{proof}
	
	\subsection{Properties of Moreau Envelope}
	\label{propertiesME}
	
	By invoking Theorem 1 of \cite{rockafellar1974conjugate}, we have that when the LL problem is fully convex, the Moreau Envelope $v_{\gamma}(x,y)$ is also convex. We further generalize this finding in the subsequent lemma, showing that $v_{\gamma}(x,y)$ retains weak convexity when the LL problem exhibits weak convexity. The foundation for this proof draws inspiration from Theorem 2 in \cite{gao2023moreau}.
	
	\begin{lemma}\label{Lem1-a}
		Suppose that $\varphi(x, y)$ is  $(\rho_{\varphi_1},\rho_{\varphi_2})$-weakly convex on $X\times Y$. 
		Then for $\gamma \in (0, \frac{1}{2\rho_{\varphi_2}})$, $\rho_{v_1} \ge \rho_{\varphi_1}$ and $\rho_{v_2}  \ge \frac{1}{ \gamma} $,
		the function 
		$$v_{\gamma}(x,y) 
		+ \frac{\rho_{v_1}}{2}\|x\|^2 + \frac{\rho_{v_2}}{2}\|y\|^2$$ is convex on $X \times \mathbb{R}^m$. 
	\end{lemma}
	\begin{proof}
		We first extend the definition of the Moreau envelope $ v_{\gamma}(x,y) $ from $x\in X$  to $x\in \mathbb{R}^n$ by
		$$
		v_{\gamma}(x,y) := 
		\inf_{\theta \in \mathbb{R}^m}  \left\{ \varphi(x,\theta) + \frac{1}{2\gamma}\|\theta - y\|^2 + \delta_{X \times Y}(x,\theta) \right\}
		\quad \forall x\in \mathbb{R}^n, y \in \mathbb{R}^m.
		$$
		It follows that $v_\gamma(x,y)=+\infty$ for $x\not \in X$.
		For any $\rho_{v_1},  \rho_{v_2} > 0$, the function $v_{\gamma}(x,y) + \frac{\rho_{v_1}}{2}\|x\|^2 + \frac{\rho_{v_2}}{2}\|y\|^2 $ can be rewritten as 
		\[
		\begin{aligned}
			&v_{\gamma}(x,y) + \frac{\rho_{v_1}}{2}\|x\|^2 + \frac{\rho_{v_2}}{2}\|y\|^2 \\ 
			=\, & \inf_{\theta \in \mathbb{R}^m}  \left\{ \phi_{\gamma, \rho_v}(x, y, \theta) : = \varphi(x,\theta) + \frac{\rho_{v_1}}{2}\|x\|^2 + \frac{\rho_{v_2}}{2}\|y\|^2 + \frac{1}{2\gamma}\|\theta - y\|^2 + \delta_{X \times Y}(x,\theta) \right\}.
		\end{aligned}
		\]
		By direct computations, we obtain the following equation,
		\[
		\begin{aligned}
			&\phi_{\gamma, \rho_v}(x, y, \theta) \\ 
			= \,& \varphi(x,\theta) + \frac{\rho_{v_1}}{2}\|x\|^2 + \frac{\rho_{\varphi_2}}{2}\| \theta\|^2 + \delta_{X \times Y}(x,\theta) + \left( \frac{1}{2\gamma} - \frac{\rho_{\varphi_2}}{2} \right)\|\theta\|^2 + \frac{1 + \gamma\rho_{v_2}}{2\gamma}\|y\|^2 - \frac{1}{\gamma} \langle \theta,y\rangle.
		\end{aligned}
		\]
		Given that $\rho_{v_1} \ge \rho_{\varphi_1}$, 
		the convexity of $\varphi(x,\theta) + \frac{\rho_{v_1}}{2}\|x\|^2 + \frac{\rho_{\varphi_2}}{2}\| \theta\|^2 + \delta_{X \times Y}(x,\theta) $ can be immediately inferred, given that $\varphi(x,y)$ is $( \rho_{\varphi_1},\rho_{\varphi_2} )$-weakly convex on $X \times Y$.
		
		Further, when $\gamma \in (0, \frac{1}{2\rho_{\varphi_2}})$ and $\rho_{v_2}  \ge \frac{1}{ \gamma} $, it can be shown that both conditions,
		$ \frac{1}{4\gamma} - \frac{\rho_{\varphi_2}}{2}>0$ and  
		$\frac{1 + \gamma\rho_{v_2}}{2} \geq 1$, hold. This implies that the function
		\begin{align*}
			&\left( \frac{1}{2\gamma} - \frac{\rho_{\varphi_2}}{2} \right)\|\theta\|^2  + \frac{1 + \gamma\rho_{v_2}}{2\gamma}\|y\|^2 - \frac{1}{\gamma} \langle \theta,y\rangle \\
			=& \left( \frac{1}{4\gamma} - \frac{\rho_{\varphi_2}}{2} \right)\|\theta\|^2 +  \frac{1}{\gamma} \left( \frac{1}{4}\|\theta\|^2 + \frac{1 + \gamma\rho_{v_2}}{2}\|y\|^2 -  \langle \theta,y\rangle  \right),
		\end{align*}
		is convex with respect to $(y, \theta)$. 
		Therefore, under the conditions $\gamma \in (0, \frac{1}{2\rho_{\varphi_2}})$, $\rho_{v_1} \ge \rho_{\varphi_1}$ and $\rho_{v_2}  \ge \frac{1}{ \gamma} $, the extended-valued function $\phi_{\gamma, \rho_v}(x, y, \theta)$ is convex with respect to $(x, y, \theta)$ over $\mathbb{R}^n \times \mathbb{R}^m \times \mathbb{R}^m$. 
		This, in turn , establishes the convexity of 
		$$
		v_{\gamma}(x,y) + \frac{\rho_{v_1}}{2}\|x\|^2 + \frac{\rho_{v_2}}{2}\|y\|^2  = \inf_{\theta \in \mathbb{R}^m} \phi_{\gamma, \rho_v}(x,y,\theta)
		$$ 
		over $X \times \mathbb{R}^m$ by leveraging Theorem 1 of \cite{rockafellar1974conjugate}.
	\end{proof}
	
	Next we develop a calculus for the Moreau Envelope $v_{\gamma}(x,y)$, providing formulas for its gradient. These results immediately give insights into the proposed algorithm. The proof closely follows that of Theorem 5 in \cite{gao2023moreau}.
	
	\begin{lemma}\label{lem2-a}
		Under Assumption of Lemma \ref{Lem1-a},
		suppose that the gradient $\nabla_{x} g(x,y)$ exists and is continuous on $X \times Y$.
		The for $\gamma \in (0, \frac{1}{2\rho_{\varphi_2}})$,  
		${S}_{\gamma}( {x}, {y})=\{\theta_{\gamma}^*(x,y)\}$ is a singleton. Furthermore, 
		\begin{equation}\label{valuefinc-a}
			\begin{aligned}
				\nabla v_\gamma({x},{y}) = \Big( \nabla_x f( {x},   \theta_{\gamma}^*(x,y) ) + \nabla_{x} g( {x},   \theta_{\gamma}^*(x,y) ),\,
				\left( {y} -  \theta_{\gamma}^*(x,y)\right)/ \gamma \Big).
			\end{aligned}
		\end{equation}
	\end{lemma}
	\begin{proof}
		Considering $\gamma \in (0, \frac{1}{2\rho_{\varphi_2}})$ and the weakly convexity of $\varphi(x,y)$, the function 
		$\varphi(x,\theta) + \frac{1}{2\gamma}\|\theta - y\|^2 + \delta_{Y}(\theta)$ is shown to be $(\frac{1}{\gamma}-\frac{1}{\rho_{\varphi_2}})$-strongly convex with respect to $\theta$. Consequently, ${S}_{\gamma}( {x}, {y})=\{\theta_{\gamma}^*(x,y)\}$ is a singleton.
		
		Further, for $\gamma \in (0, \frac{1}{2\rho_{\varphi_2}})$, we have $\rho_{v_1} \ge \rho_{\varphi_1}$ and $\rho_{v_2}  \ge \frac{1}{ \gamma} $.
		Leveraging Lemma \ref{Lem1-a} and its subsequent proof, the function $v_{\gamma}(x,y) + \frac{\rho_{v_1}}{2}\|x\|^2 + \frac{\rho_{v_2}}{2}\|y\|^2$ is established as convex, and for any $(x,y)\in X\times Y$, the following holds
		\begin{equation*}
			v_{\gamma}(x,y) + \frac{\rho_{v_1}}{2}\|x\|^2 + \frac{\rho_{v_2}}{2}\|y\|^2  
			= \inf_{\theta \in Y} 
			\left\{
			\varphi(x,\theta) + \frac{\rho_{v_1}}{2}\|x\|^2 + \frac{\rho_{v_2}}{2}\|y\|^2 + \frac{1}{2\gamma}\|\theta - y\|^2
			\right\}, 
		\end{equation*}
		where $\varphi(x,\theta) + \frac{\rho_{v_1}}{2}\|x\|^2 + \frac{\rho_{v_2}}{2}\|y\|^2 + \frac{1}{2\gamma}\|\theta - y\|^2$ is convex with respect to $(x,y,\theta)$. By applying Theorem 3 of \cite{ye2023difference} and exploiting the continuously differentiable property of $g(x,y)$ with respect to $x$, the desired formulas are derived.
	\end{proof}

	\subsection{Auxiliary lemmas} \label{lemmas}
	
	In this section, we present auxiliary lemmas crucial for the non-asymptotic convergence analysis.
	
	\begin{lemma}\label{lem3}
		Let $\gamma \in (0, \frac{1}{2\rho_{\varphi_2}})$,  $(\bar{x}, \bar{y}) \in X \times \mathbb{R}^m$. Then for any $\rho_{v_1} \ge \rho_{\varphi_1}$ , $\rho_{v_2}  \ge \frac{1}{ \gamma} $ and $(x,y)$ on $X\times \mathbb{R}^m$, the following inequality holds:
		\begin{equation}
			-v_{\gamma}(x,y)  \le - v_\gamma(\bar{x},\bar{y}) - \langle \nabla v_\gamma(\bar{x},\bar{y}), (x,y) - (\bar{x},\bar{y})\rangle + \frac{\rho_{v_1}}{2} \| x - \bar{x} \|^2 + \frac{\rho_{v_2}}{2} \|  y - \bar{y} \|^2.
		\end{equation}
	\end{lemma}
	\begin{proof}
		According to Lemma \ref{Lem1-a}, $v_{\gamma}(x,y) + \frac{\rho_{v_1}}{2}\|x\|^2 + \frac{\rho_{v_2}}{2}\|y\|^2$ is convex on $X \times \mathbb{R}^m$. As a result, for any $(x,y)$ on $X \times \mathbb{R}^m$,
		\[
		\begin{aligned}
			&v_{\gamma}(x,y) + \frac{\rho_{v_1}}{2}\|x\|^2 + \frac{\rho_{v_2}}{2}\|y\|^2\\
			\ge \, &v_\gamma(\bar{x},\bar{y}) + \frac{\rho_{v_1}}{2}\|\bar{x}\|^2 + \frac{\rho_{v_2}}{2}\|\bar{y}\|^2  + \langle \nabla v_\gamma(\bar{x},\bar{y}) + (\rho_{v_1}\bar{x}, \rho_{v_2}\bar{y}), (x,y) - (\bar{x},\bar{y})\rangle.
		\end{aligned}
		\]
		Consequently, the conclusion follows directly.
	\end{proof}
	
	\begin{lemma}\label{lem_Proxg}
		For any $0 < s < 1/\rho_{g_2}$, and $\theta, \theta' \in \mathbb{R}^m$, the following inequality is satisfied:
		\begin{equation}
			\| \mathrm{Prox}_{s \tilde{g}(x,\cdot)} (\theta) - \mathrm{Prox}_{s \tilde{g}(x,\cdot)}(\theta')\| \le 1/(1-s\rho_{g_2}) \| \theta - \theta' \|.
		\end{equation}
	\end{lemma}
	\begin{proof}
		Let us denote $\mathrm{Prox}_{s \tilde{g}(x,\cdot)} (\theta) $ and $\mathrm{Prox}_{s \tilde{g}(x,\cdot)}(\theta')$ by $\theta^+$ and ${\theta'}^+$, respectively. From the definitions, we have
		\[
		0 \in \partial_y \tilde{g}(x,\theta^+) + \frac{1}{s} (\theta^+ - \theta ),
		\]
		and 
		\[
		0 \in \partial_y \tilde{g}(x, {\theta'}^+) + \frac{1}{s} ({\theta'}^+ - \theta' ).
		\]
		Given the $\rho_{g_2}$-weakly convexity of $\tilde{g}(x, \cdot)$, it implies
		\[
		\left\langle -  \frac{1}{s} (\theta^+ - \theta ) + \frac{1}{s} ({\theta'}^+ - \theta' ),  \theta^+ -{\theta'}^+ \right\rangle \ge - \rho_{g_2} \|  \theta^+ -{\theta'}^+ \|^2.
		\]
		From the above, the desired conclusion follows directly.
	\end{proof}
	\begin{lemma}\label{lem_theta}
		Let $\gamma \in (0, \frac{1}{\rho_{f_2} + 2\rho_{g_2} })$. Then, there exists $L_\theta > 0$ such that for any $({x}, {y}), (x', y') \in X \times \mathbb{R}^m$, the following inequality holds:
		\begin{equation}
			\| \theta_\gamma^*(x,y) - \theta_\gamma^*(x',y') \| \le L_\theta \| (x,y) - (x',y')\|.
		\end{equation}
	\end{lemma}	
	\begin{proof}
		Given that $\theta_\gamma^*(x,y)$ is optimal for the convex optimization problem $\min_{\theta \in Y} \, \varphi(x, \theta) + \frac{1}{2\gamma}\|\theta-y\|^2 $, we have
		\begin{equation*}
			\begin{aligned}
				0 &\in\nabla_{y} f(x, \theta_\gamma^*(x,y)) + \partial_{y} g(x, \theta_\gamma^*(x,y))  + (\theta_\gamma^*(x,y) - y)/\gamma + \mathcal{N}_Y (\theta_\gamma^*(x,y)), \\
				0 &\in \nabla_{y} f(x', \theta_\gamma^*(x',y')) + \partial_{y} g(x', \theta_\gamma^*(x',y')) + (\theta_\gamma^*(x',y') - y')/\gamma + \mathcal{N}_Y (\theta_\gamma^*(x',y')).
			\end{aligned}
		\end{equation*}
		Due to the $\rho_{g_2}$-weakly convexity of $\tilde{g}(x,y) := g(x,y) + \delta_{Y}(y)$ with respect to $y$, we obtain
		\begin{equation}\label{lem_theta_eq1}
			\begin{aligned}
				\theta_\gamma^*(x,y) &= \mathrm{Prox}_{s \tilde{g}(x,\cdot)} \left( \theta_\gamma^*(x,y) - s\left(\nabla_{y} f(x, \theta_\gamma^*(x,y)) + (\theta_\gamma^*(x,y) - y)/\gamma \right) \right),  \\
				\theta_\gamma^*(x',y') &= \mathrm{Prox}_{s \tilde{g}(x',\cdot)} \left( \theta_\gamma^*(x',y') - s\left(\nabla_{y} f(x', \theta_\gamma^*(x',y')) + (\theta_\gamma^*(x',y') - y')/\gamma \right) \right),
			\end{aligned}
		\end{equation}
		when $0 < s < 1/\rho_{g_2}$.
		Consequently, we deduce that
		\begin{equation}\label{lem_theta_eq3}
			\begin{aligned}
				&\| \theta_\gamma^*(x,y) - \theta_\gamma^*(x',y')\| \\ = \, & \big\|  \mathrm{Prox}_{s \tilde{g}(x,\cdot)} \left( \theta_\gamma^*(x,y) - s\left(\nabla_{y} f(x, \theta_\gamma^*(x,y)) + (\theta_\gamma^*(x,y) - y)/\gamma \right) \right) \\
				& \quad -  \mathrm{Prox}_{s \tilde{g}(x',\cdot)} \left( \theta_\gamma^*(x',y') - s\left(\nabla_{y} f(x', \theta_\gamma^*(x',y')) + (\theta_\gamma^*(x',y') - y')/\gamma \right) \right) \big\| \\
				\le \, & \big\|  \mathrm{Prox}_{s \tilde{g}(x,\cdot)} \left( \theta_\gamma^*(x,y) - s\left(\nabla_{y} f(x, \theta_\gamma^*(x,y)) + (\theta_\gamma^*(x,y) - y)/\gamma \right) \right) \\
				& \quad -  \mathrm{Prox}_{s \tilde{g}(x,\cdot)} \left( \theta_\gamma^*(x',y') - s\left(\nabla_{y} f(x', \theta_\gamma^*(x',y')) + (\theta_\gamma^*(x',y') - y')/\gamma \right) \right) \big\| \\
				& +  \big\|  \mathrm{Prox}_{s \tilde{g}(x,\cdot)} \left( \theta_\gamma^*(x',y') - s\left(\nabla_{y} f(x', \theta_\gamma^*(x',y')) + (\theta_\gamma^*(x',y') - y')/\gamma \right) \right) \\
				& \qquad -  \mathrm{Prox}_{s \tilde{g}(x',\cdot)} \left( \theta_\gamma^*(x',y') - s\left(\nabla_{y} f(x', \theta_\gamma^*(x',y')) + (\theta_\gamma^*(x',y') - y')/\gamma \right) \right) \big\| \\
				\le \, & \big\|  \mathrm{Prox}_{s \tilde{g}(x,\cdot)} \left( \theta_\gamma^*(x,y) - s\left(\nabla_{y} f(x, \theta_\gamma^*(x,y)) + (\theta_\gamma^*(x,y) - y)/\gamma \right) \right) \\
				& \quad -  \mathrm{Prox}_{s \tilde{g}(x,\cdot)} \left( \theta_\gamma^*(x',y') - s\left(\nabla_{y} f(x, \theta_\gamma^*(x',y')) + (\theta_\gamma^*(x',y') - y)/\gamma \right) \right) \big\| \\
				& +  \big\|  \mathrm{Prox}_{s \tilde{g}(x,\cdot)} \left( \theta_\gamma^*(x',y') - s\left(\nabla_{y} f(x, \theta_\gamma^*(x',y')) + (\theta_\gamma^*(x',y') - y)/\gamma \right) \right)  \\
				& \qquad -  \mathrm{Prox}_{s \tilde{g}(x,\cdot)} \left( \theta_\gamma^*(x',y') - s\left(\nabla_{y} f(x', \theta_\gamma^*(x',y')) + (\theta_\gamma^*(x',y') - y')/\gamma \right) \right) \big\| \\
				& + L_{\tilde{g}} \|x - x' \|,
			\end{aligned}
		\end{equation}
		where the second inequality is a consequence of Assumption \ref{assump-LL} (iv), which states that $\| \mathrm{Prox}_{s \tilde{g}(x,\cdot)} (\theta) - \mathrm{Prox}_{s \tilde{g}(x',\cdot)}(\theta)\| \le L_{\tilde{g}} \|x - x' \|$ for any $\theta \in Y$ and $s \in (0, \bar{s}]$. 
		Invoking Lemma \ref{lem_Proxg}, for $0 < s < 1/\rho_{g_2}$, we derive 
		\begin{equation}\label{lem_theta_eq2}
			\| \mathrm{Prox}_{s \tilde{g}(x,\cdot)} (\theta) - \mathrm{Prox}_{s \tilde{g}(x,\cdot)}(\theta')\| \le 1/(1-s\rho_{g_2}) \| \theta - \theta' \| \quad \forall \theta, \theta' \in \mathbb{R}^m.
		\end{equation}
		Given that $f(x,\theta) + \frac{1}{2\gamma}\|\theta - y\|^2$ is $(\frac{1}{\gamma}- \rho_{f_2})$-strongly convex with respect to $\theta$ on $Y$, we have
		\[
		\begin{aligned}
			&\left\langle  \nabla_{y} f(x, \theta_\gamma^*(x,y)) + (\theta_\gamma^*(x,y) - y)/\gamma -  \nabla_{y} f(x, \theta_\gamma^*(x',y')) - (\theta_\gamma^*(x',y') - y)/\gamma, \theta_\gamma^*(x,y) - \theta_\gamma^*(x',y') \right\rangle \\
			& \ge \left(\frac{1}{\gamma}- \rho_{f_2} \right) \| \theta_\gamma^*(x,y) - \theta_\gamma^*(x',y') \|^2,
		\end{aligned}
		\]
		which implies that when $0 < s \le (1/\gamma - \rho_{f_2})/(L_f+1/\gamma)^2 $,
		\[
		\begin{aligned}
			&	\big\| \theta_\gamma^*(x,y) - s\left(\nabla_{y} f(x, \theta_\gamma^*(x,y)) + (\theta_\gamma^*(x,y) - y)/\gamma\right) - \theta_\gamma^*(x',y') \\
			&\ + s\left(\nabla_{y} f(x, \theta_\gamma^*(x',y')) + (\theta_\gamma^*(x',y') - y)/\gamma \right) \big\| ^2\\
			\le \, & 
			\left[ 
			1- 2s\left(1/\gamma - \rho_{f_2}\right) + s^2(L_f+1/\gamma)^2 
			\right] 
			\| \theta_\gamma^*(x,y) - \theta_\gamma^*(x',y') \|^2  \\
			\le \, & \left[ 1- s\left(1/\gamma - \rho_{f_2}\right) \right] \| \theta_\gamma^*(x,y) - \theta_\gamma^*(x',y') \|^2 .
		\end{aligned}
		\]
		Combining this with (\ref{lem_theta_eq2}), we infer that
		\begin{equation}\label{lem_theta_eq4}
			\begin{aligned}
				&\big\|  \mathrm{Prox}_{s \tilde{g}(x,\cdot)} \left( \theta_\gamma^*(x,y) - s\left(\nabla_{y} f(x, \theta_\gamma^*(x,y)) + (\theta_\gamma^*(x,y) - y)/\gamma \right) \right) \\
				& \quad -  \mathrm{Prox}_{s \tilde{g}(x,\cdot)} \left( \theta_\gamma^*(x',y') - s\left(\nabla_{y} f(x, \theta_\gamma^*(x',y')) + (\theta_\gamma^*(x',y') - y)/\gamma \right) \right) \big\| \\ 
				\le \, & 1/(1-s\rho_{g_2}) 	\big\| \theta_\gamma^*(x,y) - s\left(\nabla_{y} f(x, \theta_\gamma^*(x,y)) + (\theta_\gamma^*(x,y) - y)/\gamma\right) \\
				& \qquad \qquad \qquad - \theta_\gamma^*(x',y') + s\left(\nabla_{y} f(x, \theta_\gamma^*(x',y')) + (\theta_\gamma^*(x',y') - y)/\gamma \right) \big\| \\
				\le \, & \sqrt{1- s\left(1/\gamma - \rho_{f_2}\right)} /(1-s\rho_{g_2})  \| \theta_\gamma^*(x,y) - \theta_\gamma^*(x',y') \|.
			\end{aligned}
		\end{equation}
		Next, utilizing Lemma \ref{lem_Proxg}, for $0 < s < 1/\rho_{g_2}$, it follows that
		\begin{equation}\label{lem_theta_eq5}
			\begin{aligned}
				&\big\|  \mathrm{Prox}_{s \tilde{g}(x,\cdot)} \left( \theta_\gamma^*(x',y') - s\left(\nabla_{y} f(x, \theta_\gamma^*(x',y')) + (\theta_\gamma^*(x',y') - y)/\gamma \right) \right)  \\
				& \qquad -  \mathrm{Prox}_{s \tilde{g}(x,\cdot)} \left( \theta_\gamma^*(x',y') - s\left(\nabla_{y} f(x', \theta_\gamma^*(x',y')) + (\theta_\gamma^*(x',y') - y')/\gamma \right) \right) \big\| \\
				\le \, & 1/(1-s\rho_{g_2}) 	\big\| \theta_\gamma^*(x',y') - s\left(\nabla_{y} f(x, \theta_\gamma^*(x',y')) + (\theta_\gamma^*(x',y') - y)/\gamma\right) \\
				& \qquad \qquad \qquad - \theta_\gamma^*(x',y') + s\left(\nabla_{y} f(x', \theta_\gamma^*(x',y')) + (\theta_\gamma^*(x',y') - y')/\gamma \right) \big\| \\
				\le \, & s/(1-s\rho_{g_2}) \left( \| \nabla_{y} f(x, \theta_\gamma^*(x',y'))  -  \nabla_{y} f(x', \theta_\gamma^*(x',y')) \| + \|y - y'\|/\gamma \right) \\
				\le \, & s/(1-s\rho_{g_2})\left(L_f \| x - x'\| + \frac{1}{\gamma} \|y - y'\| \right).
			\end{aligned}
		\end{equation}
		From estimates (\ref{lem_theta_eq3}), (\ref{lem_theta_eq4}) and (\ref{lem_theta_eq5}), we deduce that, for any $s > 0$ satisfying $s \le (1/\gamma - \rho_{f_2})/(L_f+1/\gamma)^2$, $s \le \bar{s}$ and $s < 1/\rho_{g_2}$, the following condition holds
		\begin{equation}
			\begin{aligned}
				\| \theta_\gamma^*(x,y) - \theta_\gamma^*(x',y')\| \le \, &\sqrt{1- s\left(1/\gamma - \rho_{f_2}\right)} /(1-s\rho_{g_2})  \| \theta_\gamma^*(x,y) - \theta_\gamma^*(x',y') \| \\
				&+ s/(1-s\rho_{g_2})\left(L_f \| x - x'\| + \frac{1}{\gamma} \|y - y'\| \right) +  L_{\tilde{g}} \|x - x' \|.
			\end{aligned}
		\end{equation}
		Given that $\gamma < \frac{1}{\rho_{f_2} + 2\rho_{g_2} }$, it can be inferred that $1/\gamma - \rho_{f_2} > 2\rho_{g_2}$. This implies $1 - 2s\rho_{g_2} > 1 - s(1/\gamma - \rho_{f_2})$, leading to $1 - s(1/\gamma - \rho_{f_2}) < (1-s\rho_{g_2})^2$. Consequently, we deduce $\sqrt{1- s\left(1/\gamma - \rho_{f_2}\right)} /(1-s\rho_{g_2}) < 1$.
		From these derivations, the desired conclusion is evident.
	\end{proof}

	\begin{lemma}\label{lem4}
		Suppose $\gamma \in (0, \frac{1}{\rho_{f_2} + 2\rho_{g_2} })$ and $\eta_k \in (0,(1/\gamma - \rho_{f_2})/(L_f+1/\gamma)^2 ] \cap (0,1/\rho_{g_2})$, the sequence of $(x^k, y^k, \theta^k)$ generated by Algorithm \ref{MEHA} satisfies
		\begin{equation}\label{lem4_eq}
			\| \theta^{k+1} - \theta_{\gamma}^*(x^k,y^k)\| \le\sigma_k \| \theta^{k} - \theta_{\gamma}^*(x^k,y^k)\|,
		\end{equation}
		where $\sigma_k:= \sqrt{1- \eta_k\left(1/\gamma - \rho_{f_2}\right)} /(1- \eta_k\rho_{g_2}) < 1 $. 
	\end{lemma}
	\begin{proof}
		Recalling (\ref{lem_theta_eq1}) from Lemma \ref{lem_theta} that when $\eta_k  < 1/\rho_{g_2}$, 
		\[
		\theta_\gamma^*(x^k, y^k) = \mathrm{Prox}_{\eta_k \tilde{g}(x^k,\cdot)} \left( \theta_\gamma^*(x^k, y^k) - \eta_k\left(\nabla_{y} f(x^k, \theta_\gamma^*(x,y)) + (\theta_\gamma^*(x^k, y^k) - y^k)/\gamma \right) \right).
		\]
		Considering the update rule for $\theta^{k+1}$ as defined in (\ref{update_z}) and using arguments analogous to those in the derivation of (\ref{lem_theta_eq4}) from Lemma \ref{lem_theta}, 
		when $ \eta_k \le (1/\gamma - \rho_{f_2})/(L_f+1/\gamma)^2 $, it follows
		\[
		\| \theta^{k+1} - \theta_{\gamma}^*(x^k,y^k)\| \le \sigma_k  \| \theta^{k} - \theta_{\gamma}^*(x^k,y^k)\|,
		\]
		where $\sigma_k:= \sqrt{1- \eta_k\left(1/\gamma - \rho_{f_2}\right)} /(1- \eta_k\rho_{g_2}) $. Notably, $\sigma_k < 1$ is a consequence of $\gamma < \frac{1}{\rho_{f_2} + 2\rho_{g_2} }$.
	\end{proof}
	
	The update of variables $(x,y)$ in (\ref{update_x}) and (\ref{update_y}) can be interpreted as inexact alternating proximal gradient from $(x^{k}, y^{k})$ on $\min_{(x,y)\in X \times Y} \phi_{c_k}(x,y)$, in which $\phi_{c_k}$ is defined in (\ref{def_varphi}) as 
	\[
	\phi_{c_k}(x,y) := \frac{1}{ c_k}\big(F(x,y) - \underline{F} \big) +   f(x,y) + g(x,y) - v_\gamma (x,y).
	\]
	The subsequent lemma illustrates that the function 
	$\phi_{c_k}(x,y)$ exhibits a monotonic decreasing behavior with errors at each iteration.

	\begin{lemma}\label{lem6}
		Under Assumptions \ref{assump-UL} and \ref{assump-LL}, 
		suppose $\gamma \in (0, \frac{1}{2\rho_{f_2} + 2\rho_{g_2} })$ and $\beta_k < 1/ \rho_{g_2}$, the sequence of $(x^k, y^k, \theta^k)$ generated by Algorithm \ref{MEHA} satisfies
		\begin{equation}\label{lem6_eq}
			\begin{aligned}
				\phi_{c_k}(x^{k+1},y^{k+1}) \le \,& \phi_{c_k}(x^k, y^{k}) - \left( \frac{1}{2\alpha_k} - \frac{L_{\phi_k}}{2} -  \frac{\beta_k L_\theta^2}{\gamma^2}  \right)\| x^{k+1} - x^k \|^2 \\
				& - \left( \frac{1}{2\beta_k} - \frac{\rho_{g_2}}{2} - \frac{L_{\phi_k}}{2} \right)\| y^{k+1} - y^k \|^2 \\
				& + \left(\frac{\alpha_k}{2} (L_f + L_g)^2 + \frac{\beta_k}{\gamma^2} \right) \left\| \theta^{k+1} - \theta_{\gamma}^*(x^k,y^{k}) \right\|^2,
			\end{aligned}
		\end{equation}
		where $L_{\phi_k} := L_F/c_k + L_f + L_g + \max\{ \rho_{\varphi_1} , 1/\gamma \}$.
	\end{lemma}
	\begin{proof}
		Under the conditions of Assumption \ref{assump-UL}, the functions $F$ and $f$ exhibit $L_F$- and $L_f$-smooth on $X \times Y$, respectively. Further, according to Assumption \ref{assump-LL}, the function $g(\cdot,y^k)$ is $L_g$-smooth on $X$. Leveraging these assumptions and invoking Lemma \ref{lem3}, we deduce
		\begin{equation}\label{lem6_eq1}
			\begin{aligned}
				\phi_{c_k}(x^{k+1},y^{k}) \le \,& \phi_{c_k}(x^k, y^{k}) + \langle \nabla_x \phi_{c_k}(x^k, y^{k}), x^{k+1} - x^k \rangle + \frac{L_{\phi_k}}{2} \| x^{k+1} - x^k \|^2,
			\end{aligned}
		\end{equation}
		with $L_{\phi_k} := L_F/c_k + L_f + L_g + \max\{ \rho_{\varphi_1} , 1/\gamma \}$. Considering the update rule for the variable $x$ as defined in (\ref{update_x}) and leveraging the property of the projection operator $\mathrm{Proj}_{X}$, it follows that
		\[
		\langle x^k - \alpha_k  d_x^k -  x^{k+1}, x^k - x^{k+1} \rangle \le 0,
		\]
		leading to
		\[
		\langle  d_x^k, x^{k+1} - x^k \rangle \le -\frac{1}{\alpha_k} \| x^{k+1} - x^k \|^2.
		\]
		Combining this inequality with (\ref{lem6_eq1}), if can be deduced that
		\begin{equation}\label{lem6_eq2}
			\begin{aligned}
				\phi_{c_k}(x^{k+1},y^{k}) \le \,& \phi_{c_k}(x^k, y^{k}) - \left( \frac{1}{\alpha_k} - \frac{L_{\phi_k}}{2} \right)\| x^{k+1} - x^k \|^2 \\
				& + \left\langle  \nabla_x \phi_{c_k}(x^k, y^{k}) - d_x^k, x^{k+1} - x^k \right\rangle.
			\end{aligned}
		\end{equation}
		Given the expression for $\nabla v_\gamma({x},{y})$ as derived in Lemma \ref{lem2-a} and the definition of $d_x^k$ provided in (\ref{def_dx}), we obtain
		\begin{equation}\label{lem6_eq3}
			\begin{aligned}
				&\left\|  \nabla_x \phi_{c_k}(x^k, y^{k}) - d_x^k \right\|^2 \\
				= \, &   \left\|  \nabla_x f( x^k,   \theta_{\gamma}^*(x^k,y^{k}) ) + \nabla_x g( x^k,   \theta_{\gamma}^*(x^k,y^{k}) ) -  \nabla_x f( x^k,   \theta^{k+1} ) - \nabla_x g( x^k,   \theta^{k+1}  ) \right\|^2 \\
				\le \, & (L_f + L_g)^2 \left\| \theta^{k+1} - \theta_{\gamma}^*(x^k,y^{k}) \right\|^2. 
			\end{aligned}
		\end{equation}
		This yields
		\[
		\begin{aligned}
			& \left\langle  \nabla_x \phi_{c_k}(x^k, y^{k}) - d_x^k, x^{k+1} - x^k \right\rangle \\
			\le \,  & \frac{\alpha_k}{2} (L_f + L_g)^2 \left\| \theta^{k+1} - \theta_{\gamma}^*(x^k,y^{k}) \right\|^2 + \frac{1}{2\alpha_k} \| x^{k+1} - x^k \|^2,
		\end{aligned}
		\]
		which combining with (\ref{lem6_eq2}) leads to
		\begin{equation}\label{lem6_eq4}
			\begin{aligned}
				\phi_{c_k}(x^{k+1},y^{k}) \le \,& \phi_{c_k}(x^k, y^{k}) - \left( \frac{1}{2\alpha_k} - \frac{L_{\phi_k}}{2} \right)\| x^{k+1} - x^k \|^2\\
				& + \frac{\alpha_k}{2} (L_f + L_g)^2 \left\| \theta^{k+1} - \theta_{\gamma}^*(x^k,y^{k}) \right\|^2.
			\end{aligned}
		\end{equation}
		Considering the update rule for variable $y$ given by (\ref{update_y}), and the $\rho_{g_2}$-weakly convex property of $g(x^{k+1}, \cdot)$ over $Y$, it follows that for $\beta_k < 1/\rho_{g_2}$,
		\begin{equation}\label{lem6_eq5}
			\langle  d_y^k, y^{k+1} - y^k \rangle +  g(x^{k+1}, y^{k+1})  + \left( \frac{1}{\beta_k} - \frac{\rho_{g_2}}{2}\right) \| y^{k+1} - y^k \|^2 \le g(x^{k+1},y^k).
		\end{equation}
		Under Assumption \ref{assump-UL}, where $F$ and $f$ are $L_F$- and $l_f$-smooth on $X \times Y$, respectively, and invoking Lemma \ref{lem3}, we deduce
		\begin{equation}
			\begin{aligned}
				&\phi_{c_k}(x^{k+1},y^{k+1}) - g(x^{k+1}, y^{k+1}) \\
				\le \,& \phi_{c_k}(x^{k+1}, y^{k}) - g(x^{k+1}, y^{k})+ \langle \nabla_y \left(\phi_{c_k} - g \right)(x^{k+1}, y^{k}), y^{k+1} - y^k \rangle + \frac{L_{\phi_k}}{2} \| y^{k+1} - y^k \|^2.
			\end{aligned}
		\end{equation}
		Combining this inequality with (\ref{lem6_eq5}), we obtain
		\begin{equation}\label{lem6_eq6}
			\begin{aligned}
				&\phi_{c_k}(x^{k+1},y^{k+1}) \\
				\le \,& \phi_{c_k}(x^{k+1}, y^{k}) - \left( \frac{1}{\beta_k} - \frac{\rho_{g_2}}{2} - \frac{L_{\phi_k}}{2} \right)\| y^{k+1} - y^k \|^2 \\
				& + \left\langle  \nabla_y \left(\phi_{c_k} - g \right)(x^{k+1}, y^{k}) - d_y^k, y^{k+1} - y^k \right\rangle.
			\end{aligned}
		\end{equation}
		Given the expression for $\nabla v_\gamma({x},{y})$ as derived in Lemma \ref{lem2-a} and the definition of $d_y^k$ from (\ref{def_dy}), we deduce
		\begin{equation}\label{lem6_eq7}
			\begin{aligned}
				\left\|  \nabla_y \left(\phi_{c_k} - g \right)(x^{k+1}, y^{k}) - d_y^k \right\|^2
				= \, &   \left\| ( {y}^k -  \theta_{\gamma}^*(x^{k+1},y^k))/ \gamma - ( {y}^k -  \theta^{k+1})/ \gamma \right\|^2 \\
				= \, & \frac{1}{\gamma^2} \left\| \theta^{k+1} - \theta_{\gamma}^*(x^{k+1},y^{k}) \right\|^2, 
			\end{aligned}
		\end{equation}
		and thus
		\[
		\begin{aligned}
			\left\langle  \nabla_y \left(\phi_{c_k} - g \right)(x^{k+1}, y^{k}) - d_y^k, y^{k+1} - y^k \right\rangle \le & \frac{\beta_k}{2\gamma^2} \left\| \theta^{k+1} - \theta_{\gamma}^*(x^{k+1},y^{k}) \right\|^2 + \frac{1}{2\beta_k} \| y^{k+1} - y^k \|^2.
		\end{aligned}
		\]
		Consequently, we have from (\ref{lem6_eq6}) that
		\begin{equation}
			\begin{aligned}
				&\phi_{c_k}(x^{k+1},y^{k+1}) \\
				\le \,& \phi_{c_k}(x^{k+1}, y^{k}) - \left( \frac{1}{2\beta_k} - \frac{\rho_{g_2}}{2} - \frac{L_{\phi_k}}{2} \right)\| y^{k+1} - y^k \|^2 + \frac{\beta_k}{2\gamma^2} \left\| \theta^{k+1} - \theta_{\gamma}^*(x^{k+1},y^{k}) \right\|^2 \\
				\le \, &\phi_{c_k}(x^{k+1}, y^{k}) - \left( \frac{1}{2\beta_k} - \frac{\rho_{g_2}}{2} - \frac{L_{\phi_k}}{2} \right)\| y^{k+1} - y^k \|^2 + \frac{\beta_k}{\gamma^2} \left\| \theta^{k+1} - \theta_{\gamma}^*(x^{k},y^{k}) \right\|^2 \\
				& + \frac{\beta_k L_\theta^2}{\gamma^2} \left\| x^{k+1} - x^{k} \right\|^2,
			\end{aligned}
		\end{equation}
		where the last inequality follows from Lemma \ref{lem_theta}. The conclusion follows by combining this with (\ref{lem6_eq4}).
	\end{proof}
	
	\subsection{Proof of Lemma~\ref{lem9}}\label{proofV}
	
	With the auxiliary lemmas from the preceding section, we demonstrate the decreasing property of the merit function $V_k$.  
	
	\begin{lemma}\label{lem9-a}
		Under Assumptions \ref{assump-UL} and \ref{assump-LL}, 
		suppose $\gamma \in (0, \frac{1}{2\rho_{f_2} + 2\rho_{g_2} })$, $c_{k+1} \ge c_k$ and $\eta_k \in [ \underline{\eta}, (1/\gamma - \rho_{f_2})/(L_f+1/\gamma)^2] \cap  [ \underline{\eta}, 1/\rho_{g_2})$ with $ \underline{\eta} > 0$, 
		then there exists $c_{\alpha}, c_\beta, c_{\theta}> 0$ such that when $0<\alpha_k \le {c_\alpha}$ and $0<\beta_k \le c_\beta$, the sequence of $(x^k, y^k, \theta^k)$ generated by MEHA (Algorithm \ref{MEHA}) satisfies
		\begin{equation}\label{lem9_eq-a}
			\begin{aligned}
				V_{k+1} - V_k 
				\le \, & -  \frac{1}{4\alpha_k} \| x^{k+1} - x^k \|^2 - \frac{1}{4\beta_k} \| y^{k+1} - y^k \|^2 
				- c_{\theta}
				\left\| \theta^{k} - \theta_{\gamma}^*(x^{k},y^{k}) \right\|^2,
			\end{aligned}
		\end{equation}
		where $c_{\theta}=\frac{1}{2}\left(\frac{\underline{\eta}\rho_{g_2}}{1-\underline{\eta}\rho_{g_2}} \right)^2 \left( (L_f + L_g)^2+ 1/\gamma^2  \right)$.
	\end{lemma}
	\begin{proof}
		Let us first recall (\ref{lem6_eq}) from Lemma \ref{lem6}, which states that
		\begin{equation}\label{lem9_eq1}
			\begin{aligned}
				\phi_{c_k}(x^{k+1},y^{k+1}) \le \,& \phi_{c_k}(x^k, y^{k}) - \left( \frac{1}{2\alpha_k} - \frac{L_{\phi_k}}{2} -  \frac{\beta_k L_\theta^2}{\gamma^2}  \right)\| x^{k+1} - x^k \|^2 \\
				& - \left( \frac{1}{2\beta_k} - \frac{\rho_{g_2}}{2} - \frac{L_{\phi_k}}{2} \right)\| y^{k+1} - y^k \|^2 \\
				& + \left(\frac{\alpha_k}{2} (L_f + L_g)^2 + \frac{\beta_k}{\gamma^2} \right) \left\| \theta^{k+1} - \theta_{\gamma}^*(x^k,y^{k}) \right\|^2,
			\end{aligned}
		\end{equation}
		when $\beta_k < 1/ \rho_{g_2}$.
		Since $c_{k+1} \ge c_k$,
		we can infer that $ (F(x^{k+1},y^{k+1})- \underline{F} )/c_{k+1} \le (F(x^{k+1},y^{k+1})- \underline{F} )/c_{k}$. 
		Combining this with (\ref{lem9_eq1}) leads to
		\begin{equation}\label{lem9_eq4}
			\begin{aligned}
				V_{k+1} - V_k =\, & \phi_{c_{k+1}}(x^{k+1},y^{k+1}) - \phi_{c_k}(x^k, y^{k})   \\
				&+ \left( (L_f + L_g)^2+ 1/\gamma^2 \right) \left\| \theta^{k+1} - \theta_{\gamma}^*(x^{k+1},y^{k+1}) \right\|^2 \\
				& - \left( (L_f + L_g)^2+ 1/\gamma^2 \right) \left\| \theta^{k} - \theta_{\gamma}^*(x^{k},y^{k}) \right\|^2\\
				\le\, & \phi_{c_{k}}(x^{k+1},y^{k+1}) - \phi_{c_k}(x^k, y^{k}) \\
				& +\left( (L_f + L_g)^2+ 1/\gamma^2 \right)  \left\| \theta^{k+1} - \theta_{\gamma}^*(x^{k+1},y^{k+1}) \right\|^2 \\
				& - \left( (L_f + L_g)^2+ 1/\gamma^2 \right)  \left\| \theta^{k} - \theta_{\gamma}^*(x^{k},y^{k}) \right\|^2 \\
				\le \, &- \left( \frac{1}{2\alpha_k} - \frac{L_{\phi_k}}{2} -  \frac{\beta_k L_\theta^2}{\gamma^2}  \right)\| x^{k+1} - x^k \|^2 \\
				&- \left( \frac{1}{2\beta_k} - \frac{\rho_{g_2}}{2} - \frac{L_{\phi_k}}{2} \right)\| y^{k+1} - y^k \|^2 \\
				& + \left( (L_f + L_g)^2+ 1/\gamma^2 \right)  \left\| \theta^{k+1} - \theta_{\gamma}^*(x^{k+1},y^{k+1}) \right\|^2\\
				& -  \left( (L_f + L_g)^2+ 1/\gamma^2 \right)  \left\| \theta^{k} - \theta_{\gamma}^*(x^{k},y^{k}) \right\|^2\\
				& + \left(\frac{\alpha_k}{2} (L_f + L_g)^2 + \frac{\beta_k}{\gamma^2} \right) \left\| \theta^{k+1} - \theta_{\gamma}^*(x^k,y^{k}) \right\|^2.
			\end{aligned}
		\end{equation}
		We can demonstrate that
		\begin{equation*}
			\begin{aligned}
				& \left\| \theta^{k+1} - \theta_{\gamma}^*(x^{k+1},y^{k+1}) \right\|^2 -  \left\| \theta^{k} - \theta_{\gamma}^*(x^{k},y^{k}) \right\|^2 + \frac{\alpha_k}{2} \left\| \theta^{k+1} - \theta_{\gamma}^*(x^k,y^{k}) \right\|^2 \\
				\le \, & (1+\epsilon_k+  \frac{\alpha_k}{2} )\left\| \theta^{k+1} - \theta_{\gamma}^*(x^k,y^{k}) \right\|^2 -  \left\| \theta^{k} - \theta_{\gamma}^*(x^{k},y^{k}) \right\|^2 \\&+ (1+ \frac{1}{\epsilon_k}) \| \theta_{\gamma}^*(x^{k+1},y^{k+1}) - \theta_{\gamma}^*(x^k,y^{k}) \|^2 \\
				\le \, & (1+\epsilon_k+  \frac{\alpha_k}{2} )\sigma_k^2\| \theta^{k} - \theta_{\gamma}^*(x^k,y^{k})\|^2 -  \left\| \theta^{k} - \theta_{\gamma}^*(x^{k},y^{k}) \right\|^2 \\& + (1+ \frac{1}{\epsilon_k}) L_\theta^2 \left\| (x^{k+1}, y^{k+1}) - (x^k, y^k) \right\|^2,
			\end{aligned}
		\end{equation*}
		for any $\epsilon_k >0$, where the second inequality is a consequence of Lemmas \ref{lem_theta} and \ref{lem4}. 
		Since $\gamma < \frac{1}{\rho_{f_2} + 2\rho_{g_2} }$, we have  $1 - 2\eta_k\rho_{g_2} > 1 - \eta_k(1/\gamma - \rho_{f_2})$, and thus $\sigma_k^2  = (1- \eta_k\left(1/\gamma - \rho_{f_2}\right)) /(1- \eta_k\rho_{g_2})^2\le  1 - \left(\frac{\eta_k\rho_{g_2}}{1-\eta_k\rho_{g_2}} \right)^2$.
		By setting $\epsilon_k = \frac{1}{4} \left(\frac{\eta_k\rho_{g_2}}{1-\eta_k\rho_{g_2}} \right)^2$ in the above inequality, we deduce that when $\alpha_k \le \frac{1}{2} \left(\frac{\eta_k\rho_{g_2}}{1-\eta_k\rho_{g_2}} \right)^2$, it follows that
		\begin{equation}\label{lem9_eq5}
			\begin{aligned}
				&  \left\| \theta^{k+1} - \theta_{\gamma}^*(x^{k+1},y^{k+1}) \right\|^2 -  \left\| \theta^{k} - \theta_{\gamma}^*(x^{k},y^{k}) \right\|^2 + \frac{\alpha_k}{2} \left\| \theta^{k+1} - \theta_{\gamma}^*(x^k,y^{k}) \right\|^2  \\
				\le \, & - \frac{1}{2}\left(\frac{\eta_k\rho_{g_2}}{1-\eta_k\rho_{g_2}} \right)^2 \left\| \theta^{k} - \theta_{\gamma}^*(x^{k},y^{k}) \right\|^2  \\
				& + \left(1+ 4 \left(\frac{1-\eta_k\rho_{g_2}}{\eta_k\rho_{g_2}} \right)^2 \right)  L_\theta^2  \left\| (x^{k+1}, y^{k+1}) - (x^k, y^k) \right\|^2.
			\end{aligned}
		\end{equation}
		Similarly, we can show that, when $\beta_k \le \frac{1}{4} \left(\frac{\eta_k\rho_{g_2}}{1-\eta_k\rho_{g_2}} \right)^2$, it holds that
		\begin{equation}\label{lem9_eq5.5}
			\begin{aligned}
				&  \left\| \theta^{k+1} - \theta_{\gamma}^*(x^{k+1},y^{k+1}) \right\|^2 -  \left\| \theta^{k} - \theta_{\gamma}^*(x^{k},y^{k}) \right\|^2 +\beta_k \left\| \theta^{k+1} - \theta_{\gamma}^*(x^k,y^{k}) \right\|^2  \\
				\le \, & - \frac{1}{2}\left(\frac{\eta_k\rho_{g_2}}{1-\eta_k\rho_{g_2}} \right)^2 \left\| \theta^{k} - \theta_{\gamma}^*(x^{k},y^{k}) \right\|^2  \\
				& + \left(1+ 4 \left(\frac{1-\eta_k\rho_{g_2}}{\eta_k\rho_{g_2}} \right)^2 \right)  L_\theta^2  \left\| (x^{k+1}, y^{k+1}) - (x^k, y^k) \right\|^2.
			\end{aligned}
		\end{equation}
		Combining (\ref{lem9_eq4}), (\ref{lem9_eq5}) and (\ref{lem9_eq5.5}), we have
		\begin{equation}\label{lem9_eq6}
			\begin{aligned}
				&	V_{k+1} - V_k \\
				\le \, & - \left[ \frac{1}{2\alpha_k} - \frac{L_{\phi_k}}{2} -  \frac{\beta_k L_\theta^2}{\gamma^2}  - \left(1+ 4 \left(\frac{1-\eta_k\rho_{g_2}}{\eta_k\rho_{g_2}} \right)^2 \right)  L_\theta^2 \left( (L_f + L_g)^2+ 1/\gamma^2 \right) \right] \| x^{k+1} - x^k \|^2 \\
				&  - \left[ \frac{1}{2\beta_k} - \frac{\rho_{g_2}}{2} - \frac{L_{\phi_k}}{2} - \left(1+ 4 \left(\frac{1-\eta_k\rho_{g_2}}{\eta_k\rho_{g_2}} \right)^2 \right)  L_\theta^2 \left( (L_f + L_g)^2+ 1/\gamma^2  \right)  \right]\| y^{k+1} - y^k \|^2 \\
				&-  \frac{1}{2}\left(\frac{\eta_k\rho_{g_2}}{1-\eta_k\rho_{g_2}} \right)^2 \left( (L_f + L_g)^2+ 1/\gamma^2  \right)  \left\| \theta^{k} - \theta_{\gamma}^*(x^{k},y^{k}) \right\|^2 .
			\end{aligned}
		\end{equation}
		When $c_{k+1} \ge c_k$, $\eta_k \ge \underline{\eta} > 0$, $\alpha_k \le \frac{1}{2} \left(\frac{\underline{\eta}\rho_{g_2}}{1-\underline{\eta}\rho_{g_2}} \right)^2$ and $\beta_k \le \frac{1}{4} \left(\frac{\underline{\eta}\rho_{g_2}}{1-\underline{\eta}\rho_{g_2}} \right)^2$ and , it holds that, for any $k$, $\alpha_k \le \frac{1}{2} \left(\frac{\eta_k\rho_{g_2}}{1-\eta_k\rho_{g_2}} \right)^2$, $\beta_k \le \frac{1}{4} \left(\frac{\eta_k\rho_{g_2}}{1-\eta_k\rho_{g_2}} \right)^2$,
		\begin{equation}\label{c-alpha}
			\begin{aligned}
				& \frac{L_{\phi_k}}{2} + \frac{\beta_k L_\theta^2}{\gamma^2}  + \left(1+ 4 \left(\frac{1-\eta_k\rho_{g_2}}{\eta_k\rho_{g_2}} \right)^2 \right)  L_\theta^2 \left( (L_f + L_g)^2+ 1/\gamma^2 \right) \\
				\le \, & \frac{L_{\phi_0}}{2} + \frac{L_\theta^2}{4\gamma^2} \left(\frac{\underline{\eta}\rho_{g_2}}{1-\underline{\eta}\rho_{g_2}} \right)^2  + \left(1+ 4 \left(\frac{1-\underline{\eta}\rho_{g_2}}{\underline{\eta}\rho_{g_2}} \right)^2 \right)  L_\theta^2 \left( (L_f + L_g)^2+ 1/\gamma^2 \right) =: C_\alpha, 
			\end{aligned}
		\end{equation}
		and
		\begin{equation}\label{c-beta}
			\begin{aligned}
				& \frac{\rho_{g_2}}{2} + \frac{L_{\phi_k}}{2} + \left(1+ 4 \left(\frac{1-\eta_k\rho_{g_2}}{\eta_k\rho_{g_2}} \right)^2 \right)  L_\theta^2 \left( (L_f + L_g)^2+ 1/\gamma^2  \right)   \\
				\le \, & \frac{\rho_{g_2}}{2}  + \frac{L_{\phi_0}}{2} + \left(1+ 4 \left(\frac{1-\underline{\eta}\rho_{g_2}}{\underline{\eta}\rho_{g_2}} \right)^2 \right)  L_\theta^2 \left( (L_f + L_g)^2+ 1/\gamma^2 \right) =: C_\beta, 
			\end{aligned}
		\end{equation}
		Consequently, since $\frac{1}{4C_\beta}<\frac{1}{2\rho_{g_2}}$, if $c_\alpha, c_\beta >0$ satisfies
		\begin{equation}
			c_{\alpha} 
			\leq
			\min \left\{ \frac{1}{2} \left(\frac{\underline{\eta}\rho_{g_2}}{1-\underline{\eta}\rho_{g_2}} \right)^2, \frac{1}{4C_\alpha}  \right\}, 
			\quad
			c_{\beta} 
			\leq
			\min \left\{  \frac{1}{4} \left(\frac{\underline{\eta}\rho_{g_2}}{1-\underline{\eta}\rho_{g_2}} \right)^2, \frac{1}{4C_\beta}  \right\}, 
		\end{equation}
		then, when $0 < \alpha_k \le c_\alpha$ and $0 < \beta_k \le c_\beta$, it holds that
		\[
		\frac{1}{2\alpha_k} - \frac{L_{\phi_k}}{2} -  \frac{\beta_k L_\theta}{\gamma^2}  - \left(1+ 4 \left(\frac{1-\eta_k\rho_{g_2}}{\eta_k\rho_{g_2}} \right)^2 \right)  L_\theta^2 \left( (L_f + L_g)^2+ 1/\gamma^2 \right)  \ge \frac{1}{4\alpha_k},
		\] 
		and
		\[
		\frac{1}{2\beta_k} - \frac{\rho_{g_2}}{2} - \frac{L_{\phi_k}}{2} - \left(1+ 4 \left(\frac{1-\eta_k\rho_{g_2}}{\eta_k\rho_{g_2}} \right)^2 \right)  L_\theta^2 \left( (L_f + L_g)^2+ 1/\gamma^2  \right)  \ge \frac{1}{4\beta_k}.
		\]
		Consequently, the conclusion follows from (\ref{lem9_eq6}).
	\end{proof}
	
	\subsection{Proof of Theorem~\ref{prop1}}
	\label{proofofthm}
	
	By leveraging the monotonically decreasing property of the merit function $V_k$, we can establish the non-asymptotic convergence for the sequence ${(x^k, y^k, \theta^k)}$ generated by the proposed MEHA.
	
	\begin{theorem}\label{prop1-a}
		Under Assumptions \ref{assump-UL} and \ref{assump-LL}, 
		suppose $\gamma \in (0, \frac{1}{2\rho_{f_2} + 2\rho_{g_2} })$, $c_k =\underline{c} (k+1)^p$ with $p \in [0,1/2)$, $\underline{c} >0$ and $\eta_k \in [ \underline{\eta}, (1/\gamma - \rho_{f_2})/(L_f+1/\gamma)^2] \cap  [ \underline{\eta}, 1/\rho_{g_2})$ with $ \underline{\eta} > 0$, 
		then there exists $c_{\alpha}, c_\beta > 0$ such that when $\alpha_k \in ( \underline{\alpha}, {c_\alpha})$ and $\beta_k \in ( \underline{\beta}, c_\beta)$ with $\underline{\alpha}, \underline{\beta} > 0$, the sequence of $(x^k, y^k, \theta^k)$ generated by Algorithm \ref{MEHA} satisfies
		\[
		\min_{0 \le k \le K} \left\| \theta^{k} - \theta_{\gamma}^*(x^{k},y^{k}) \right\| = O\left(\frac{1}{K^{1/2}} \right),
		\]
		and 
		\[
		\min_{0 \le k \le K} R_k(x^{k+1}, y^{k+1}) = O\left(\frac{1}{K^{(1-2p)/2}} \right).
		\]
		Furthermore, if $p \in (0,1/2)$ and there exists $M > 0$ such that $ \psi_{c_k}(x^k, y^k) \le M$ for any $k$, the sequence of $(x^k, y^k)$ satisfies
		\[
		\varphi(x^K, y^K) - v_\gamma(x^K, y^K) = O\left(\frac{1}{K^{p}} \right).
		\]
	\end{theorem}
	\begin{proof}
		First, Lemma \ref{lem9} ensures the existence of $c_\alpha, c_\beta > 0$ for which (\ref{lem9_eq}) is valid under the conditions $\alpha_k \le c_\alpha$, $\beta_k \le c_\beta$.
		Upon telescoping (\ref{lem9_eq}) over the range $k = 0, 1, \ldots, K-1$, we derive
		\begin{equation}\label{prop1_eq1}
			\begin{aligned}
				&\sum_{k = 0}^{K-1} 
				\Big( 
				\frac{1}{4\alpha_k} \| x^{k+1} - x^k \|^2 + \frac{1}{4\beta_k} \| y^{k+1} - y^k \|^2 
				\\
				&\quad\qquad+  \frac{1}{2}\left(\frac{\underline{\eta}\rho_{g_2}}{1-\underline{\eta}\rho_{g_2}} \right)^2 \left( (L_f + L_g)^2+ 1/\gamma^2  \right)  \left\| \theta^{k} - \theta_{\gamma}^*(x^{k},y^{k}) \right\|^2 
				\Big) \\
				\le \, & V_0 - V_K \le  V_0,
			\end{aligned}
		\end{equation}
		where the last inequality is valid because $V_K$ is nonnegative. Thus, we have
		\[
		\sum_{k = 0}^{\infty}  \left\| \theta^{k} - \theta_{\gamma}^*(x^{k},y^{k}) \right\|^2  < \infty,
		\]
		and then
		\[
		\min_{0 \le k \le K} \left\| \theta^{k} - \theta_{\gamma}^*(x^{k},y^{k}) \right\| = O\left(\frac{1}{K^{1/2}} \right).
		\]
		
		According to the update rule of variables $(x,y)$ as defined in (\ref{update_x}) and (\ref{update_y}), we have that 
		\begin{equation}
			\begin{aligned}
				0 &\in c_kd_x^k + \mathcal{N}_X(x^{k+1}) + \frac{c_k}{\alpha_k} \left(x^{k+1} - x^k\right), \\
				0 &\in c_k d_y^k + c_k \partial_{y} g(x^{k+1}, y^{k+1}) + \mathcal{N}_Y(y^{k+1}) + \frac{c_k}{\beta_k} \left(y^{k+1} - y^k\right).
			\end{aligned}
		\end{equation}
		From the definitions of $d_x^k$ and $d_y^k$ provided in (\ref{def_dx}) and (\ref{def_dy}), and given $\nabla_{x} g(x^{k+1}, y^{k+1})\times \partial_{y} g(x^{k+1}, y^{k+1})\subseteq \partial g(x^{k+1}, y^{k+1})$, a result stemming from the weakly convexity of $g$ and its continuously differentiable property with respect to $x$ as outlined in Assumption \ref{assump-LL} and corroborated by Theorem 5 of \cite{gao2023moreau}, we deduce
		\begin{align*}
			(e^k_x, e^k_y) \in & \nabla F(x^{k+1}, y^{k+1}) + c_k \left(\nabla f(x^{k+1}, y^{k+1}) + \partial g(x^{k+1}, y^{k+1}) - \nabla v_{\gamma} (x^{k+1}, y^{k+1})\right) \\
			&+ \mathcal{N}_{X \times Y}(x^{k+1}, y^{k+1}),
		\end{align*}
		with
		\begin{equation}
			\begin{aligned}
				e^k_x &:= \nabla_x \psi_{c_k}(x^{k+1}, y^{k+1}) - c_k d_x^k - \frac{c_k}{\alpha_k} \left(x^{k+1} - x^k\right), \\
				e^k_y & := \nabla_y \left(\psi_{c_k} -  c_k g \right)(x^{k+1}, y^{k+1}) - c_kd_y^k - \frac{c_k}{\beta_k} \left(y^{k+1} - y^k\right).
			\end{aligned}
		\end{equation}
		Next, we estimate $\|e^k_x \|$. We have
		\[
		\|e^k_x \| \le \| \nabla_x \psi_{c_k}(x^{k+1}, y^{k+1})  -  \nabla_x \psi_{c_k}(x^{k}, y^{k}) \| + \|  \nabla_x \psi_{c_k}(x^{k}, y^{k})  - c_k d_x^k\| + \frac{c_k}{\alpha_k} \left\|x^{k+1} - x^k\right\|.
		\]
		Considering the first term on the right hand side of the preceding inequality, and invoking Assumptions \ref{assump-UL} and \ref{assump-LL} alongside Lemma \ref{lem2-a}, \ref{lem3}, \ref{lem_theta}, we establish the existence of $L_{\psi_1} > 0$ such that
		\[
		\| \nabla_x \psi_{c_k}(x^{k+1}, y^{k+1})  -  \nabla_x \psi_{c_k}(x^{k}, y^{k}) \| \le c_k L_{\psi_1} \| (x^{k+1}, y^{k+1}) - (x^{k}, y^{k}) \|.
		\]
		Using (\ref{lem6_eq3}) and Lemma \ref{lem4}, we deduce
		\begin{equation}
			\begin{aligned}
				\|  \nabla_x \psi_{c_k}(x^{k}, y^{k})  - c_k d_x^k\|  = c_k\left\|  \nabla_x \phi_{c_k}(x^k, y^{k}) - d_x^k \right\|
				\le \, & c_k(L_f + L_g) \left\| \theta^{k} - \theta_{\gamma}^*(x^k,y^{k}) \right\|. 
			\end{aligned}
		\end{equation}
		Hence, we have
		\[
		\|e^k_x \| \le c_k L_{\psi_1} \| (x^{k+1}, y^{k+1}) - (x^{k}, y^{k}) \| + \frac{c_k}{\alpha_k} \left\|x^{k+1} - x^k\right\| + c_k(L_f + L_g) \left\| \theta^{k} - \theta_{\gamma}^*(x^k,y^{k}) \right\|.
		\]
		For $\|e^k_y\|$, it follows that
		\[
		\begin{aligned}
			\|e^k_y \| \le \, & \| \nabla_y \left(\psi_{c_k} -  c_k g \right)(x^{k+1}, y^{k+1}) -  \nabla_y \left(\psi_{c_k} -  c_k g \right)(x^{k+1}, y^{k}) \| + \frac{c_k}{\beta_k} \left\|y^{k+1} - y^k\right\|\\
			&+ \| \nabla_y \left(\psi_{c_k} -  c_k g \right)(x^{k+1}, y^{k})  - c_k d_y^k\|.
		\end{aligned}
		\]
		Analogously, invoking Assumptions \ref{assump-UL} and \ref{assump-LL} together with Lemmas \ref{lem2-a}, \ref{lem3}, and \ref{lem_theta},  we have the existence of $L_{\psi_2}:= L_{F}+L_{f}+\frac{1}{\gamma}+L_{\theta} > 0$ such that
		\[
		\| \nabla_y \left(\psi_{c_k} - c_k g \right)(x^{k+1}, y^{k+1}) -  \nabla_y \left(\psi_{c_k} -  c_k g \right)(x^{k+1}, y^{k}) \|  \le c_k L_{\psi_2} \|  y^{k+1} - y^{k} \|.
		\]
		Using (\ref{lem6_eq7}), Lemma \ref{lem_theta} and Lemma \ref{lem4}, we obtain
		\begin{align*}
			\| \nabla_y \left(\psi_{c_k} -  c_k g \right)(x^{k+1}, y^{k})  - c_k d_y^k\| = & c_k 	\left\|  \nabla_y \left(\phi_{c_k} -  g \right)(x^{k+1}, y^{k}) - d_y^k \right\| \\
			\le & \frac{c_k}{\gamma} \left(\left\| \theta^{k} - \theta_{\gamma}^*(x^{k},y^{k}) \right\| 
			+L_{\theta} \| x^{k+1} - x^{k} \|\right).
		\end{align*}
		Therefore, we have
		\[
		\|e^k_y \| \le c_k L_{\psi_2} \|  y^{k+1} - y^{k} \| + \frac{c_k}{\beta_k} \left\|y^{k+1} - y^k\right\| + \frac{c_k}{\gamma} \left(\left\| \theta^{k} - \theta_{\gamma}^*(x^{k},y^{k}) \right\| 
		+L_{\theta} \| x^{k+1} - x^{k} \|\right).
		\]
		With the estimations of $\|e^k_x \| $ and $\|e^k_y \| $, we obtain the existence of $L_\psi > 0$ such that
		\[
		\begin{aligned}
			R_k(x^{k+1}, y^{k+1}) \le \, &c_k L_{\psi} \| (x^{k+1}, y^{k+1}) - (x^{k}, y^{k}) \| + \left( \frac{c_k}{\alpha_k} 
			+\frac{c_k L_{\theta}}{\gamma}\right)\left\|x^{k+1} - x^k\right\| 
			\\ &+ \frac{c_k}{\beta_k} \left\|y^{k+1} - y^k\right\| 
			+ c_k(L_f + L_g + \frac{1}{\gamma}) \left\| \theta^{k} - \theta_{\gamma}^*(x^k,y^{k}) \right\|.
		\end{aligned}
		\]
		Employing the aforementioned inequality and given that 
		$\alpha_k \ge \underline{\alpha}$ and $\beta_k \ge \underline{\beta}$ for some positive constants $\underline{\alpha}, \underline{\beta}$, we demonstrate the existence of $C_R > 0$ such that 
		\begin{equation}
			\begin{aligned}
				&\frac{1}{c_k^2} 	R_k(x^{k+1}, y^{k+1})^2 \\
				\le \,& C_R
				\Big(
				\frac{1}{4\alpha_k} \| x^{k+1} - x^k \|^2 + \frac{1}{4\beta_k} \| y^{k+1} - y^k \|^2 
				\\
				&+  \frac{1}{2}\left(\frac{\underline{\eta}\rho_{g_2}}{1-\underline{\eta}\rho_{g_2}} \right)^2 \left( (L_f + L_g)^2+ 1/\gamma^2  \right)  \left\| \theta^{k} - \theta_{\gamma}^*(x^{k},y^{k}) \right\|^2 
				\Big).
			\end{aligned}
		\end{equation}
		Combining this with (\ref{prop1_eq1}) implies that
		\begin{equation}\label{prop1_eq2}
			\sum_{k = 0}^{\infty}	\frac{1}{c_k^2} R_k(x^{k+1}, y^{k+1})^2 < \infty.
		\end{equation}
		Because $2p < 1$, it holds that 
		\[
		\sum_{k = 0}^{K} \frac{1}{c_k^2}  = \frac{1}{\underline{c}^2} \sum_{k = 0}^{K} \left(\frac{1}{k+1} \right)^{2p} 
		\geq   \frac{1}{\underline{c}^2}\int_{1}^{K+2}\frac{1}{t^{^{2p}}}dt 
		\geq \frac{(K+2)^{1-2p}-1}{(1- 2p)\underline{c}^2 },
		\]
		and we can conclude from (\ref{prop1_eq2}) that
		\[
		\min_{0 \le k \le K} R_k(x^{k+1}, y^{k+1}) = O\left(\frac{1}{K^{(1-2p)/2}} \right).
		\]
		
		In fact, since 
		\begin{equation*}
			\sum_{k = \lfloor K/2 \rfloor-1}^{K} \frac{1}{c_k^2}  = \frac{1}{\underline{c}^2} \sum_{k = \lfloor K/2 \rfloor-1}^{K} \left(\frac{1}{k+1} \right)^{2p} 
			\geq   \frac{1}{\underline{c}^2}\int_{\lfloor K/2 \rfloor}^{2 \lfloor K/2 \rfloor}\frac{1}{t^{^{2p}}}dt 
			\geq \frac{2^{1-2p}-1}{(1- 2p)\underline{c}^2 } \lfloor K/2 \rfloor^{1-2p},
		\end{equation*}
		we have 
		\begin{equation}
			\min_{\lfloor K/2 \rfloor-1 \le k \le K} R_k(x^{k+1}, y^{k+1}) = O\left(\frac{1}{K^{(1-2p)/2}} \right).
		\end{equation}

		Finally, since $ \psi_{c_k}(x^k, y^k) \le M$ and $F(x^k, y^k) \ge \underline{F}$ for any $k$, we have
		\[
		c_k \Big( \varphi (x^k, y^k)- v_{\gamma}(x^k, y^k) \Big) \le M - \underline{F},  \quad \forall k,
		\]
		and we can obtain from $c_k =\underline{c} (k+1)^p$ that
		\[
		\varphi(x^K, y^K) - v_\gamma(x^K, y^K) = O\left(\frac{1}{K^{p}} \right).
		\]
	\end{proof}
	
	\subsection{Proof of Theorem \ref{prop2}}
	\label{relation}
	In the typical scenario where $\varphi(x,y) = f(x,y)$ being a strongly convex smooth function in $y$ for each $x$, and $X=\mathbb{R}^n$, $Y=\mathbb{R}^m$, it can be established that the residual function $R_k(x,y)$ is closely related to the norm of the hyper-gradient $\nabla \Phi(x)$ associated with the hyper-objective $\Phi(x):=F(x, y^*(x))$, where $y^*(x)$ is the unique LL optimal solution. 
	
	Recall that 
	this residual function is a stationarity measure for the following penalized version of the constrained problem (\ref{wVP}), with $c_k$ serving as the penalty parameter:
	\begin{equation*}
		\min_{(x,y)\in X \times Y} \psi_{c_k}(x,y) := F(x,y) + c_k \big( f(x,y)  - v_{\gamma} (x,y)\big).
	\end{equation*}	
	By the expression \eqref{valuefinc} of $\nabla v_{\gamma}(x,y)$ and the optimality condition of $\theta_{\gamma}^*:=\theta_{\gamma}^*(x,y)$ given as 
	\begin{equation*}
		\nabla_{y} f(x,\theta_{\gamma}^*)+\frac{\theta_{\gamma}^*-y}{\gamma}=0,
	\end{equation*}
	the residual function $R_k(x,y)$ in \eqref{residual} becomes 
	\begin{align*}
		R_k(x,y) &=  \left\| \begin{pmatrix}
			\nabla_{x}F(x,y) + c_k \left[ \nabla_{x}f(x,y) - \nabla_{x} f(x,\theta_\gamma^*(x,y)) \right] \\
			\nabla_{y}F(x,y) + c_k \left[\nabla_{y}f(x,y) - \nabla_{y} f(x,\theta_\gamma^*(x,y)) \right]
		\end{pmatrix} \right\|
		=: \left\| \begin{pmatrix}
			R_k^{(1)}(x,y) \\
			R_k^{(2)}(x,y)
		\end{pmatrix} \right\| .
	\end{align*}

	\begin{lemma}\label{hypergradient-approx}
		Under Assumptions \ref{assump-UL} and \ref{assump-LL}, suppose that $X=\mathbb{R}^n$, $Y=\mathbb{R}^m$, and the lower-level objective $\varphi(x,\cdot) = f(x, \cdot)$ is a $\mu$-strongly convex smooth function for each $x$. Let $\gamma > 1/\mu$, then 
		\begin{equation}\label{estimate1}
			\| y - y^*(x) \|  \leq 
			\frac{2L_{F,0} + 4\|R_k^{(2)}(x,y) \|}{c_k \mu},
		\end{equation}
		where $L_{F,0}$ is the upper bound of $\|\nabla_{y}F(x,y)\|$ on $X\times Y$. 
		Additionally, suppose $\|R_k^{(2)}(x,y) \|\leq L_{F,0}$,
		then $c_k\| y - y^*(x) \|  \leq 6L_{F,0}/\mu$. 
		If further $\nabla_{xy}^2 f(x,\cdot)$, $\nabla_{yy}^2 f(x,\cdot)$ are $L_{f,2}$-Lipschitz continuous on $X\times Y$, then 
		\begin{equation}\label{estimate0}
			\| 
			\nabla \Phi(x) - 
			R_k^{(1)}(x,y)
			\|
			\leq   \frac{L_{\mu}}{c_k}
			+\frac{L_{f}}{\mu} \| R_k^{(2)}(x,y)  \|,
		\end{equation}
		where $L_{\mu}:=\frac{6 L_{F,0}  }{\mu}
		\left( 1+\frac{L_f}{\mu} \right) 
		\left(
		L_F
		+\frac{6 L_{f,2} L_{F,0} }{\mu}
		\right)+\frac{6 L_{\Phi} L_{F,0}}{\mu^2 \gamma} $.
	\end{lemma}
	\begin{proof}
		First, by the optimality condition of $\theta_{\gamma}^*:=\theta_{\gamma}^*(x,y)$, we have  
		\begin{equation*}
			\left[\nabla_{yy}^2 f(x,\theta_{\gamma}^*) \right]
			\nabla_{y}\theta_{\gamma}^*
			+\frac{1}{\gamma}\left(\nabla_{y}\theta_{\gamma}^* - I\right)=0,
		\end{equation*}
		which implies that $\nabla_{y}\theta_{\gamma}^*=\frac{1}{\gamma}\left[\nabla_{yy}^2 f(x,\theta_{\gamma}^*)+\frac{1}{\gamma}I \right]^{-1}$. Note that $L_f I \succeq\nabla_{yy}^2 f\succeq \mu I$. Then $\frac{1}{\gamma} \left(\mu +\frac{1}{\gamma}\right)^{-1} I \succeq \nabla_{y}\theta_{\gamma}^* \succeq \frac{1}{\gamma} \left(L_f+\frac{1}{\gamma}\right)^{-1} I$. Thus 
		$\nabla_{yy}^2 v_{\gamma} = \frac{1}{\gamma}(I-\nabla_{y}\theta_{\gamma}^*)\preceq \frac{\mu}{\mu\gamma+1} I$. This implies that $f(x,\theta)-v_{\gamma}(x,\theta)$ is $\mu/2$-strongly convex in $\theta$ when $\gamma > 1/\mu$. 
		Thus, $\arg\min_{\theta} \big( f(x,\theta)  - v_{\gamma} (x,\theta)\big)$ has a unique solution, denoted by $y_\gamma^*(x)$. 
		We claim that $y_\gamma^*(x)=y^*(x)=:y^*$. Here, we denote $y^*(x)$ as $y^*$ for brevity. 
		Indeed, the expression \eqref{valuefinc} of $\nabla v_{\gamma}(x,y)$ shows that
		\begin{equation*}
			\nabla_{y} (f-v_{\gamma})(x, y^*)
			= \nabla_{y} f(x, y^*)- \frac{1}{\gamma}\big( y^* -\theta_\gamma^*(x, y^*) \big)=0,
		\end{equation*}
		by using the facts $\nabla_{y} f(x, y^*)=0$ and $\theta_\gamma^*(x, y^*) =y^*$.

		Second, we demonstrate the validity of estimate \eqref{estimate1}. 
		Since $f(x,\theta) - v_{\gamma}(x,\theta)$ is a $\mu/2$-strongly convex smooth function in $\theta$,  
		the function $\psi_{c_k}(x,\theta) =F(x,\theta)+ c_k   
		\big( f(x,\theta)  - v_{\gamma} (x,\theta)\big)$ is $(-L_F+c_k \mu/2)$-strongly convex in $\theta$. Let $c_k>4L_F/\mu$, then it is $c_k\mu/4$-strongly convex. Hence, $\arg\min_{\theta\in\mathbb{R}^m} \psi_{c_k}(x,\theta) $ has a unique solution, denoted by $y^*_{c_k}(x)$. By the first-order optimality condition, we have 
		\begin{equation*}
			\nabla_y F(x, y^*_{c_k}(x) )+ c_k  \big( \nabla_yf(x, y^*_{c_k}(x) ) - 
			\nabla_y v_\gamma(x, y^*_{c_k}(x) ) \big)
			=0.
		\end{equation*} 
		By the coercivity property of $\mu/2$-strongly convex functions $f(x,\cdot)-v_{\gamma}(x,\cdot)$, 
		$$
		\left\langle \nabla_{y} (f-v_{\gamma})(x, y^*(x)) - \nabla_{y} (f-v_{\gamma})(x, y^*_{c_k}(x)), 
		\theta_\gamma^*(x)- y^*_{c_k}(x) \right\rangle\geq \frac{\mu}{2} 
		\| \theta_\gamma^*(x)- y^*_{c_k}(x)\|^2.
		$$ 
		Recall that $\nabla_{y} (f-v_{\gamma})(x, y^*(x))=0$, we have 
		\begin{align}\label{estimate2}
			\| y^*(x)- y^*_{c_k}(x)\|
			\leq \frac{2}{\mu}  \| \nabla_y (f-v_{\gamma})(x, y^*_{c_k}(x) ) \|
			\leq \frac{2}{c_k \mu}  \| \nabla_y F(x, y^*_{c_k}(x) ) \|
			\leq \frac{2 L_{F,0}}{c_k \mu} .
		\end{align}
		Similarly, by the coercivity property of $c_k\mu/4$-strongly convex functions $F(x,\cdot)+c_k \big( f(x,\cdot)  - v_{\gamma} (x,\cdot)\big)$, 
		$$
		\left\langle \nabla_y F(x, y )+ c_k  \big( \nabla_y f(x,y)  - \nabla_y v_{\gamma} (x,y)\big),
		y- y^*_{c_k}(x) \right\rangle\geq \frac{c_k \mu }{4}
		\| y- y^*_{c_k}(x)\|^2.
		$$ 
		This implies that
		\begin{align}\label{estimate3}
			\| y- y^*_{c_k}(x)\|
			\leq \frac{4}{c_k \mu}  \| \nabla_y F(x, y )+ c_k  \big( \nabla_y f(x,y)  - \nabla_y v_{\gamma} (x,y)\big) \|
			=\frac{4\|R_k^{(2)}(x,y) \| }{c_k \mu}.
		\end{align}
		By combining the estimates \eqref{estimate2} and \eqref{estimate3}, we have successfully established estimate \eqref{estimate1}. 
		Hence, when $\|R_k^{(2)}(x,y) \|\leq L_{F,0}$, we get $c_k\| y - y^*(x) \|  \leq 6L_{F,0}/\mu$. 
		
		Next, we shall validate the estimate \eqref{estimate0}. Recall that 
		\begin{equation}
			\nabla\Phi(x)=\nabla_x F(x, y^*)
			-\nabla_{xy}^2 f(x,y^*)\left[\nabla_{yy}^2 f(x,y^*)\right]^{-1} \nabla_y F(x, y^*).
		\end{equation}
		Define 
		\begin{equation}
			\bar\nabla\Phi(x,y)=\nabla_x F(x, \theta_\gamma^*)
			-\nabla_{xy}^2 f(x,\theta_\gamma^*)\left[\nabla_{yy}^2 f(x,\theta_\gamma^*)\right]^{-1} \nabla_y F(x, \theta_\gamma^*).
		\end{equation}
		Then by Lemma 2.2 of \cite{ghadimi2018approximation}, there is a positive constant $L_{\Phi}$ such that 
		\begin{equation}
			\| \nabla\Phi(x) - \bar\nabla\Phi(x,y) \|
			\leq L_{\Phi} \| y^*(x)-\theta_\gamma^*(x,y)\|. 
		\end{equation}
		By $\nabla_{y} f(x,\theta_{\gamma}^*)+\frac{1}{\gamma}(\theta_{\gamma}^*-y)=0$ and the strongly convexity of $f(x,\cdot)$, we get 
		\begin{equation*}
			\left\langle \frac{1}{\gamma}(\theta_{\gamma}^*-y), y^*-\theta_\gamma^* \right\rangle \geq \mu \| y^* -\theta_\gamma^* \|^2,
		\end{equation*}
		which implies that 
		\begin{equation}
			\| y^* - \theta_\gamma^*\| \leq \frac{1}{\mu \gamma} \| \theta_\gamma^* - y\|.
		\end{equation}
		We claim that $\| \theta_\gamma^* - y\| \leq \|y-y^*\|$. In fact, by the optimality of $\theta_\gamma^*$, we have 
		\begin{equation*}
			f(x, \theta_\gamma^*)+\frac{1}{2\gamma} \| \theta_\gamma^* -y\|^2
			\leq f(x,y^*)+\frac{1}{2\gamma} \| y^*-y\|^2.
		\end{equation*}
		The desired result follows from the fact that $f(x,\theta_\gamma^*)-f(x,y^*)\geq \frac{\mu}{2} \| \theta_\gamma^* -y^*\|^2\geq0$. 
		Hence, 
		\begin{equation}
			\| \nabla\Phi(x) - \bar\nabla\Phi(x,y) \|
			\leq \frac{L_{\Phi}}{\mu \gamma} \|y-y^*(x)\|. 
		\end{equation}

		Since $R_k^{(1)}(x,y)=\nabla_{x}F(x,y) + c_k \left[ \nabla_{x}f(x,y) - \nabla_{x} f(x,\theta_\gamma^*(x)) \right]$, we have 
		\begin{align*}
			&\bar\nabla \Phi(x,y) - 
			R_k^{(1)}(x,y)\\
			=&\nabla_x F(x, \theta_\gamma^*)-\nabla_x F(x,y)\\
			&-\nabla_{xy}^2 f(x,\theta_\gamma^*)\left[\nabla_{yy}^2 f(x,\theta_\gamma^*)\right]^{-1} 
			\left(
			\nabla_y F(x, \theta_\gamma^*)
			-\nabla_y F(x, y)
			\right)\\
			&-\nabla_{xy}^2 f(x,\theta_\gamma^*)\left[\nabla_{yy}^2 f(x,\theta_\gamma^*)\right]^{-1} 
			\left(
			\nabla_{y}F(x,y) + c_k \big(\nabla_{y}f(x,y) - \nabla_{y}f(x,\theta_\gamma^*) \big) 
			\right)\\
			&+c_k\nabla_{xy}^2 f(x,\theta_\gamma^*)\left[\nabla_{yy}^2 f(x,\theta_\gamma^*)\right]^{-1} 
			\left(
			\nabla_y f(x,y)-\nabla_y f(x, \theta_\gamma^*)- \nabla_{yy}^2 f(x,\theta_\gamma^*)(y-\theta_\gamma^*)
			\right)\\
			&-c_k
			\left(
			\nabla_x f(x,y)-\nabla_x f(x, \theta_\gamma^*)- \nabla_{xy}^2 f(x,\theta_\gamma^*)(y-\theta_\gamma^*)
			\right).
		\end{align*}
		By Assumptions \ref{assump-UL} and \ref{assump-LL} (i), we have 
		\begin{align*}
			\| \nabla_x F(x, y^*)-\nabla_x F(x,y) \|
			\leq & L_F \| y^*-y\|, \\
			\| \nabla_y F(x, y^*)-\nabla_y F(x,y) \|
			\leq & L_F \| y^*-y\|.
		\end{align*}
		If further $\nabla_{xy}^2 f(x,\cdot)$, $\nabla_{yy}^2 f(x,\cdot)$ are $L_{f,2}$-Lipschitz continuous on $X\times Y$, then 
		\begin{align*}
			\| 
			\nabla_y f(x,y)-\nabla_y f(x, y^*)- \nabla_{yy}^2 f(x,y^*)(y-y^*)
			\|
			\leq &L_{f,2}\| y-y^* \|^2,
			\\
			\|
			\nabla_x f(x,y)-\nabla_x f(x, y^*)- \nabla_{xy}^2 f(x,y^*)(y-y^*)
			\|
			\leq & L_{f,2}\| y-y^* \|^2.
			\\
		\end{align*}
		Therefore, by the $\mu$-strongly convex of $f(x,\cdot)$ and $c_k\| y - y^*(x) \|  \leq 6L_{F,0}/\mu$, we have
		\begin{align*}
			\| \bar\nabla \Phi(x,y) - 
			R_k^{(1)}(x,y)\|
			\leq &L_F \| y^*-y\|
			+ \frac{L_f}{\mu} L_F \| y^*-y\|
			+\frac{L_{f}}{\mu} \| R_k^{(2)}(x,y)  \| \\
			& + c_k  \frac{L_f}{\mu}  L_{f,2}\| y-y^* \|^2
			+ c_k L_{f,2}\| y-y^* \|^2
			\\
			\leq &
			\left( 1+\frac{L_f}{\mu} \right) \| y^*-y\|
			\left(
			L_F
			+L_{f,2} c_k \| y^*-y\|
			\right)
			+\frac{L_{f}}{\mu} \| R_k^{(2)}(x,y) \|
			\\
			\leq &
			\left( 1+\frac{L_f}{\mu} \right) \| y^*-y\|
			\left(
			L_F
			+\frac{6 L_{f,2} L_{F,0} }{\mu}
			\right)
			+\frac{L_{f}}{\mu} \| R_k^{(2)}(x,y) \|\\
			\leq &
			\frac{6 L_{F,0}  }{\mu}
			\left( 1+\frac{L_f}{\mu} \right) 
			\left(
			L_F
			+\frac{6 L_{f,2} L_{F,0} }{\mu}
			\right)
			\frac{1}{c_k}
			+\frac{L_{f}}{\mu} \| R_k^{(2)}(x,y) \|.
		\end{align*}
		Therefore, 
		\begin{align*}
			\| \nabla \Phi(x) - 
			R_k^{(1)}(x,y)\|
			\leq &  \| \bar\nabla \Phi(x,y) - 
			R_k^{(1)}(x,y)\| + \| \nabla\Phi(x) - \bar\nabla\Phi(x,y) \|\\
			\leq & 
			\frac{6 L_{F,0}  }{\mu}
			\left( 1+\frac{L_f}{\mu} \right) 
			\left(
			L_F
			+\frac{6 L_{f,2} L_{F,0} }{\mu}
			\right)
			\frac{1}{c_k}
			+\frac{L_{f}}{\mu} \| R_k^{(2)}(x,y) \|
			+\frac{6 L_{\Phi} L_{F,0}}{\mu^2 \gamma} \frac{1}{c_k},
		\end{align*}
		which implies the desired outcome and thereby concludes the proof.
	\end{proof}
	
	We can obtain Theorem \ref{prop2} by combining Lemma \ref{hypergradient-approx} and Theorem \ref{prop1}.

	\subsection{Verifying Assumptions \ref{assump-LL}(ii) and (iii) are special cases of Assumption \ref{assump-LL}(iv)} \label{A9-assump}
	
	In this section, we prove that Assumptions \ref{assump-LL}(ii) and (iii) are special cases of Assumption \ref{assump-LL}(iv).
	
	\begin{lemma}
		The function $g(x,y)= x \| y \|_1$ satisfies Assumption \ref{assump-LL}(iv) when $X=\mathbb{R}_+$ and $Y = \mathbb{R}^m$.
	\end{lemma}
	\begin{proof}
		Initially, as depicted in Section 6.1 of \cite{gao2023moreau}, for $x\in\mathbb{R}_+$ and $y\in\mathbb{R}^m$, 
		\begin{equation}
			x \| y \|_1 + \frac{\sqrt{p}}{2} x^2
			+\frac{\sqrt{p}}{2} \| y \|^2
			=\sum_{i=1}^m \frac{1}{2\sqrt{p}} \left(
			x+\sqrt{p} | y_i |
			\right)^2,
		\end{equation}
		which is convex with respect to $(x,y)\in\mathbb{R}_+\times\mathbb{R}^m$. 
		Consequently, $g(x,y)$ is $\sqrt{p}$-weakly convex. 
		Further, for any given $s \in (0, \bar{s}]$, we have 
		\begin{equation}
			\mathrm{Prox}_{s \tilde{g}(x,\cdot)} (\theta)
			=\mathrm{Prox}_{s x \| \cdot \|_1} (\theta)
			=\mathcal{T}_{sx}(\theta) =\left(\mathcal{T}_{sx}(\theta_i)\right)_{i=1}^m
			= \left( [|\theta_i|-sx] _+ \cdot \mathrm{sgn}(\theta_i) \right)_{i=1}^m.
		\end{equation}
		This results in
		\begin{equation}\label{ass-LL-eq-a}
			\left\| 
			\mathrm{Prox}_{s \tilde{g}(x,\cdot)} (\theta) - \mathrm{Prox}_{s \tilde{g}(x',\cdot)}(\theta)
			\right\| 
			\leq s \|x - x' \|
			\leq \bar{s} \|x - x' \|.
		\end{equation}
		In summary, Assumption \ref{assump-LL}(iv) is satisfied by $g(x,y)= x \| y \|_1$ when $X=\mathbb{R}_+$ and $Y = \mathbb{R}^m$.
	\end{proof}
	
	\begin{lemma}
		The function $g(x,y) = \sum_{j = 1}^J x_j\| y^{(j)}\|_2$, where $\{1, \ldots, m\}$ is divided into $J$ groups, $y^{(j)}$ denotes the corresponding $j$-th group of $y$, satisfies Assumption \ref{assump-LL}(iv) when $X=\mathbb{R}^J_+$ and $Y = \mathbb{R}^m$.
	\end{lemma}
	\begin{proof}
		Initially, as depicted in Section 6.2 of \cite{gao2023moreau}, for $x\in\mathbb{R}^J_+$ and $y\in\mathbb{R}^m$, 
		\begin{equation}
			\sum_{j = 1}^J x_j\| y^{(j)}\|_2 + \sum_{j = 1}^J \frac{1}{2} x_j^2
			+\frac{1}{2} \| y \|_2^2
			=\sum_{j=1}^J \frac{1}{2} \left(
			x_j+ \| y^{(j)} \|_2
			\right)^2,
		\end{equation}
		which is convex with respect to $(x,y)\in\mathbb{R}^J_+\times\mathbb{R}^m$. 
		Consequently, $g(x,y)$ is $1$-weakly convex. 
		Further, for any given $s \in (0, \bar{s}]$, we have 
		\begin{equation}
			\begin{aligned}
				\mathrm{Prox}_{s \tilde{g}(x,\cdot)} (\theta)
				&=\mathrm{Prox}_{s\sum_{j = 1}^J x_j\| \cdot^{(j)}\|_2} (\theta) \\
				&= \bigtimes_{j=1}^J \mathrm{Prox}_{s x_j\| \cdot \|_2} (\theta^{(j)} )\\
				&= \bigtimes_{j=1}^J \left\{ \begin{array}{ll}
					[ \|\theta^{(j)}\|_2 - sx_j] _+  \frac{\theta^{(j)}}{\|\theta^{(j)}\|_2}, \quad &\theta^{(j)} \neq 0, \\
					0, & \theta^{(j)} = 0.
				\end{array} \right. 
			\end{aligned}
		\end{equation}
		This results in
		\begin{equation}\label{ass-LL-eq-a}
			\left\| 
			\mathrm{Prox}_{s \tilde{g}(x,\cdot)} (\theta) - \mathrm{Prox}_{s \tilde{g}(x',\cdot)}(\theta)
			\right\| 
			\leq  s \|x - x' \|
			\leq \bar{s} \|x - x' \|.
		\end{equation}
		In summary, Assumption \ref{assump-LL} (iv) is satisfied by $g(x,y)= \sum_{j = 1}^J x_j\| y^{(j)}\|_2$ when $X=\mathbb{R}^J_+$ and $Y = \mathbb{R}^m$.
	\end{proof}
		\begin{algorithm}[h]
		\caption{Single-loop Moreau Envelope based Hessian-
			free Algorithm (Smooth Case)}\label{MEHA-SC}
		\hspace*{0.02in} {\bf Initialize:} 
		$x^0,y^0,\theta^0$, learning rates $\alpha_k, \beta_k, \eta_k$, proximal parameter $\gamma$, penalty parameter $c_k$;
		\begin{algorithmic}[1]
			\FOR{$k=0,1,\dots,K-1$} 
			\STATE update 
			\[
			\begin{aligned}
				\theta^{k+1} &=  \mathrm{Proj}_Y \left( \theta^k - \eta_k \left(\nabla_y f(x^k,\theta^k)  + \frac{1}{\gamma}(\theta^k - y^k)\right) \right),\\
				x^{k+1} &=  \mathrm{Proj}_X \left( x^k - \alpha_k \left(
				\frac{1}{c_k}\nabla_{x}F(x^k, y^k)
				+\nabla_{x} f(x^k, y^k) - \nabla_{x} f(x^k, \theta^{k+1})
				\right) \right) ,\\
				y^{k+1} &= \mathrm{Proj}_Y \left( y^{k}  - \beta_k  
				\left(
				\frac{1}{c_k}\nabla_{y}F(x^{k+1}, y^k)
				+\nabla_{y}f(x^{k+1}, y^k)
				-\frac{1}{\gamma}(y^k - \theta^{k+1})
				\right) \right).
			\end{aligned}
			\]
			\ENDFOR
		\end{algorithmic}
	\end{algorithm}
	\subsection{Single-loop Moreau Envelope based Hessian-
		free Algorithm (Smooth Case)}\label{smoothcase}
	
	For smooth BLO problems, specifically when $g(x,y)\equiv0$ in (\ref{general_problem}), Algorithm \ref{MEHA} is specialized to Algorithm \ref{MEHA-SC}.

\end{document}